\RequirePackage{fix-cm}
\documentclass[12pt]{article}
\newcommand{\blind}{1}
\usepackage{placeins}
\usepackage{booktabs}    
\usepackage{multirow}    
\usepackage{graphicx}    
\usepackage{amssymb,amsthm}
\usepackage{amsmath}
\usepackage{natbib}
\usepackage{bibunits}
\usepackage{longtable}
\defaultbibliographystyle{plainnat}  
\setcitestyle{round} 
\defaultbibliography{bibfile_main} 

\usepackage{graphics,enumitem,subcaption,physics}
\usepackage{tikz}
\usetikzlibrary{calc}
\usepackage{tikz}
\usetikzlibrary{arrows.meta, positioning}
\usetikzlibrary{positioning} 
\usetikzlibrary{arrows.meta}
\usetikzlibrary{shadows}
\usepackage{algorithm}   
\usepackage{algpseudocode} 
\usepackage{amsmath}
\usepackage{enumitem}
\usepackage{booktabs}
\usepackage{multirow}
\usepackage{setspace}
\newtheorem{theorem}{Theorem}[section]
\newtheorem{remark}{Remark}[section]

\newtheorem{assumption}{Assumption}[section]

\newtheorem{lemma}{Lemma}[section]
\newenvironment{keywords}{\noindent\textbf{Keywords:} }{}

\usepackage{chngcntr}
\counterwithin{table}{section}
\counterwithin{algorithm}{section}
\counterwithin{figure}{section}
\numberwithin{equation}{section}
\usepackage{fontspec}

\usepackage[colorlinks=true, linkcolor=blue, urlcolor=blue, citecolor=blue]{hyperref}
\usepackage{xcolor}

\allowdisplaybreaks[4]
\usepackage{float}
\begin{document}	
\begin{bibunit}
	\if1\blind
	{
		\title{\vspace{-0.5cm}Learning the Average Dose-Response Function with a Right-Censored Treatment}
      \author{
Sihan Wu and Yifan Cui\\
Zhejiang University
}
\date{}

        \date{}
		\maketitle
	} \fi

	\if0\blind
	{
		\bigskip
		\bigskip
		\bigskip
		\begin{center}
			{\LARGE\bf Learning the Average Dose-response Function with a Right-Censored Treatment}
		\end{center}
		\medskip
	} \fi

\begin{abstract}
This paper considers estimation and inference for the average dose-response function under right-censored treatment duration and missing outcomes for censored units. Using the duration lower bounds provided by censored records, we construct a novel multiply robust estimator through augmented local estimating equations. Under suitable regularity conditions, the estimator remains consistent when any two of the outcome regression, treatment density, and censoring survival function are consistently estimated. We establish its convergence rate and asymptotic distribution within a cross-fitting framework under appropriate conditions on nuisance estimation. To address local smoothing bias that can distort confidence interval coverage, we further develop a debiased estimator that supports asymptotically valid pointwise inference without undersmoothing. Under additional uniform conditions, we construct simultaneous confidence bands using a Gaussian multiplier bootstrap. Simulation studies evaluate estimation accuracy, the validity of pointwise and uniform inference, and robustness to nuisance model misspecification. A semi-synthetic study calibrated to NHANES illustrates the proposed methods in a setting with a known dose-response function.
\end{abstract}


\section{Introduction}
\label{sec:introduction}

In many intervention studies, researchers seek to understand not only whether an intervention is effective, but also how its benefits vary with treatment duration. In inpatient rehabilitation, length of stay is often considered alongside functional recovery and the use of healthcare resources. If the benefits of rehabilitation vary with duration, simply dividing patients into short- and long-stay groups, or summarizing the relationship through a single linear effect, may obscure how benefits accumulate or where they begin to level off. These questions motivate estimation of the average dose-response function (ADRF), which describes the population mean potential outcome under interventions assigning different treatment durations \citep{GillRobins2001, kennedy2017, doss2024}.

However, the complete treatment duration and the corresponding outcome
may only become available when treatment is completed. Consider an
inpatient rehabilitation cohort with ongoing enrollment, in which
the outcome of interest is measured at discharge. At the data cutoff,
a patient may have received 10 days of rehabilitation but still be
undergoing treatment. Although the complete course might ultimately
last 20 or 40 days, the available record shows only that treatment
has lasted at least 10 days. Treatment duration is therefore
right-censored \citep{coxoakes1984}, and the outcome has not yet been measured.

Even when treatment and outcome are fully observed, estimating the ADRF for a continuous treatment is challenging. Observations generally
do not occur at a prespecified treatment level, so local smoothing methods borrow information from nearby treatment levels. A narrow neighborhood provides limited information, whereas a wider neighborhood may smooth over important variation in the response function and introduce smoothing bias \citep{FanGijbels1996,kennedy2017}. Right-censoring further complicates this local estimation because both the exact treatment duration and the outcome at completion
are unavailable for censored individuals. Without knowing the eventual duration, we cannot directly determine how much weight an individual should receive when estimating the response at a particular treatment level. Even if this weight were known, the missing outcome would still prevent direct calculation of the corresponding response contribution.

Existing methods for fully observed continuous treatments provide
important foundations. Propensity score methods
\citep{imbens2000role,ImaiVanDyk2004} and weighting approaches
\citep{ai2021unified,huling2024} adjust for observed confounding through
treatment assignment modeling or covariate balancing.
\citet{kennedy2017} develop doubly robust ADRF estimators using
pseudo-outcomes that combine outcome regression with treatment density
adjustment. Subsequent work investigates convergence rates and
higher-order corrections \citep{bonvinikennedy2026}. Related kernel
methods \citep{su2019,colangelolee2026} and series methods
\citep{SemenovaChernozhukov2021} accommodate high-dimensional covariates
or flexible nuisance estimation. Other developments include global
tests of a constant ADRF \citep{doss2024} and alternative identification
strategies based on instrumental variables or proxies for unmeasured
confounding \citep{chen2026iv,wu2024proximal,ham2026proximal}.

The fully observed-data ADRF estimators reviewed above cannot be applied directly to
censored records, and restricting them to complete records without
selection adjustment can introduce bias. Under our identification
assumptions, outcome regression fitted to complete records and
standardized over the full population covariate distribution can
identify the ADRF, but relies on consistent estimation of the outcome regression.
Semiparametric theory provides a general framework for coarsened data
\citep{rotnitzkyrobins1995,tsiatis2006,wangrotnitzkylin2010}.
For censored outcomes, \citet{FanGijbels1994} develop local linear regression based on transformed responses, and \citet{rubinvanderlaan2007} introduce
a doubly robust censoring-unbiased transformation. Furthermore, \citet{cui2023survival} develop causal survival forests
to estimate heterogeneous treatment effects, and \citet{westling2024inference} propose cross-fitted estimators of treatment-specific survival curves with machine learning nuisance estimation and establish pointwise and uniform inference. When treatment itself
is censored, however, both local weights and response contributions
are unavailable for censored individuals. We address this difficulty by augmenting both the
local response vector and design matrix using observed treatment
duration lower bounds. Under suitable regularity conditions, we show that the resulting estimator is multiply robust.

Consistency alone does not ensure valid inference: at a bandwidth of mean-squared-error-optimal order, smoothing bias can be of the same order as the standard error. Confidence intervals that ignore this bias may therefore fail to attain their nominal coverage. Undersmoothing reduces the bias at the cost of increased variance.
Shape-constrained methods provide another route to estimation and pointwise inference. For example, \citet{westling2020causal} develop causal isotonic regression under monotonicity, and \citet{doss2026monotonic} study doubly robust pointwise confidence intervals for a monotone dose-response function. Bias correction provides an alternative without imposing shape constraints \citep{calonico2018}. \citet{takatsuwestling2025} develop debiased pointwise inference and simultaneous confidence bands with bandwidths of order $n^{-1/5}$. Building on this approach, we incorporate bias correction into our
augmented local estimating equations. Under suitable smoothness and nuisance estimation conditions, the resulting estimator supports valid pointwise inference without undersmoothing and, under additional uniform conditions, simultaneous confidence bands.

The paper makes a number of contributions to both causal inference and survival analysis literature. First, we derive the observed data efficient influence function for an auxiliary integral functional of the dose-response function and use it to construct local multiply robust estimating equations. The resulting estimator is consistent whenever any two of the outcome regression, treatment density, and censoring survival function
are correctly specified. Second, we characterize the convergence rate and asymptotic distribution of the resulting estimator under cross-fitting. Third, we develop a debiased estimator that supports asymptotically valid pointwise inference with bandwidths of order $n^{-1/5}$, as well as simultaneous confidence bands under additional uniform conditions.

The remainder of the paper is organized as follows. Section~\ref{sec:setup} introduces the observed data structure, target parameter, and identification assumptions. Section~\ref{sec:mr} develops the multiply robust estimator and its asymptotic properties. Section~\ref{sec:debiased} presents the debiased estimator, pointwise inference, and simultaneous confidence bands. Section~\ref{sec:simulation} reports simulation studies, Section~\ref{sec:application} presents an NHANES-calibrated semi-synthetic study, and Section~\ref{sec:discussion} concludes with a discussion. Proofs and additional results are provided in the Appendix.

\section{Setup}
\label{sec:setup}

Let \(Z=(X,A,Y)\) denote the complete data, where \(X\in\mathcal X\subseteq\mathbb R^d\) is a vector of baseline covariates, \(A\in\mathcal A\subseteq\mathbb R\) is the treatment duration, and \(Y\in\mathbb R\) is the outcome measured at treatment completion. The conditional distribution of \(A\) given \(X\) has common support \(\mathcal A=[0,T]\), where \(0<T<\infty\). Let \(C\in[0,\infty)\) denote the time from treatment initiation until observation ends. We assume that \(A\) and \(C\) have conditional densities given \(X\). We observe \(n\) independent and identically distributed copies of
\[
O=(X,U,\Delta,\Delta Y),\qquad
U=\min(A,C),\qquad \Delta=I(A\le C),
\]
where \(I(\cdot)\) denotes the indicator function. When \(\Delta=1\), both \(A=U\) and \(Y\) are observed. When \(\Delta=0\), neither \(A\) nor \(Y\) is observed; we only know that \(A>U=C\).

Let \(Y(a)\) denote the potential outcome measured at completion under a well-defined intervention assigning treatment duration \(a\in\mathcal A\). Our target is the average dose-response function
$
\theta(a)=E\{Y(a)\}, a\in\mathcal A,
$
which describes the population mean potential outcome. 
Define the outcome regression \(\mu(x,a)=E(Y\mid X=x,A=a)\), the conditional treatment density \(f_A(a\mid x)\), and the conditional censoring survival function \(S_C(t\mid x)=P(C\ge t\mid X=x)\). 

We impose the following identification assumptions.

\begin{assumption}[Consistency]
\label{ass: consistency}
Potential outcomes are well defined, there is no interference between individuals, and \(Y=Y(A)\) almost surely.
\end{assumption}

\begin{assumption}[Unconfoundedness]
\label{ass: unconfoundedness}
Conditional on \(X\), treatment duration is independent of each potential outcome: \(Y(a)\perp\!\!\!\perp A\mid X\) for every \(a\in\mathcal A\).
\end{assumption}

\begin{assumption}[Conditionally independent censoring]
\label{ass: independent_censoring}
Censoring is independent of the complete treatment-outcome pair conditional on \(X\): \(C\perp\!\!\!\perp(A,Y)\mid X\).
\end{assumption}

\begin{assumption}[Positivity]
\label{ass: positivity}
There exist constants \(c_1,c_2>0\) such that \(f_A(a\mid X)\ge c_1\) and \(S_C(a\mid X)\ge c_2\) for all \(a\in\mathcal A\), almost surely.
\end{assumption}

Assumption~\ref{ass: consistency} links the realized outcome to the potential outcome under the treatment duration actually received. Assumption~\ref{ass: unconfoundedness} requires the baseline covariates to suffice for confounding adjustment: conditional on \(X\), treatment duration is independent of the potential outcome at each treatment level \citep{imbens2000role}. Assumption~\ref{ass: independent_censoring} allows the censoring time \(C\) to depend on \(X\), but requires \(C\) to be independent of \((A,Y)\) conditional on \(X\). In particular, \(P(\Delta=1\mid X,A,Y)=S_C(A\mid X)\), which provides the basis for inverse probability of censoring weighting and augmentation \citep{robinsrotnitzky1992,tsiatis2006}. The two parts of Assumption~\ref{ass: positivity} serve different roles. The treatment positivity condition ensures that, given \(X\), treatment durations near each level in \(\mathcal A\) have positive probability. The censoring positivity condition ensures that the probability of observing treatment completion at each level is bounded away from zero.

Write \(m(a)=E_X\{\mu(X,a)\}\), where \(E_X\) denotes expectation over the marginal covariate distribution in the full population. Under Assumptions~\ref{ass: consistency}–\ref{ass: positivity} and the regularity condition for continuous treatments stated in Appendix~\ref{app:pointwise_identification},
the target function is identified from the observed data distribution by
\[
\theta(a)
=m(a)
=E_X\!\left[E\{Y\mid X,U=a,\Delta=1\}\right].
\]
Conditionally independent censoring ensures that the outcome regression among individuals whose treatment is completed equals $\mu(X,a)$. Averaging this regression over the full population distribution of $X$ identifies the target function.

For the subsequent derivations, also define the marginal treatment density \(p_A(a)=E_X\{f_A(a\mid X)\}\) and the conditional treatment survival function \(S_A(t\mid x)=P(A\ge t\mid X=x)\). Write \(f_C(t\mid x)=-\partial_t S_C(t\mid x)\) and \(\lambda_C(t\mid x)=f_C(t\mid x)/S_C(t\mid x)\) for the conditional censoring density and hazard, respectively, with the hazard defined wherever \(S_C(t\mid x)>0\). The observed censoring counting process and at-risk process are \(N_C(t)=I(U\le t,\Delta=0)\) and \(R_C(t)=I(U\ge t)\), respectively. Under Assumption~\ref{ass: independent_censoring}, the corresponding censoring martingale has increments \(dM_C(t)=dN_C(t)-R_C(t)\lambda_C(t\mid X)\,dt\) \citep{andersen1993,fleming2013counting}.

Finally, define the auxiliary integrated functional
$
\Phi=\int_{\mathcal A}m(a)p_A(a)\,da
$
for the influence function construction in Section~\ref{sec:mr}. Its role is to guide the construction of estimating equations; the function \(\theta(a)\) remains the target of estimation and inference.

\section{Multiply Robust Estimation}
\label{sec:mr}

\subsection{Multiply Robust Estimator}

For a continuously distributed treatment, point evaluation of the dose-response function is generally not pathwise differentiable in a nonparametric model. We therefore use the integrated functional \(\Phi\), defined in Section~\ref{sec:setup}, as an auxiliary parameter for constructing a pseudo-outcome and local estimating equations. Semiparametric efficiency theory provides the efficient influence function for this regular functional \citep{bickel1993,newey1990}. \citet{kennedy2017} show that the full data efficient influence function of \(\Phi\) is
\begin{equation}
\label{eq:full_eif}
\varphi^F(Z)=\xi(X,A,Y)-\Phi+\int_{\mathcal A}\{\mu(X,t)-m(t)\}p_A(t)\,dt,
\end{equation}
where
\begin{equation}
\label{pseudo1}
\xi(X,A,Y)=\frac{Y-\mu(X,A)}{f_A(A\mid X)}p_A(A)+m(A).
\end{equation}
The pseudo-outcome is doubly robust: its conditional mean given \(A=a\) equals \(\theta(a)\) whenever either the outcome regression or the treatment density is correct.  This property motivates local linear regression of the pseudo-outcome on treatment. For later derivations, write \(R(X)=\int_{\mathcal A}\{\mu(X,t)-m(t)\}p_A(t)\,dt-\Phi\).

Let \(K_h(u)=h^{-1}K(u/h)\), where \(K\) is a continuous symmetric probability density supported on \([-1,1]\). Write \( z_{t,h}(a)=(1,(t-a)/h)^\top\), \(e_1=(1,0)^\top\), and let \(\mathbb P_n\) denote the sample average. With complete data, the local linear estimator of \citet{kennedy2017} can be written as
\[
\hat\theta_h^F(a)=e_1^\top\{\mathbb P_n[K_h(A-a) z_{A,h}(a)\{ z_{A,h}(a)\}^\top]\}^{-1}\mathbb P_n[K_h(A-a) z_{A,h}(a)\hat\xi],
\]
where \(\hat\xi\) is the estimated pseudo-outcome. For censored individuals, the exact treatment duration \(A\) and the outcome \(Y\) at treatment completion are both unobserved; the observed duration \(U\) supplies the lower bound \(A>U\). We therefore cannot directly evaluate the pseudo-outcome, the kernel weight \(K_h(A-a)\), or the local basis vector \( z_{A,h}(a)\). Define \(V(t,X)=E\{\varphi^F(Z)\mid A\ge t,X\}\). We adapt the full data influence function in \eqref{eq:full_eif} using censoring augmentation \citep{vanderlaanrobins2003,tsiatis2006}. 

\begin{theorem}[Observed Data Efficient Influence Function]
\label{thm:observed_eif}
Under Assumptions~\ref{ass: consistency}-\ref{ass: positivity}, consider the nonparametric model in which the complete data distribution and the conditional censoring distribution given $X$ are otherwise unrestricted within a dominated class. If $E[\{\varphi^F(Z)\}^2]<\infty$, the observed data efficient influence function of $\Phi$ is
\begin{equation}
\label{eq:observed_eif}
\varphi^O(O)=\frac{\Delta\varphi^F(X,U,Y)}{S_C(U\mid X)}+\int_0^T\frac{V(t,X)}{S_C(t\mid X)}\,dM_C(t).
\end{equation}
\end{theorem}

We apply the augmentation in \eqref{eq:observed_eif} to both parts of these local estimating equations. To isolate the robustness properties of the construction, we first hold the working functions \(\eta^*=(S_C^*,f_A^*,\mu^*)\) fixed. Here \(f_A^*\) is a conditional density on \(\mathcal A\), and \(S_C^*\) is an absolutely continuous survival function with \(S_C^*(0\mid x)=1\). Define \(m^*(t)=E_X\{\mu^*(X,t)\}\), \(p_A^*(t)=E_X\{f_A^*(t\mid X)\}\), \(S_A^*(t\mid x)=\int_t^\infty f_A^*(s\mid x)\,ds\), and \(\lambda_C^*(t\mid x)=-\partial_t\log S_C^*(t\mid x)\). All treatment densities are extended by zero outside \(\mathcal A\). Let \(\xi^*\) denote \eqref{pseudo1} evaluated at these working functions, and define
\begin{align*}
& H_N^*(t,X;a,h)=\int_t^\infty K_h(s-a) z_{s,h}(a) m^*(s)\frac{f_A^*(s\mid X)}{S_A^*(t\mid X)}\,ds,\\
& H_D^*(t,X;a,h)=\int_t^\infty K_h(s-a) z_{s,h}(a)\{ z_{s,h}(a)\}^\top\frac{f_A^*(s\mid X)}{S_A^*(t\mid X)}\,ds.
\end{align*}
The factor \(f_A^*(s\mid X)/S_A^*(t\mid X)\), for \(s\ge t\), is the working conditional density of the final treatment duration given that treatment lasts at least until \(t\). Thus \( H_N^*\) averages the local mean response contribution \(K_h(s-a) z_{s,h}(a) m^*(s)\), whereas \( H_D^*\) averages the local design contribution \(K_h(s-a) z_{s,h}(a)\{ z_{s,h}(a)\}^\top\), over these possible durations. These tail averages incorporate the known lower bound on treatment duration into both parts of the estimating equations. Write \(dM_C^*(t)=dN_C(t)-R_C(t)\lambda_C^*(t\mid X)\,dt\) for the working compensated process. The observed data contributions are
\begin{align*}
&\Gamma_N^*(O,a;h)=\frac{\Delta K_h(U-a) z_{U,h}(a)\xi^*(X,U,Y)}{S_C^*(U\mid X)}+\int_0^\infty\frac{ H_N^*(t,X;a,h)}{S_C^*(t\mid X)}\,dM_C^*(t),\\
&\Gamma_D^*(O,a;h)=\frac{\Delta K_h(U-a) z_{U,h}(a)\{ z_{U,h}(a)\}^\top}{S_C^*(U\mid X)}+\int_0^\infty\frac{ H_D^*(t,X;a,h)}{S_C^*(t\mid X)}\,dM_C^*(t).
\end{align*}
We set point contributions to zero when \(\Delta K_h(U-a)=0\),
and tail contributions to zero when \(t\ge a+h\) or
\(S_A^*(t\mid X)=0\). We also write these contributions as \(\Gamma_N(O,a;h,\eta^*)\) and \(\Gamma_D(O,a;h,\eta^*)\).

Define the sample numerator vector and design matrix by
\[
 N_h^*(a)=\mathbb P_n\{\Gamma_N(O,a;h,\eta^*)\},\qquad  D_h^*(a)=\mathbb P_n\{\Gamma_D(O,a;h,\eta^*)\}.
\]
The estimator based on the fixed working functions is
\begin{equation}
\label{eq:mr_estimator}
\beta_h^*(a)=( D_h^*(a))^{-1} N_h^*(a),\qquad \theta_h^*(a)=e_1^\top\beta_h^*(a).
\end{equation}
Set \(\beta_h^*(a)=0\) when \( D_h^*(a)\) is singular. 

For the following result, fix \(0<a<\tau<T\) and a neighborhood
\(I\subset(0,\tau)\) of \(a\). Assume that \(\theta\) and \(p_A\)
are continuous near \(a\), with \(p_A(a)>0\), and that, for some
constants \(0<c<M<\infty\), uniformly over \(x\in\mathcal X\),
\(f_A\) and \(f_A^*\) lie in \([c,M]\) on \([0,T]\),
\(S_C\) and \(S_C^*\) are bounded below by \(c\) on \([0,\tau]\),
and both \(|\mu^*(x,t)|\) and \(E(Y^2\mid X=x,A=t)\)
are bounded above by \(M\) for \(t\in I\).

\begin{theorem}[Triple Robustness]
\label{thm:triple_robustness}
Under Assumptions~\ref{ass: consistency}-\ref{ass: positivity} and the regularity conditions stated above, suppose that any two of \(S_C^*\), \(f_A^*\), and \(\mu^*\) equal their true counterparts. If \(h=h_n\to0\) and \(nh\to\infty\), the estimator in \eqref{eq:mr_estimator} satisfies
\[
\theta_h^*(a)\xrightarrow{p}\theta(a).
\]
\end{theorem}

The three pairs of correct working functions have distinct interpretations. If the censoring survival function is correct, censoring adjustment recovers the expectations of the full data local contributions, and the pseudo-outcome has the required conditional mean when either the outcome regression or the treatment density is correct. If both the outcome regression and the treatment density are correct, the working tail averages are correct, allowing augmentation to compensate for an incorrect censoring survival function. Proofs are provided in the Appendix.

\subsection{Convergence Rate and Asymptotic Normality}
\label{sec:mr_asymptotics}

We now estimate the nuisance functions using a fixed number \(K\) of cross-fitting folds of equal size. Each estimating equation is evaluated on observations excluded from the corresponding nuisance fits. This separation builds on sample splitting ideas \citep{schick1986} and is central to modern inference with estimated nuisance functions \citep{chernozhukov2018dml}. For fully observed continuous treatments, it has also been combined with kernel localization \citep{colangelolee2026}. Here, the analysis must additionally account for censoring augmentation in both the local response vector and the design matrix. Let \(\mathbb P_{n,k}\) denote the empirical average over fold \(k\), and let \(\hat\eta_{-k}=(\hat S_{C,-k},\hat f_{A,-k},\hat\mu_{-k})\) denote functions fitted using the corresponding training sample. Evaluation at \(\hat\eta_{-k}\) also uses the fitted marginal functions \(\hat m_{-k}\) and \(\hat p_{A,-k}\). Define
\begin{align*}
\hat{ N}_h(a)=\frac1K\sum_{k=1}^K\mathbb P_{n,k}\{\Gamma_N(O,a;h,\hat\eta_{-k})\},\quad \hat{ D}_h(a)=\frac1K\sum_{k=1}^K\mathbb P_{n,k}\{\Gamma_D(O,a;h,\hat\eta_{-k})\}.
\end{align*}
The cross-fitted estimator is \(\hat\beta_h(a)=\{\hat{ D}_h(a)\}^{-1}\hat{ N}_h(a)\), with \(\hat\theta_h(a)=e_1^\top\hat\beta_h(a)\). Set \(\hat\beta_h(a)=0\) on the event that \(\hat{ D}_h(a)\) is singular.

Fix \(0<a<\tau<T\), and let \(I\subset(0,\tau)\) be a fixed neighborhood of \(a\). Write \(\bar\eta=(\bar S_C,\bar f_A,\bar\mu)\) for the fixed limits of the nuisance estimators. Fold indices are suppressed below; all convergence requirements apply to every fold. In addition to Assumptions~\ref{ass: consistency}-\ref{ass: positivity}, impose the following conditions.
\begin{enumerate}[label=(\alph*),ref=(\alph*),leftmargin=*,series=mrregularity]
\item \label{cond:rate:nuisance} \textit{Nuisance limits and correctness.}
\begin{enumerate}[label=(\roman*),ref=\theenumi(\roman*),leftmargin=*]
\item \label{cond:rate:limits} The fitted nuisance functions converge uniformly to their limits:
\begin{align*}
&\sup_{t\in[0,T],\,x\in\mathcal X}|\hat S_C(t\mid x)-\bar S_C(t\mid x)|=o_p(1),\\
&\sup_{t\in[0,T],\,x\in\mathcal X}|\hat f_A(t\mid x)-\bar f_A(t\mid x)|=o_p(1),\\
&\sup_{t\in I,\,x\in\mathcal X}|\hat\mu(x,t)-\bar\mu(x,t)|=o_p(1).
\end{align*}
\item \label{cond:rate:correctness} Any two of \(\bar S_C\), \(\bar f_A\), and \(\bar\mu\) equal their true counterparts.
\end{enumerate}

\item \label{cond:rate:smoothing} \textit{Kernel and smoothness.}
\begin{enumerate}[label=(\roman*),ref=\theenumi(\roman*),leftmargin=*]
\item \label{cond:rate:kernel} The kernel \(K\) is a continuous symmetric probability density supported on \([-1,1]\).
\item \label{cond:rate:smoothness} The ADRF \(\theta\) is twice continuously differentiable, and \(p_A\) is continuously differentiable near \(a\), with \(p_A(a)>0\).
\end{enumerate}

\item \label{cond:rate:bounds} \textit{Probability models and bounds.} For constants \(0<c<M<\infty\), the following bounds hold uniformly over \(x\in\mathcal X\). Bounds involving fitted functions hold with probability tending to one, uniformly over the folds.
\begin{enumerate}[label=(\roman*),ref=\theenumi(\roman*),leftmargin=*]
\item \label{cond:rate:density} Each of \(f_A\), \(\bar f_A\), and \(\hat f_A\) is a conditional density on the common support \([0,T]\), integrates to one, and is bounded between \(c\) and \(M\). The corresponding treatment survival functions satisfy
\[
\hat S_A(t\mid x)=\int_t^T\hat f_A(s\mid x)\,ds,\qquad \bar S_A(t\mid x)=\int_t^T\bar f_A(s\mid x)\,ds.
\]
\item \label{cond:rate:censoring} The functions \(S_C\), \(\bar S_C\), and \(\hat S_C\) are absolutely continuous and nonincreasing on \([0,T]\), are at least \(c\) on this interval, and equal one at zero. Their hazards are defined from these survival functions; in particular, \(\hat\lambda_C=-\partial_t\log\hat S_C\) and \(\bar\lambda_C=-\partial_t\log\bar S_C\) almost everywhere.
\item \label{cond:rate:moments} For \(t\in I\), \(|\bar\mu(x,t)|\le M\), \(|\hat\mu(x,t)|\le M\), and \(E(Y^2\mid X=x,A=t)\le M\).
\end{enumerate}

\item \label{cond:rate:marginalization} \textit{Marginal function estimation.}
\begin{enumerate}[label=(\roman*),ref=\theenumi(\roman*),leftmargin=*]
\item \label{cond:rate:marginal-definition} The marginal functions are estimated using only the corresponding training sample. Their limits are \(\bar m(t)=E_X\{\bar\mu(X,t)\}\) and \(\bar p_A(t)=E_X\{\bar f_A(t\mid X)\}\).
\item \label{cond:rate:marginal-error} The marginalization errors satisfy
\[
\sup_{|t-a|\le h}\left\{|\hat m(t)-E_X\{\hat\mu(X,t)\}|+|\hat p_A(t)-E_X\{\hat f_A(t\mid X)\}|\right\}=O_p(n^{-1/2}).
\]
Here \(E_X\) uses the true marginal distribution of \(X\), holding the fitted functions fixed.
\end{enumerate}

\item \label{cond:rate:bandwidth} \textit{Bandwidth.} The deterministic bandwidth sequence satisfies \(h=h_n\to0\) and \(nh\to\infty\).

\end{enumerate}

Condition~\ref{cond:rate:marginal-error} controls the error from computing the marginal averages and is a separate requirement from cross-fitting. Let \(q_n(a)\), \(r_n(a)\), and \(s_n(a)\) be nonnegative deterministic sequences bounding estimation errors relative to the true nuisance functions:
\begin{align*}
&\sup_{t\in[0,\tau],\,x\in\mathcal X}|\hat S_C(t\mid x)-S_C(t\mid x)|=O_p\{q_n(a)\},\\
&\sup_{t\in[0,T],\,x\in\mathcal X}|\hat f_A(t\mid x)-f_A(t\mid x)|=O_p\{r_n(a)\},\\
&\sup_{|t-a|\le h,\,x\in\mathcal X}|\hat\mu(x,t)-\mu(x,t)|=O_p\{s_n(a)\}.
\end{align*}
These bounds need not tend to zero for a misspecified nuisance component.

\begin{theorem}[Convergence Rate]
\label{thm: rate1}
Under Assumptions~\ref{ass: consistency}-\ref{ass: positivity}, conditions~\ref{cond:rate:nuisance}-\ref{cond:rate:bandwidth}, and the cross-fitting setup above, the estimator at the fixed interior point \(a\) satisfies
\[
|\hat\theta_h(a)-\theta(a)|=O_p\!\left\{(nh)^{-1/2}+h^2+q_n(a)r_n(a)+q_n(a)s_n(a)+r_n(a)s_n(a)\right\}.
\]
\end{theorem}

Theorem~\ref{thm: rate1} separates the usual local linear smoothing error, \((nh)^{-1/2}+h^2\) \citep{FanGijbels1996}, from products of nuisance estimation errors. This extends the separation between smoothing error and nuisance estimation error in fully observed dose-response estimation \citep{bonvinikennedy2026}. When two nuisance limits are correct, their errors vanish, and each pairwise product contains at least one vanishing factor, yielding consistency under the stated conditions.

For asymptotic normality, define the limiting pseudo-outcome
\[
\bar\xi(X,t,Y)=\frac{Y-\bar\mu(X,t)}{\bar f_A(t\mid X)}\bar p_A(t)+\bar m(t)
\]
and the conditional second moment function
\[
v(t)=E\!\left[\frac{S_C(t\mid X)}{\bar S_C(t\mid X)^2}\{\bar\xi(X,t,Y)-\theta(t)\}^2\,\middle|\,A=t\right].
\]
Strengthen condition~\ref{cond:rate:bandwidth} by requiring
\begin{equation}
\label{eq:clt_bandwidth}
nh^5=O(1),
\end{equation}
and impose the following additional conditions.

\begin{enumerate}[resume*=mrregularity]

\item \label{cond:clt:products} \textit{Nuisance error products.}
\[
q_n(a)r_n(a)+q_n(a)s_n(a)+r_n(a)s_n(a)=o\{(nh)^{-1/2}\}.
\]

\item \label{cond:clt:regularity} \textit{Variance and moments.}
\begin{enumerate}[label=(\roman*),ref=\theenumi(\roman*),leftmargin=*]
\item \label{cond:clt:variance} The function \(v(t)\) admits a version continuous at \(a\), with \(v(a)>0\).
\item \label{cond:clt:moments} For some \(\delta>0\), \(\sup_{t\in I}E(|Y|^{2+\delta}\mid A=t)<\infty\).
\end{enumerate}

\end{enumerate}

\begin{theorem}[Asymptotic Normality]
\label{thm:asymptotic_normality}
Under the assumptions of Theorem~\ref{thm: rate1}, the bandwidth requirement~\eqref{eq:clt_bandwidth}, and conditions~\ref{cond:clt:products}-\ref{cond:clt:regularity},
\[
\sqrt{nh}\{\hat\theta_h(a)-\theta(a)-b_h(a)\}\xrightarrow{d}\mathcal N\{0,\sigma^2(a)\},
\]
where
\[
b_h(a)=\frac12h^2\theta''(a)\nu_2,\qquad \nu_2=\int u^2K(u)\,du,\qquad \sigma^2(a)=\frac{v(a)}{p_A(a)}\int K(u)^2\,du.
\]
\end{theorem}

The normal approximation in Theorem~\ref{thm:asymptotic_normality} is centered at \(\theta(a)+b_h(a)\). With \(h\asymp n^{-1/5}\), the leading bias can be of the same order as the standard error, as in complete data ADRF estimation \citep{kennedy2017}. The next section develops a debiased estimator to address this smoothing bias.

\section{Debiased Multiply Robust Method}
\label{sec:debiased}

\subsection{Debiased Multiply Robust Estimator}

The multiply robust estimator in the preceding section does not remove the bias introduced by local smoothing. Theorem~\ref{thm:asymptotic_normality} identifies its leading term as \(\frac12h^2\nu_2\theta''(a)\), which can be of the same order as the standard error when \(h\asymp n^{-1/5}\). Following the local polynomial bias correction approach \citep{calonico2018,takatsuwestling2025}, we estimate \(\theta''(a)\) using a local quadratic fit and subtract the estimated bias from \(\hat\theta_h(a)\). Because the correction is itself estimated, its sampling variability and its covariance with the original fit must also be incorporated into inference.

For \(\ell\in\{h,b\}\), we write
\[
z_{s,\ell}(a)=\left(1,\frac{s-a}{\ell}\right)^\top,\qquad
z_{s,\ell}^{(2)}(a)=
\left(1,\frac{s-a}{\ell},\frac{(s-a)^2}{\ell^2}\right)^\top.
\]
Here \(z_{s,h}(a)\) is the local linear basis used in the
preceding section. We use bandwidth \(b\) for a local quadratic
fit and first give the definitions using the true nuisance
functions. The corresponding tail functions are
\begin{align*}
H_N^{(2)}(t,X;a,b)
&=\int_t^\infty K_b(s-a)z_{s,b}^{(2)}(a)m(s)
  \frac{f_A(s\mid X)}{S_A(t\mid X)}\,ds,\\
H_D^{(2)}(t,X;a,b)
&=\int_t^\infty K_b(s-a)z_{s,b}^{(2)}(a)
  \{z_{s,b}^{(2)}(a)\}^\top
  \frac{f_A(s\mid X)}{S_A(t\mid X)}\,ds.
\end{align*}
These functions average the local quadratic mean response
and design contributions, respectively, over treatment
durations of at least \(t\), conditional on \(X\).

We then define \(\Gamma_N^{(2)}(O,a;b)\) and
\(\Gamma_D^{(2)}(O,a;b)\) using the same censoring augmentation
as in the preceding section, with \(K_h\) replaced by \(K_b\),
each occurrence of \(z_{s,h}(a)\) replaced by
\(z_{s,b}^{(2)}(a)\), and the corresponding tail functions
replaced by \(H_N^{(2)}\) and \(H_D^{(2)}\).
The resulting contributions form a three-dimensional response
vector and a \(3\times3\) design matrix.

For the local quadratic fit, we use the same cross-fitting
folds and nuisance estimates as for the local linear estimator.
These include \(\hat m_{-k}\), \(\hat p_{A,-k}\), and the
estimated censoring hazard.
The previous conventions for unobserved outcomes and zero
tail contributions apply with bandwidth \(b\). Using the same evaluation folds as for \(\hat\theta_h(a)\), define
\begin{align*}
\hat{ N}_b^{(2)}(a)=\frac1K\sum_{k=1}^K\mathbb P_{n,k}\{\Gamma_N^{(2)}(O,a;b,\hat\eta_{-k})\},\quad \hat{ D}_b^{(2)}(a)=\frac1K\sum_{k=1}^K\mathbb P_{n,k}\{\Gamma_D^{(2)}(O,a;b,\hat\eta_{-k})\}.
\end{align*}
In the local quadratic expansion, the coefficient of \(((s-a)/b)^2\) is \(\frac12b^2\theta''(a)\). With \(e_3=(0,0,1)^\top\), the second derivative estimator is therefore
\begin{equation}
\label{eq:theta_double_prime}
\hat\theta''(a)=2b^{-2}e_3^\top(\hat{ D}_b^{(2)}(a))^{-1}\hat{ N}_b^{(2)}(a).
\end{equation}

To estimate the contribution of the quadratic term to the local linear fit, define
\[
 z_{s,h}^{B}(a)=\frac{(s-a)^2}{h^2} z_{s,h}(a)
=\left(\frac{(s-a)^2}{h^2},\frac{(s-a)^3}{h^3}\right)^\top.
\]
Its censoring-adjusted contribution is constructed from
\[
 H_B(t,X;a,h)=\int_t^\infty K_h(s-a) z_{s,h}^{B}(a)\frac{f_A(s\mid X)}{S_A(t\mid X)}\,ds
\]
and
\[
\Gamma_B(O,a;h)=\frac{\Delta K_h(U-a) z_{U,h}^{B}(a)}{S_C(U\mid X)}
+\int_0^\infty\frac{ H_B(t,X;a,h)}{S_C(t\mid X)}\,dM_C(t).
\]
Its cross-fitted average is
\[
\hat{ B}_h(a)=\frac1K\sum_{k=1}^K\mathbb P_{n,k}\{\Gamma_B(O,a;h,\hat\eta_{-k})\}.
\]
The scalar \(e_1^\top\{\hat{ D}_h(a)\}^{-1}\hat{ B}_h(a)\) measures how the quadratic term contributes to the intercept of the fitted local linear regression. Using this estimated projection makes the correction reflect the local design in the sample; its probability limit is the kernel moment \(\nu_2\), as established below.

Combining this projection with the derivative estimate in \eqref{eq:theta_double_prime}, we define the debiased multiply robust estimator by
\begin{equation}
\label{eq:debiased_estimator}
\hat\theta_{\mathrm{DB}}(a)=\hat\theta_h(a)-\frac12h^2\hat\theta''(a)\left(e_1^\top\{\hat{ D}_h(a)\}^{-1}\hat{ B}_h(a)\right).
\end{equation}
As before, an estimator is defined to be zero when a required design matrix is singular. The bandwidth conditions and asymptotic properties of \eqref{eq:debiased_estimator} are given in the following subsection.

\subsection{Pointwise and Joint Asymptotic Properties}
\label{sec:db_asymptotics}

Fix an interior point \(a\) and retain the cross-fitting setup, notation, and conditions~\ref{cond:rate:nuisance}-\ref{cond:rate:bandwidth} of Section~\ref{sec:mr_asymptotics}. To accommodate the local quadratic fit, impose the following additional requirements. All requirements involving fitted functions apply to every fold.

\begin{enumerate}[label=(\alph*),ref=(\alph*),leftmargin=*,start=8]

\item \label{cond:db:bandwidth} \textit{Bandwidth for the local quadratic fit.}
\begin{enumerate}[label=(\roman*),ref=\theenumi(\roman*),leftmargin=*]
\item \label{cond:db:bandwidth-basic} The deterministic bandwidth sequence satisfies \(b=b_n\to0\).
\item \label{cond:db:bandwidth-ratio} The bandwidth ratio satisfies \(h/b\to\kappa\in[0,\infty)\).
\end{enumerate}

\item \label{cond:db:marginal-error} \textit{Marginalization over the enlarged window.}
Strengthen condition~\ref{cond:rate:marginal-error} to
\[
\sup_{|t-a|\le h\vee b}\left\{|\hat m(t)-E_X\{\hat\mu(X,t)\}|+|\hat p_A(t)-E_X\{\hat f_A(t\mid X)\}|\right\}=O_p(n^{-1/2}),
\]
where \(h\vee b=\max(h,b)\). The definition of \(E_X\) and the training sample requirements remain as in condition~\ref{cond:rate:marginalization}.

\end{enumerate}

The error bounds \(q_n(a)\) and \(r_n(a)\) retain their definitions from Section~\ref{sec:mr_asymptotics}. Enlarge the window in the outcome regression error bound to cover both local fits:
\[
\sup_{|t-a|\le h\vee b,\,x\in\mathcal X}|\hat\mu(x,t)-\mu(x,t)|=O_p\{s_n(a)\}.
\]
These bounds are relative to the true nuisance functions and need not tend to zero for a misspecified component.

\begin{theorem}[Convergence Rate of the Debiased Estimator]
\label{thm:debiased_convergence}
Under Assumptions~\ref{ass: consistency}-\ref{ass: positivity}, conditions~\ref{cond:rate:nuisance}-\ref{cond:rate:bandwidth}, conditions~\ref{cond:db:bandwidth}-\ref{cond:db:marginal-error}, and the cross-fitting setup and error bounds above,
\[
|\hat\theta_{\mathrm{DB}}(a)-\theta(a)|=O_p\!\left\{\frac1{\sqrt{nh}}+\frac{h^2}{\sqrt{nb^5}}+q_n(a)r_n(a)+q_n(a)s_n(a)+r_n(a)s_n(a)\right\}+o(h^2).
\]
\end{theorem}

The proof of Theorem~\ref{thm:debiased_convergence} also establishes that \(e_1^\top\{\hat{ D}_h(a)\}^{-1}\hat{ B}_h(a)\xrightarrow{p}\nu_2\). The term \(h^2/\sqrt{nb^5}\) describes the sampling error introduced by estimating the second derivative. Since
\[
\frac{h^2}{\sqrt{nb^5}}=\left(\frac hb\right)^{5/2}\frac1{\sqrt{nh}},
\]
this term is of order \((nh)^{-1/2}\) or smaller under condition~\ref{cond:db:bandwidth-ratio}. The result does not require \(\hat\theta''(a)\) itself to be consistent.

\begin{remark}
If \(\theta\) is additionally four times continuously differentiable near \(a\), the deterministic remainder is \(O(h^4+h^2b^2)=O(h^2b^2)\). The rate in Theorem~\ref{thm:debiased_convergence} then becomes
\[
|\hat\theta_{\mathrm{DB}}(a)-\theta(a)|=O_p\!\left\{\frac1{\sqrt{nh}}+\frac{h^2}{\sqrt{nb^5}}+h^2b^2+q_n(a)r_n(a)+q_n(a)s_n(a)+r_n(a)s_n(a)\right\}.
\]
This stronger smoothness condition is used only for the sharper rate.
\end{remark}
For asymptotic normality, we also impose the bandwidth requirement~\eqref{eq:clt_bandwidth} and conditions~\ref{cond:clt:products}-\ref{cond:clt:regularity}, using the enlarged error bound \(s_n(a)\) above. Retain \(\bar\xi\) and \(v\) from Section~\ref{sec:mr_asymptotics}, and write \(\nu_j=\int u^jK(u)\,du\) for \(j=2,4\). The combined fluctuation of the original fit and its bias correction is described by the equivalent kernel
\begin{equation}
\label{eq:equivalent_kernel}
K_{\mathrm{DB}}(u)=K(u)-\frac{\nu_2\kappa^3}{\nu_4-\nu_2^2}(\kappa^2u^2-\nu_2)K(\kappa u).
\end{equation}
The kernel conditions imply \(\nu_4>\nu_2^2\). When \(\kappa=0\), define \(K_{\mathrm{DB}}=K\). This representation includes both the variance of the bias estimate and its covariance with the original estimator, as required for inference after bias correction \citep{calonico2018}.

\begin{theorem}[Asymptotic Normality of the Debiased Estimator]
\label{thm:debiased_asymptotic_normality}
Under the assumptions of Theorem~\ref{thm:debiased_convergence}, the bandwidth requirement~\eqref{eq:clt_bandwidth}, and conditions~\ref{cond:clt:products}-\ref{cond:clt:regularity},
\[
\sqrt{nh}\{\hat\theta_{\mathrm{DB}}(a)-\theta(a)\}\xrightarrow{d}\mathcal N\{0,\sigma_{\mathrm{DB}}^2(a)\},
\]
where
\begin{equation}
\label{eq:debiased_variance}
\sigma_{\mathrm{DB}}^2(a)=\frac{v(a)}{p_A(a)}V_{K,\kappa},\qquad V_{K,\kappa}=\int\{K_{\mathrm{DB}}(u)\}^2\,du.
\end{equation}
\end{theorem}

The limiting distribution is centered at \(\theta(a)\). Under twice continuous differentiability, the remaining bias is \(o(h^2)\), and \(nh^5=O(1)\) makes it negligible relative to the standard error. Thus Theorem~\ref{thm:debiased_asymptotic_normality} permits \(h\asymp n^{-1/5}\) without undersmoothing, as in debiased inference for fully observed continuous treatments \citep{takatsuwestling2025}. When \(\kappa=0\), the bias correction makes no first-order contribution to the variance.

Let \(\hat\sigma_{\mathrm{DB}}^2(a)\) be the variance estimator constructed from the estimated influence-function contributions for the local estimator, as detailed in Appendix~\ref{app:db_variance}. The resulting asymptotic \((1-\alpha)\) pointwise confidence interval is
\[
\hat\theta_{\mathrm{DB}}(a)\ \pm\ z_{1-\alpha/2}\frac{\hat\sigma_{\mathrm{DB}}(a)}{\sqrt{nh}},
\]
where \(z_{1-\alpha/2}\) is the corresponding standard normal quantile. This interval provides coverage at a fixed treatment level; simultaneous coverage over a treatment interval requires the construction in Section~\ref{sec:uniform_inference}.

We finally consider estimation at a fixed collection of treatment levels.

\begin{theorem}[Joint Asymptotic Normality]
\label{thm:joint_normality}
Let \(\mathcal A^*=\{a_1,\ldots,a_m\}\) be a finite, fixed collection of distinct interior points. Construct all estimators from the same sample and cross-fitting scheme, using common bandwidths \(h\) and \(b\). If the assumptions of Theorem~\ref{thm:debiased_asymptotic_normality} hold at every \(a_j\in\mathcal A^*\), then
\[
\sqrt{nh}\bigl(\hat\theta_{\mathrm{DB}}(a_1)-\theta(a_1),\ldots,\hat\theta_{\mathrm{DB}}(a_m)-\theta(a_m)\bigr)^\top\xrightarrow{d}\mathcal N_m(0,\Sigma_{\mathrm{DB}}),
\]
where \(\Sigma_{\mathrm{DB}}=\operatorname{diag}\{\sigma_{\mathrm{DB}}^2(a_1),\ldots,\sigma_{\mathrm{DB}}^2(a_m)\}\), with variances given by \eqref{eq:debiased_variance}.
\end{theorem}

The diagonal limit follows because the shrinking kernel windows around fixed distinct points eventually become disjoint. This argument applies to a fixed finite collection and does not make estimates on an increasingly dense grid independent: nearby kernel windows can overlap. Inference over an entire treatment interval therefore requires control of the estimation process across treatment levels, which we develop next.

\subsection{Uniform Inference and Confidence Bands}
\label{sec:uniform_inference}

In this section, we use Gaussian
multiplier bootstrap to construct simultaneous confidence bands \citep{Chernozhukov2014,Chernozhukov2015}.

Let \(\mathcal{A}_0\) be a fixed nondegenerate compact interval
contained in the interior of a compact interval
\(\mathcal{A}_C\subset(0,\tau)\), where \(\tau<T\). Using the estimated influence-function contributions
\(\hat\varphi_{\mathrm{DB},i}(a)\) from the pointwise analysis
(see Appendix~\ref{app:db_variance}), we write
\(\bar\varphi_{\mathrm{DB}}(a)
=n^{-1}\sum_i\hat\varphi_{\mathrm{DB},i}(a)\). Let \(\mathcal{A}_n\subset\mathcal{A}_0\) be a deterministic
nonempty finite grid. We generate random weights
\(G_1,\ldots,G_n\overset{\mathrm{iid}}{\sim}N(0,1)\),
independently of the observed data and define
\begin{equation}
\label{eq:uniform_multiplier}
Z_n^G(a)=\frac{\sqrt h}{\sqrt n\,\hat\sigma_{\mathrm{DB}}(a)}\sum_{i=1}^n G_i\{\hat\varphi_{\mathrm{DB},i}(a)-\bar\varphi_{\mathrm{DB}}(a)\},\qquad a\in \mathcal{A}_n.
\end{equation}
For fixed \(\alpha\in(0,1)\), let \(c_{1-\alpha,n}\) be the conditional \((1-\alpha)\)-quantile of \(\max_{a\in \mathcal{A}_n}|Z_n^G(a)|\), given the observed sample. The resulting band is
\begin{equation}
\label{eq:uniform_band}
\mathcal C_n(a)=\left[\hat\theta_{\mathrm{DB}}(a)-c_{1-\alpha,n}\frac{\hat\sigma_{\mathrm{DB}}(a)}{\sqrt{nh}},\quad
\hat\theta_{\mathrm{DB}}(a)+c_{1-\alpha,n}\frac{\hat\sigma_{\mathrm{DB}}(a)}{\sqrt{nh}}\right],\qquad a\in \mathcal{A}_0.
\end{equation}

The bounds \(q_n,r_n,s_n\) are defined as in
Section~\ref{sec:mr_asymptotics}, with the outcome regression
error now taken uniformly over
\(\mathcal{A}_C\times\mathcal X\).
Let \(\delta_n\) denote the sum of the uniform errors relative
to the fixed limits \(\bar S_C,\bar f_A,\bar\mu\) on the same
domains, maximized over the cross-fitting folds.

For uniform validity, we impose the following additional
requirements.
\begin{enumerate}[label=(\alph*),ref=(\alph*),leftmargin=*,start=10]
\item \label{cond:uniform:regularity} \textit{Uniform smoothness and moments.}
\begin{enumerate}[label=(\roman*),ref=\theenumi(\roman*),leftmargin=*]
\item \label{cond:uniform:smoothness} The function \(\theta''\) is H\"older continuous on \(\mathcal{A}_C\) with exponent \(\gamma\in(0,1]\). The density \(p_A\) is continuously differentiable on \(\mathcal{A}_C\), and \(v\) is continuous there. Both \(p_A\) and \(v\) are bounded above and away from zero on \(\mathcal{A}_C\).
\item \label{cond:uniform:moments} The kernel \(K\), extended by zero outside \([-1,1]\), is Lipschitz continuous. For some finite \(q\ge4\), \(\sup_{t\in \mathcal{A}_C}E(|Y|^q\mid A=t)<\infty\).
\end{enumerate}

\item \label{cond:uniform:nuisance} \textit{Uniform nuisance estimation.}
\begin{enumerate}[label=(\roman*),ref=\theenumi(\roman*),leftmargin=*]
\item \label{cond:uniform:limits} The error relative to the working limits satisfies \(\delta_n\log n=o_p(1)\).
\item \label{cond:uniform:marginalization}
For every fold, the marginalization error in
condition~\ref{cond:db:marginal-error} is \(O_p(n^{-1/2})\)
with the supremum taken over \(t\in\mathcal{A}_C\).
\item \label{cond:uniform:products} The product errors satisfy
$
\sqrt{nh\log n}\,(q_nr_n+q_ns_n+r_ns_n)\longrightarrow0.
$
\end{enumerate}

\item \label{cond:uniform:bandwidth-grid} \textit{Bandwidths and evaluation grid.}
\begin{enumerate}[label=(\roman*),ref=\theenumi(\roman*),leftmargin=*]
\item \label{cond:uniform:bandwidth}
The bandwidths satisfy
\[
\frac{nh}{(\log n)^7}\longrightarrow\infty,\qquad
\frac{n^{1-2/q}h}{(\log n)^3}\longrightarrow\infty.
\]
\item \label{cond:uniform:bias}
The bandwidths satisfy
$
\sqrt{nh^5\log n}\,(h^\gamma+b^\gamma)\to0.
$
\item \label{cond:uniform:grid} The deterministic grid \(\mathcal{A}_n\subset \mathcal{A}_0\) has size \(N_n=O(n^{d_0})\) for some fixed \(d_0>0\), and mesh \(\omega_n=\sup_{a\in \mathcal{A}_0}\min_{t\in \mathcal{A}_n}|a-t|\) satisfying \(\omega_n\log n/h\to0\).
\end{enumerate}

\end{enumerate}

\begin{theorem}[Simultaneous Confidence Bands]
\label{thm:uniform_band}
Suppose the assumptions of
Theorem~\ref{thm:debiased_asymptotic_normality}
hold with \(\mathcal{A}_C\) as the common neighborhood
and all bounds uniform over \(a\in\mathcal{A}_0\).
Under conditions~\ref{cond:uniform:regularity}-
\ref{cond:uniform:bandwidth-grid},
the band in \eqref{eq:uniform_band} satisfies,
for every fixed \(\alpha\in(0,1)\),
\[
P\{\theta(a)\in\mathcal C_n(a)\text{ for all }a\in \mathcal{A}_0\}\longrightarrow1-\alpha.
\]
\end{theorem}

Theorem~\ref{thm:uniform_band} establishes simultaneous coverage of the dose-response function over the entire interval \(\mathcal{A}_0\) with asymptotic probability \(1-\alpha\). Although the bootstrap critical value is computed on the finite
grid \(\mathcal{A}_n\), condition~\ref{cond:uniform:grid} ensures that the discretization error is asymptotically negligible, allowing the coverage guarantee to extend to the full interval. The proof is given in Appendix~\ref{app:uniform_inference}.

\section{Simulation Study}
\label{sec:simulation}

\subsection{Experimental Design}
\label{subsec:simulation_design}

We examine how accounting for censoring affects estimation, how
smoothing bias affects inference, and how bias correction changes
coverage and precision. Linear and nonlinear dose-response functions
allow these features to be assessed under different degrees of
curvature. We generate $X\sim\operatorname{Uniform}(-1,1)$
and let $\lambda_A(x)=\exp(-1+0.2x)$ and
$\lambda_C(x)=\exp(-0.9+1.2x)$. Conditional on $X=x$, treatment
follows an exponential distribution truncated to $[0,3]$, with density
\[
f_A(a\mid x)
=\frac{\lambda_A(x)\exp\{-\lambda_A(x)a\}}
       {1-\exp\{-3\lambda_A(x)\}}
  I(0\leq a\leq3).
\]
Independently of $A$ given $X$, we generate
$C\mid X=x\sim\operatorname{Exponential}\{\lambda_C(x)\}$,
with rate $\lambda_C(x)$. Outcomes follow
$Y=\theta_s(A)+(4+A)X+\varepsilon$, where
$\varepsilon\sim N(0,1)$ is independent of $(X,A,C)$ and
\[
\begin{aligned}
\text{Setting 1:}\quad &\theta_1(a)=2+1.5a,\\
\text{Setting 2:}\quad &\theta_2(a)=2+1.5a+6(a-1.2)^2.
\end{aligned}
\]
These are the true ADRFs since $E(X)=0$. We observe
$U=\min(A,C)$ and $\Delta=I(A\leq C)$, with $Y$
available only when $\Delta=1$. The censoring proportion is
approximately $37.5\%$ in both settings.

The main experiment compares the proposed multiply robust estimator
(MR) and its debiased version (MRDB) with Conventional and
ConventionalDB, which apply the corresponding local linear and
debiased procedures to complete cases. This comparison assesses
the effect of discarding censored records. We also include a regression
plug-in estimator (PI), obtained by averaging cross-fitted outcome
predictions over the covariate distribution of the full sample. PI
reuses the outcome regression fitted for MR and requires no
second-stage local smoothing. Nuisance functions are estimated using
kernel ridge regression for the outcome and neural network models
for the treatment density and censoring hazard. All methods use
two-fold cross-fitting; implementation details are provided in
Appendix~\ref{sec:sim-supplement}.

Local estimation uses the Epanechnikov kernel and two bandwidth rules,
\[
h_{\mathrm{r}}=2n^{-1/5},
\qquad h_{\mathrm{u}}=n^{-1/4},
\qquad b=1.25h,
\]
Here, $h_{\mathrm{r}}$ and $h_{\mathrm{u}}$ denote the regular
and undersmoothing bandwidths, respectively, and $b$ is the
pilot bandwidth for bias correction. Within each replication, both bandwidth rules use the same simulated data, folds, and nuisance estimates.

We use $n\in\{2000,4000,8000\}$, with $1000$ replications per
configuration and $201$ equally spaced evaluation doses in
$[0.8,1.8]$. Performance measures are pointwise bias, root mean
squared error (RMSE), coverage probability (CP), and the mean
length of nominal $95\%$ confidence intervals. PI is evaluated
only by bias and RMSE. For MRDB, we additionally assess $95\%$
simultaneous confidence bands using coverage over all evaluation
doses and mean band width. A supplementary experiment examines
robustness when one nuisance model is misspecified and the
other two are correctly specified.

\subsection{Numerical Results}
\label{subsec:simulation_results}

The complete case procedures exhibit substantial negative bias in both settings, and their coverage remains poor at all sample sizes. This bias arises because selection into the complete case sample depends on treatment and covariates, making it unrepresentative of the target population.

\begin{figure}[ht]
\centering
\linespread{1}\selectfont
\includegraphics[
  width=\textwidth,
  height=0.52\textheight,
  keepaspectratio
]{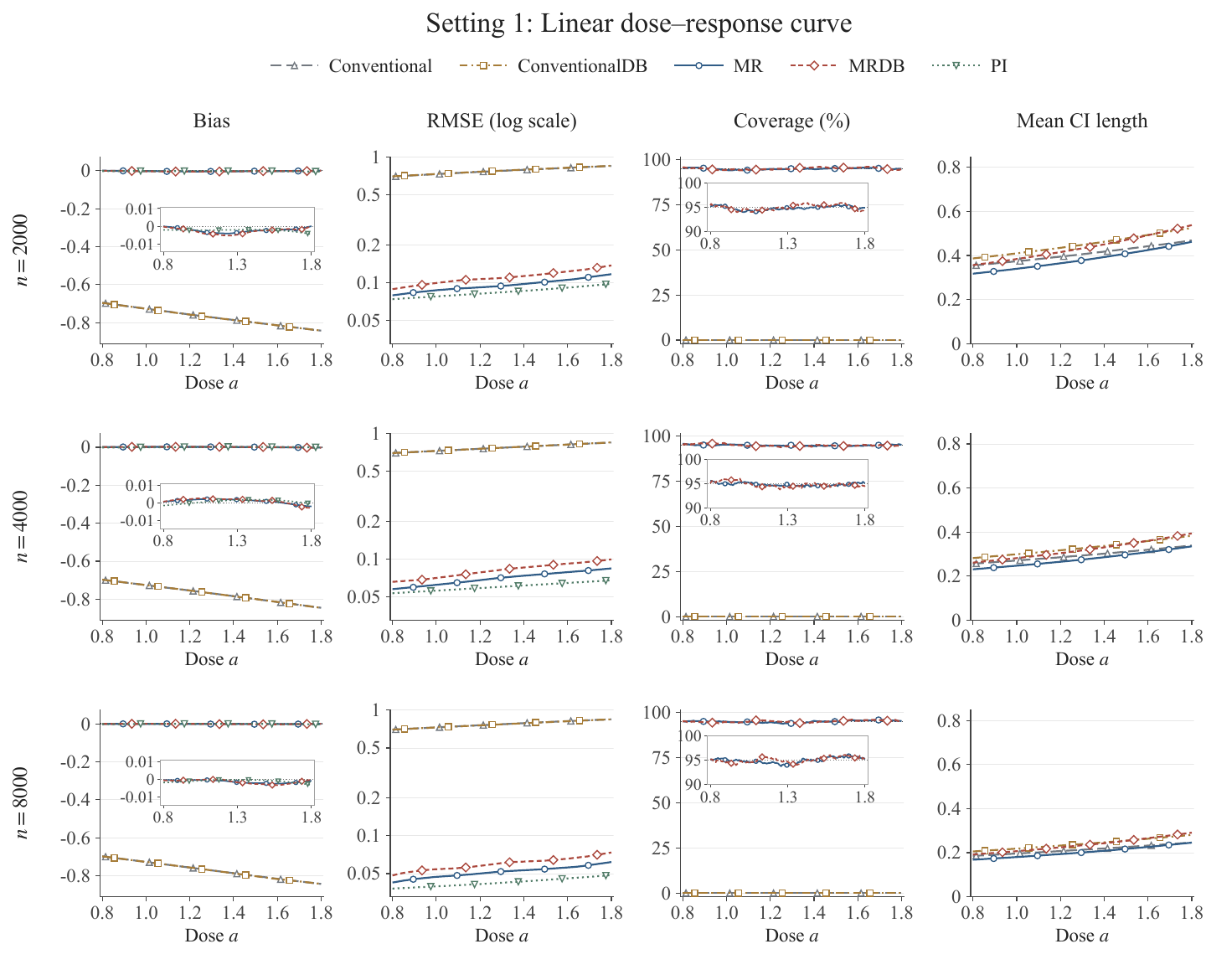}
\caption{Setting 1, regular bandwidth.
Insets enlarge near-zero bias and coverage between $90\%$ and
$100\%$. PI is shown only for bias and RMSE.}
\label{fig:sim-main-linear-regular}
\end{figure}

Figure~\ref{fig:sim-main-linear-regular} presents the
regular-bandwidth results for Setting 1. MR and MRDB both have small
bias and pointwise coverage close to $95\%$, with RMSE and interval
length decreasing as the sample size increases. Since the linear
curve introduces no curvature-induced smoothing bias, bias
correction provides little benefit here.

\begin{figure}[ht]
\centering
\linespread{1}\selectfont
\includegraphics[
  width=\textwidth,
  height=0.52\textheight,
  keepaspectratio
]{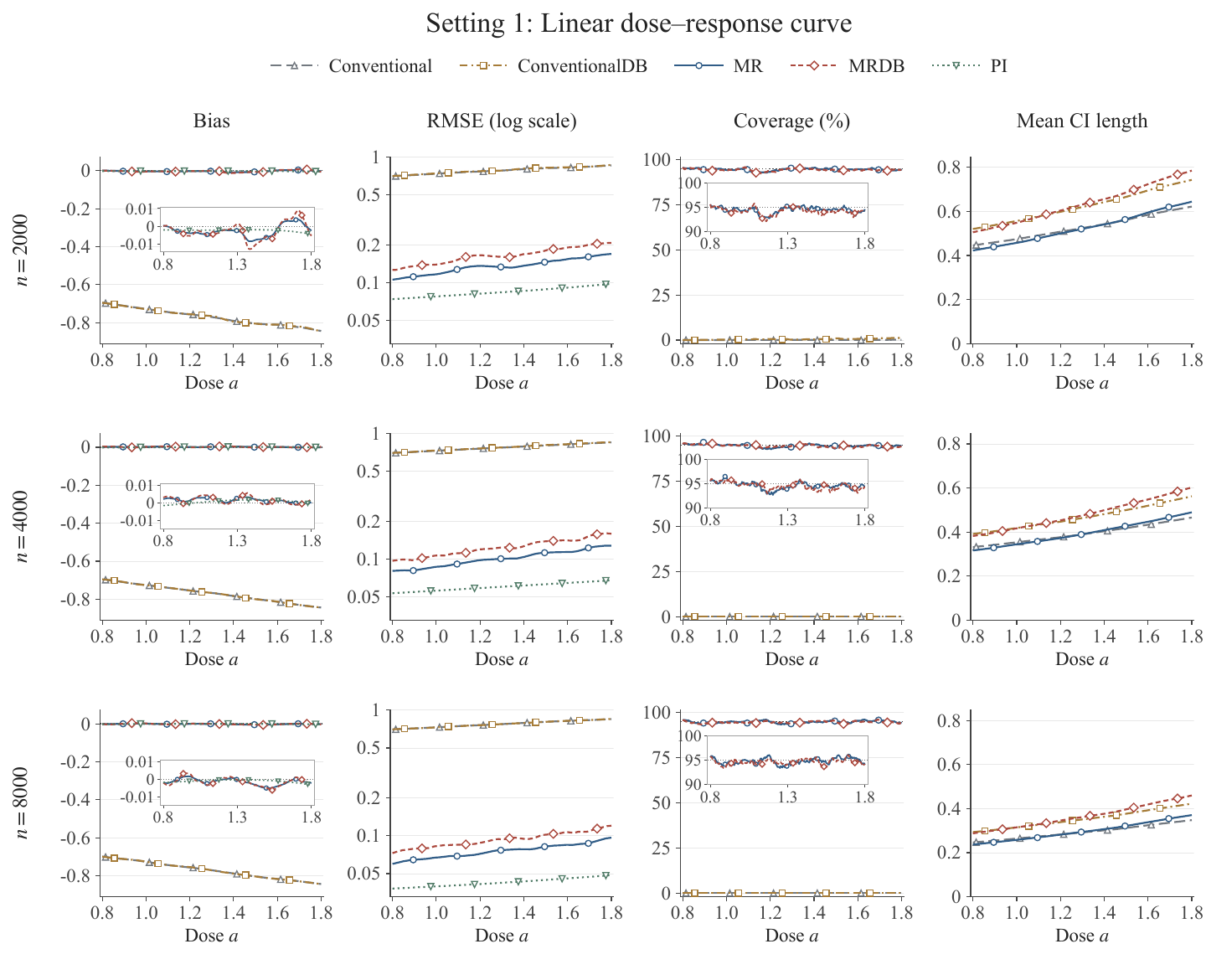}
\caption{Setting 1, undersmoothing bandwidth.}
\label{fig:sim-main-linear-undersmooth}
\end{figure}

Figure~\ref{fig:sim-main-linear-undersmooth} gives the
corresponding undersmoothing results. MRDB retains small bias
and coverage generally close to the nominal level, but its
RMSE and interval length are larger than under the regular
bandwidth. The regular bandwidth provides more precise
inference in this setting.

\begin{figure}[ht]
\centering
\linespread{1}\selectfont
\includegraphics[
  width=\textwidth,
  height=0.52\textheight,
  keepaspectratio
]{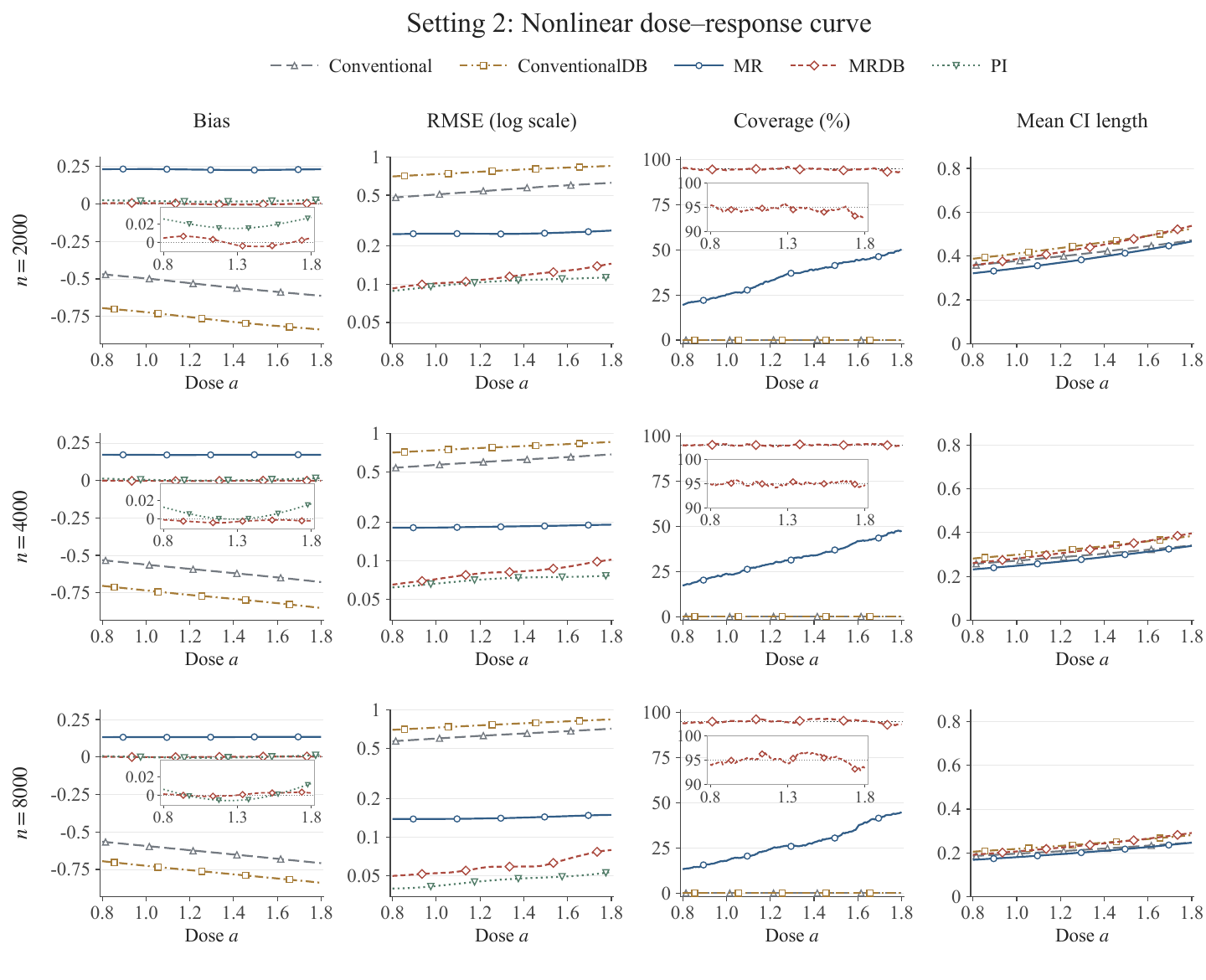}
\caption{Setting 2, regular bandwidth.}
\label{fig:sim-main-nonlinear-regular}
\end{figure}

The benefit of bias correction is most evident in Setting 2
(Figure~\ref{fig:sim-main-nonlinear-regular}). With the regular
bandwidth, MR exhibits substantial undercoverage because smoothing
bias remains appreciable relative to its standard error, even as
the sample size increases. MRDB removes most of this bias and
maintains coverage close to the nominal level.

\begin{figure}[ht]
\centering
\linespread{1}\selectfont
\includegraphics[
  width=\textwidth,
  height=0.52\textheight,
  keepaspectratio
]{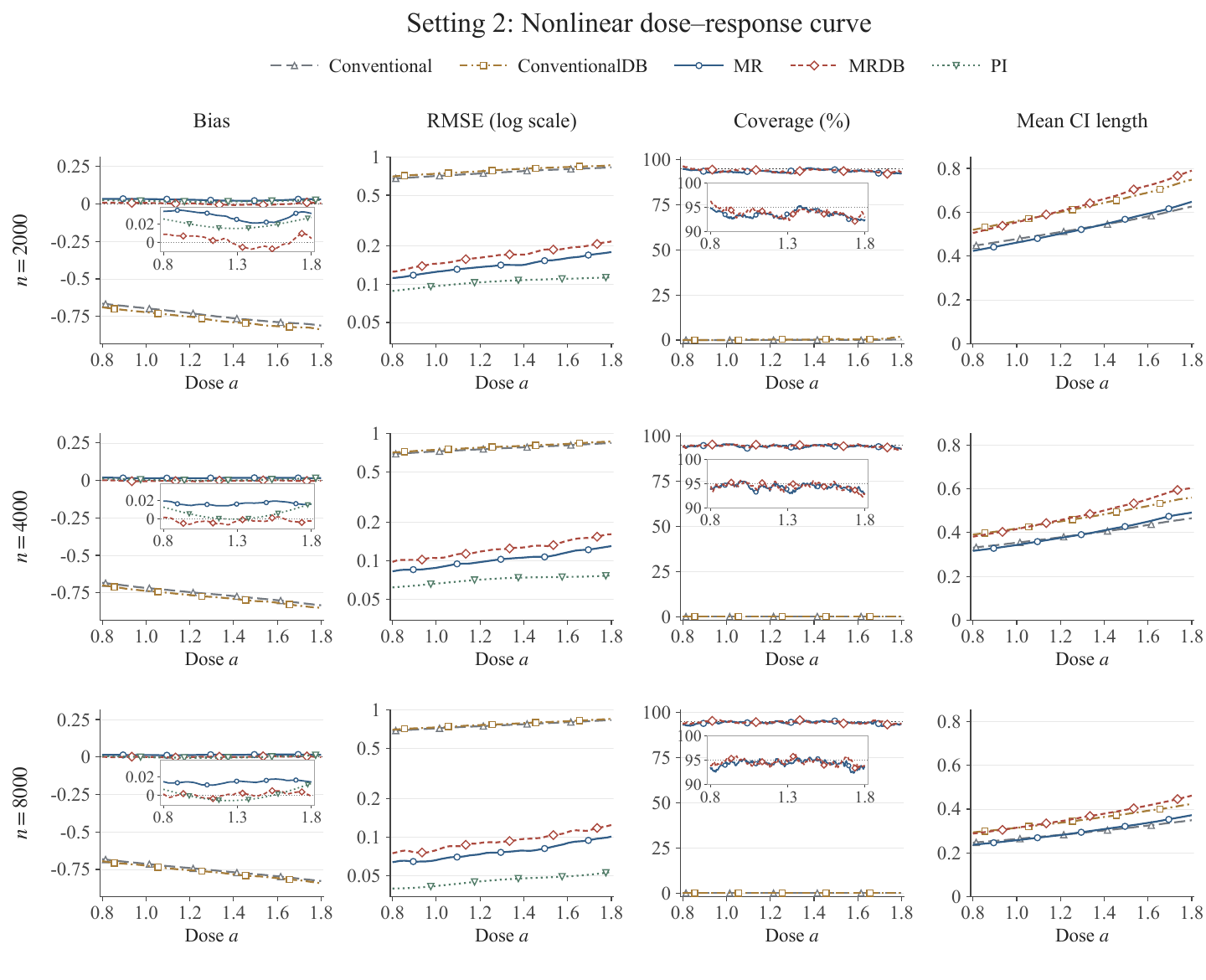}
\caption{Setting 2, undersmoothing bandwidth.}
\label{fig:sim-main-nonlinear-undersmooth}
\end{figure}

Undersmoothing substantially restores MR's coverage in
Setting 2 (Figure~\ref{fig:sim-main-nonlinear-undersmooth}),
whereas MRDB achieves comparable coverage with shorter
intervals under the regular bandwidth. For example, at
$n=8000$ and $a=1.3$, regular-bandwidth MRDB has coverage
$94.3\%$ and mean interval length $0.233$, compared with
$94.7\%$ and $0.295$ for undersmoothed MR. Debiasing therefore
supports accurate inference while retaining the precision
afforded by the regular bandwidth.

Table~\ref{tab:sim-bands} reports MRDB simultaneous band
performance. Coverage ranges from $94.1\%$ to $96.0\%$ under the
regular bandwidth and from $92.4\%$ to $94.4\%$ under
undersmoothing. Band widths decrease with $n$ under both rules,
and the regular bandwidth gives narrower bands in every
configuration.

\begin{table}[ht]
\centering
\linespread{1}\selectfont
\caption{MRDB simultaneous confidence bands.}
\label{tab:sim-bands}
\begingroup
\small
\setlength{\tabcolsep}{5pt}
\renewcommand{\arraystretch}{1.08}

\begin{tabular*}{\textwidth}{@{\extracolsep{\fill}}r*{8}{r}@{}}
\toprule
\multirow{3}{*}{$n$}
& \multicolumn{4}{c}{Linear}
& \multicolumn{4}{c}{Nonlinear} \\
\cmidrule(lr){2-5}\cmidrule(lr){6-9}
& \multicolumn{2}{c}{$h_{\mathrm{r}}$}
& \multicolumn{2}{c}{$h_{\mathrm{u}}$}
& \multicolumn{2}{c}{$h_{\mathrm{r}}$}
& \multicolumn{2}{c}{$h_{\mathrm{u}}$} \\
\cmidrule(lr){2-3}\cmidrule(lr){4-5}
\cmidrule(lr){6-7}\cmidrule(lr){8-9}
& SCP & Width & SCP & Width
& SCP & Width & SCP & Width \\
\midrule
2000
& 94.1 & 0.583 & 93.4 & 0.970
& 94.9 & 0.583 & 92.4 & 0.974 \\
4000
& 95.3 & 0.434 & 94.0 & 0.749
& 96.0 & 0.436 & 94.4 & 0.753 \\
8000
& 95.8 & 0.325 & 93.5 & 0.579
& 94.9 & 0.326 & 93.5 & 0.580 \\
\bottomrule
\end{tabular*}

\par\vspace{4pt}
\begin{minipage}{\textwidth}
\footnotesize
\textit{Notes:} SCP is simultaneous coverage (\%) over all
$201$ evaluation doses. Width is averaged over doses and replications.
\end{minipage}

\endgroup
\end{table}

Appendix~\ref{sec:sim-supplement}
provides the main experiment results at $a=0.8,1.3,1.8$ in
Tables~\ref{tab:sim-main-linear} and~\ref{tab:sim-main-nonlinear}.
It also reports a separate robustness experiment in which one
nuisance model is misspecified and the other two are correctly
specified parametric models. MRDB retains small
bias and coverage close to the nominal level under each of the
three misspecification scenarios.

\section{Semi-synthetic Data Analysis}
\label{sec:application}

We illustrate the proposed methods in a semi-synthetic study calibrated to the National Health and Nutrition Examination Survey (NHANES) \citep{nhanes20032008}. NHANES has supported methodological studies of continuous dietary exposures \citep{huangzhang2023} and exposure-response relationships for blood lipid outcomes \citep{antonelli2022causal}. Following semi-synthetic evaluations of causal methods \citep{wu2024proximal, ham2026proximal}, we combine observed covariates with generated treatments and outcomes.

We use the 2003-2008 survey cycles, with total dietary fiber intake from the first 24-hour recall (g/day) as treatment and calculated low-density lipoprotein (LDL) cholesterol (mg/dL) as outcome. This choice is motivated by randomized trials showing that certain soluble fiber interventions can lower LDL cholesterol \citep{wolever2010fiber}. The calibration sample contains 4,797 adults aged 20-79 with reliable dietary recalls, eligible fasting measurements and complete covariate information. Covariates include age, sex, body mass index, energy intake, poverty-income ratio, education, smoking status, race/ethnicity and survey cycle, giving 17 variables after categorical encoding. We sample covariates from the empirical calibration distribution. Treatments are generated from a fitted conditional truncated lognormal model with a uniform component. Outcomes follow $Y(a)=\mu_0(X,a)+\varepsilon$, where $\mu_0$ is fitted using cubic dose splines, covariate terms and dose-covariate interactions, and $\varepsilon$ is independent mean-zero Gaussian noise. The known ADRF is
$
\theta(a)=E_X\{\mu_0(X,a)\},
$
with the expectation evaluated over all 4,797 source covariate records. We generate an artificial censoring threshold $C$ independently of treatment and outcome conditional on $X$, using an exponential distribution with an age-dependent rate calibrated to 35\% expected censoring. The observed data are $O=(X,U,\Delta,\Delta Y)$, with $U=\min(A,C)$ and $\Delta=I(A\le C)$; the exact treatment and outcome are unavailable when $\Delta=0$. 

\begin{figure}[ht]
  \centering
  \includegraphics[width=12cm]{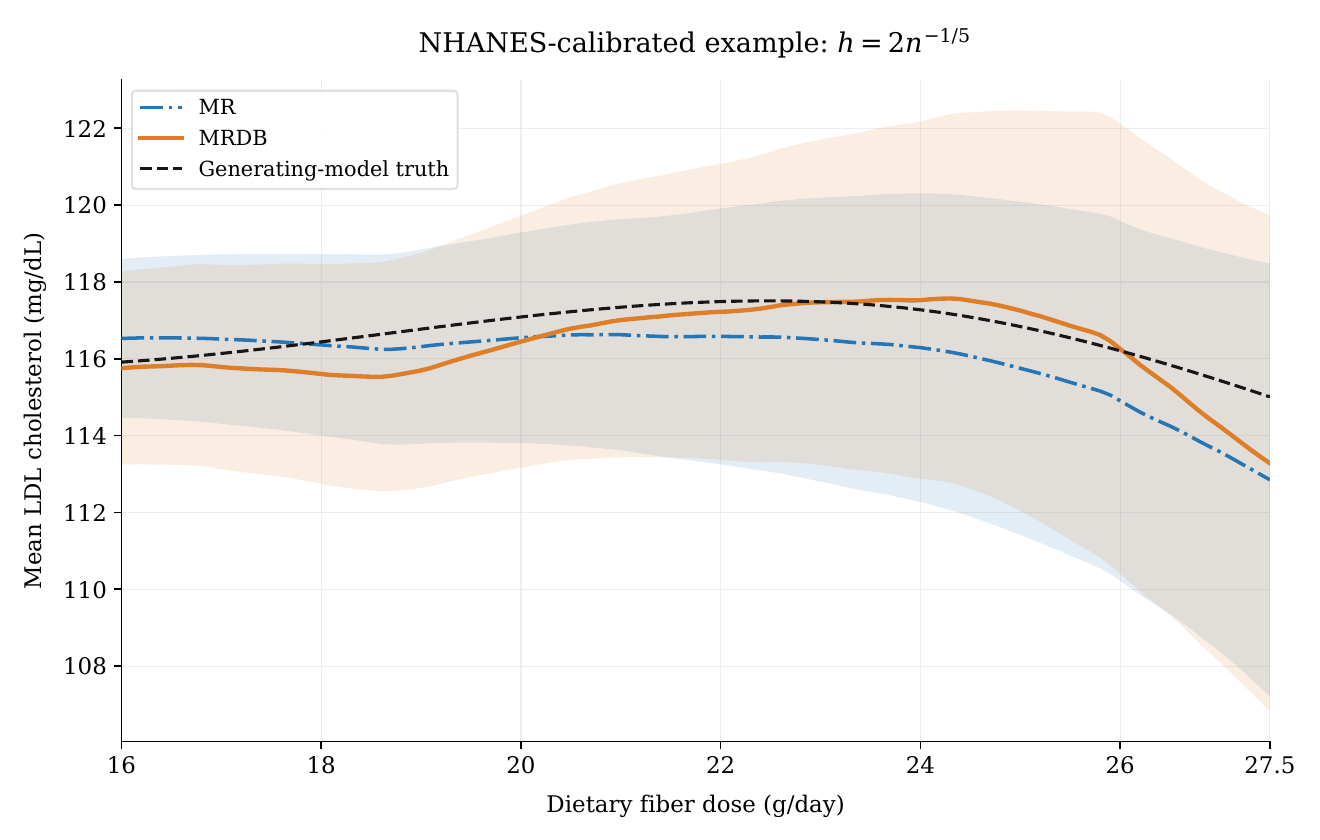}
  \caption{MR and MRDB estimates in the NHANES-calibrated
  semi-synthetic sample.}
  \label{fig:nhanes_ldl_adrf}
\end{figure}

We generate one sample of $n=4,000$, containing 2,549 complete and 1,451 censored records. MR and MRDB use the same two-fold cross-fitting scheme and nuisance estimates. As in the simulations, the outcome regression uses kernel ridge regression, and the treatment density and censoring hazard use neural network models with Bernstein bases. Both estimators use $h=2n^{-1/5}$, and MRDB uses $b=1.25h$, on the rescaled dose axis. Figure~\ref{fig:nhanes_ldl_adrf} shows that both MR and MRDB fit the true ADRF well and recover its nonlinear shape.

\section{Discussion}
\label{sec:discussion}

In this paper, we propose a novel multiply robust ADRF estimator for right-censored continuous treatments with outcomes observed only upon treatment completion, which remains consistent whenever any two of the outcome regression, treatment density, and censoring survival function are consistently estimated. We further develop a debiased estimator that removes the leading smoothing bias and enables pointwise inference and simultaneous confidence bands without undersmoothing.

Our work suggests several directions for future research. First, the outcome may remain fully observed even when treatment or exposure is censored. One possible application is to management settings where demand cannot be fully observed because of inventory limits, but revenue is still recorded \citep{tang2025offline}. Under the stated assumptions, our method remains consistent and multiply robust in this setting. However, it does not use the outcomes of censored individuals, leaving room to improve efficiency. The EIF for $\Phi$ derived in Appendix~\ref{app:always_observed_outcome} provides a starting point for constructing local estimating equations that use these additional outcomes. 
Second, our estimator is constructed through augmented local polynomial estimating equations. Future work could investigate whether this construction can be adapted to alternative regression learners, such as the Forster-Warmuth learner \citep{yang2024forster}, and characterize the resulting estimation and inference properties.
Third, treatment censoring is one form of incomplete treatment information. Drawing on recent work on hidden variables \citep{zhoutchetgen2024hidden,bian2026hidden,guo2026proximal}, we could explore settings with hidden exposures. Such extensions would require suitable proxy variables and assumptions that ensure identification. Finally, we could extend our method to policy learning and dose-finding \citep{chen2016personalized}. The goal would be to choose treatment durations based on individual characteristics to optimize expected outcomes while accounting for treatment censoring.

 \subsection*{Conflict of interest}

    None declared. 

    \subsection*{Funding}

    The authors were partially supported by the National Key R\&D Program of China (2024YFA1015600), the National Natural Science Foundation of China (12471266 and U23A2064). 

    \subsection*{Data availability}

The dataset is available at \href{https://www.cdc.gov/nchs/nhanes/}{https://www.cdc.gov/nchs/nhanes/}. The computer code is available at \href{https://github.com/SihanWu03/ADRF-with-censored-treatment}{https://github.com/SihanWu03/ADRF-with-censored-treatment}.

\putbib[bibfile_main] 
\end{bibunit}	
\defaultbibliographystyle{plainnat}  
\setcitestyle{round} 
\defaultbibliography{bibfile_supp} 
\begin{bibunit}

\newpage
\appendix
\section*{Appendix}
\section{Setup}
\subsection{Pointwise Identification}
\label{app:pointwise_identification}

For each \(a\in\mathcal A\), we assume that the conditional mean
\(t\mapsto E\{Y(a)\mid X,A=t\}\) has a version continuous at \(t=a\),
with value \(\mu(X,a)\) there, almost surely. This specifies conditional
mean consistency at continuous treatment levels; see
\citet{GillRobins2001} for a detailed discussion of continuous-treatment
identification.

Unconfoundedness implies that this conditional mean equals
\(E\{Y(a)\mid X\}\) for almost every \(t\). Treatment positivity and
continuity extend this equality to \(t=a\). Consequently,
\(E\{Y(a)\mid X\}=\mu(X,a)\), and averaging over \(X\) gives
\(\theta(a)=m(a)\). Conditional independent censoring then yields
the observed-data identification formula in Section~\ref{sec:setup}.
Under the smoothness conditions imposed in the main text, pointwise
values of the marginal curve \(m\) are taken from its continuous version.
\section{Multiply Robust Estimation}
\subsection{The Observed-Data Efficient Influence Function}
\label{app:observed_eif}

\begin{proof}[Proof of Theorem~\ref{thm:observed_eif}]
We use G\^ateaux derivatives along dominated paths and the tangent-space characterization of efficient influence functions \citep{bickel1993,tsiatis2006}. All expectations are under the true distribution. Terms multiplied by $\Delta$ use $Y$ only when $\Delta=1$ and are zero otherwise. Tail integrands are zero wherever $S_A(t\mid X)=0$.

Consider bounded scores $S_Z(Z)$ and $S_{\mathrm{cen}}(C,X)$ satisfying $E(S_Z)=0$ and $E(S_{\mathrm{cen}}\mid X)=0$. For sufficiently small $\epsilon$, consider the dominated paths
\[
\frac{dP_{Z,\epsilon}}{dP_Z}=1+\epsilon S_Z,
\qquad
\frac{dP_{C\mid X,\epsilon}}{dP_{C\mid X}}=1+\epsilon S_{\mathrm{cen}}.
\]
The joint path retains conditional independent censoring. Its induced observed-data score is $S(O)=E\{S_Z(Z)+S_{\mathrm{cen}}(C,X)\mid O\}$. Write $\Phi_\epsilon=\int_{\mathcal A}m_\epsilon(t)p_{A,\epsilon}(t)\,dt$. Differentiating the conditional outcome mean, the marginal law of $X$, and the marginal treatment density gives the G\^ateaux derivative
\begin{equation}
\label{eq:eif_gateaux_full}
\begin{aligned}
\left.\frac{d\Phi_\epsilon}{d\epsilon}\right|_{\epsilon=0}
={}&E\!\left[\frac{p_A(A)}{f_A(A\mid X)}\{Y-\mu(X,A)\}S_Z(Z)\right]\\
&+E\!\left[\left\{m(A)+\int_{\mathcal A}\mu(X,t)p_A(t)\,dt-2\Phi\right\}S_Z(Z)\right]\\
={}&E\{\varphi^F(Z)S_Z(Z)\}.
\end{aligned}
\end{equation}
The centering constant $2\Phi$ does not affect the derivative because $E(S_Z)=0$. This recovers the full-data influence function in \eqref{eq:full_eif}, as in \citet{kennedy2017}. There is no censoring-score term because $\Phi$ depends only on $P_Z$.

Let $\varphi^O$ be the candidate in \eqref{eq:observed_eif}. Since the outcome residual has conditional mean zero given $(X,A)$,
\[
V(t,X)=\frac{\int_t^T m(s)f_A(s\mid X)\,ds}{S_A(t\mid X)}+R(X)
\]
on the positive tail support. Expanding the censoring martingale gives
\[
\varphi^O(O)=\frac{\Delta\varphi^F(X,U,Y)+(1-\Delta)V(U,X)}{S_C(U\mid X)}-\int_0^U\frac{V(t,X)\lambda_C(t\mid X)}{S_C(t\mid X)}\,dt.
\]
Thus, the candidate is a function of the observed data.

We record the integrability needed for the conditional calculations below. Absolute continuity and positivity imply
\[
\int_0^A\frac{\lambda_C(t\mid X)}{S_C(t\mid X)}\,dt=S_C(A\mid X)^{-1}-1\le c_2^{-1}-1.
\]
Conditional Jensen's inequality and Fubini's theorem therefore give
\[
E\!\left[\int_0^A V(t,X)^2\frac{\lambda_C(t\mid X)}{S_C(t\mid X)}\,dt\right]
\le(c_2^{-1}-1)E[\{\varphi^F(Z)\}^2].
\]
The martingale isometry bounds the second moment of the augmentation by the same quantity. The inverse-probability term has second moment at most $c_2^{-1}E[\{\varphi^F(Z)\}^2]$, so $\varphi^O\in L_2(P_O)$. Weighted Cauchy-Schwarz also bounds the second moment of $\int_0^A V\lambda_C/S_C\,dt$ by $(c_2^{-1}-1)^2E[\{\varphi^F(Z)\}^2]$.

Conditional on $Z$, the jump and compensator terms in the augmentation have the same expectation, while $E(\Delta\mid Z)=S_C(A\mid X)$. Hence
\begin{equation}
\label{eq:eif_conditional_full_data}
E\{\varphi^O(O)\mid Z\}=\varphi^F(Z),
\qquad E\{\varphi^O(O)\}=0.
\end{equation}
To account for changes in the censoring distribution, condition instead on $(X,C=r)$. Using $S_A(t\mid X)V(t,X)=E\{I(A\ge t)\varphi^F(Z)\mid X\}$ and the expanded candidate, we obtain
\[
\begin{aligned}
&E\{\varphi^O(O)\mid X,C=r\}\\
&\quad=E\!\left[\left.\varphi^F(Z)\left\{\frac{1}{S_C(A\wedge r\mid X)}-\int_0^{A\wedge r}\frac{\lambda_C(t\mid X)}{S_C(t\mid X)}\,dt\right\}\right|X\right]\\
&\quad=E\{\varphi^F(Z)\mid X\}.
\end{aligned}
\]
It follows that $E\{\varphi^O(O)S_{\mathrm{cen}}(C,X)\}=0$. Combining this with \eqref{eq:eif_conditional_full_data} and \eqref{eq:eif_gateaux_full} yields
\[
E\{\varphi^O(O)S(O)\}=E\{\varphi^F(Z)S_Z(Z)\}
=\left.\frac{d\Phi_\epsilon}{d\epsilon}\right|_{\epsilon=0}.
\]
The candidate therefore represents the derivative along both full-data and censoring directions.

Finally, we establish efficiency by showing that the candidate belongs to the closed observed-data tangent space. Define
\[
\widetilde S_Z(Z)=\frac{\varphi^F(Z)}{S_C(A\mid X)}-\int_0^A\frac{V(t,X)\lambda_C(t\mid X)}{S_C(t\mid X)}\,dt.
\]
Positivity and the preceding second-moment bound imply $\widetilde S_Z\in L_2(P_Z)$. We claim that
\begin{equation}
\label{eq:eif_preimage_identity}
E\{\widetilde S_Z(Z)\mid O\}=\varphi^O(O).
\end{equation}
This is immediate when $\Delta=1$, because $Z$ is observed. For a censored observation with $U=u$, the conditional law of $(A,Y)$ is its law given $(X,A>u)$. On this event, write
\[
\begin{aligned}
\widetilde S_Z(Z)
={}&\frac{\varphi^F(Z)}{S_C(u\mid X)}-\int_0^u\frac{V(t,X)\lambda_C(t\mid X)}{S_C(t\mid X)}\,dt\\
&+\int_u^A\frac{\{\varphi^F(Z)-V(t,X)\}\lambda_C(t\mid X)}{S_C(t\mid X)}\,dt.
\end{aligned}
\]
The last integral has conditional mean zero given $(X,A>u)$, since $E[I(A\ge t)\{\varphi^F(Z)-V(t,X)\}\mid X]=0$ for $t\ge u$. Thus,
\[
E\{\widetilde S_Z(Z)\mid X,A>u\}=\frac{V(u,X)}{S_C(u\mid X)}-\int_0^u\frac{V(t,X)\lambda_C(t\mid X)}{S_C(t\mid X)}\,dt,
\]
which equals the candidate on censored observations and proves \eqref{eq:eif_preimage_identity}. Conditional absolute continuity of $A$ justifies replacing $A>u$ by $A\ge u$ in these tail expectations. Taking expectations in \eqref{eq:eif_preimage_identity} gives $E(\widetilde S_Z)=0$.

Because the full-data model is nonparametric, $\widetilde S_Z\in L_2^0(P_Z)$ can be approximated in $L_2(P_Z)$ by bounded mean-zero full-data scores. Holding the conditional censoring distribution fixed, their observed-data scores are the corresponding conditional expectations given $O$. By conditional-expectation contraction and \eqref{eq:eif_preimage_identity}, these scores converge in $L_2(P_O)$ to $\varphi^O$. Consequently, $\varphi^O$ belongs to the closed observed-data tangent space. Together with the derivative identity, this proves that it is the canonical gradient of $\Phi$.
\end{proof}

\subsection{Triple Robustness}
\label{app: mr}
\begin{proof}[Proof of Theorem~\ref{thm:triple_robustness}]
All expectations below are taken under the true distribution, with \(\eta^*\) fixed. Define
\[
 N_h(a)=E\{K_h(A-a) z_{A,h}(a)\theta(A)\},\qquad  D_h(a)=E\{K_h(A-a) z_{A,h}(a)\{ z_{A,h}(a)\}^\top\}.
\]
We first establish the moment identities
\begin{equation}
\label{eq:mr_moment_identity}
E\{\Gamma_N^*(O,a;h)\}= N_h(a),\qquad E\{\Gamma_D^*(O,a;h)\}= D_h(a).
\end{equation}

Suppose first that \(S_C^*=S_C\) and either \(\mu^*=\mu\) or \(f_A^*=f_A\). Then \(M_C^*\) is the true censoring martingale, so the augmentation terms have mean zero. Conditional independent censoring gives
\[
E\{\Gamma_N^*(O,a;h)\}=E\{K_h(A-a) z_{A,h}(a)\xi^*(X,A,Y)\},\qquad E\{\Gamma_D^*(O,a;h)\}= D_h(a).
\]
By the definitions of \(m^*\) and \(p_A^*\),
\[
E\{\xi^*(X,A,Y)\mid A=t\}-\theta(t)=E_X\!\left[\{\mu^*(X,t)-\mu(X,t)\}\left\{1-\frac{p_A^*(t)f_A(t\mid X)}{p_A(t)f_A^*(t\mid X)}\right\}\right].
\]
The right-hand side is zero when \(\mu^*=\mu\), or when \(f_A^*=f_A\), in which case \(p_A^*=p_A\). Thus \(E(\xi^*\mid A=t)=\theta(t)\), proving \eqref{eq:mr_moment_identity} in both cases.

Now suppose that \(\mu^*=\mu\) and \(f_A^*=f_A\), while \(S_C^*\) may be incorrect. The outcome residual has conditional mean zero given \((X,A)\), so \(E(\xi^*\mid X,A)=m(A)=\theta(A)\). Conditional independent censoring therefore allows us to replace \(\xi^*\) by \(\theta(A)\) when computing the expectation of the weighted response contribution. Write
\[
 F(t,X)=S_A(t\mid X) H_N^*(t,X;a,h)=\int_t^\infty K_h(s-a) z_{s,h}(a)\theta(s)f_A(s\mid X)\,ds.
\]
Combining the conditional expectations of the censoring jump and compensator terms gives
\[
\begin{aligned}
E\{\Gamma_N^*(O,a;h)\mid X\}
&=\int_0^{a+h}\frac{S_C(t\mid X)}{S_C^*(t\mid X)}K_h(t-a) z_{t,h}(a)\theta(t)f_A(t\mid X)\,dt\\
&\quad-\int_0^{a+h} F(t,X)\,d\!\left\{\frac{S_C(t\mid X)}{S_C^*(t\mid X)}\right\}.
\end{aligned}
\]
Here \( F(a+h,X)=0\), \(S_C(0\mid X)=S_C^*(0\mid X)=1\), and \(d_t F(t,X)=-K_h(t-a) z_{t,h}(a)\theta(t)f_A(t\mid X)\,dt\). Integration by parts yields
\[
-\int_0^{a+h} F(t,X)\,d\!\left\{\frac{S_C(t\mid X)}{S_C^*(t\mid X)}\right\}= F(0,X)-\int_0^{a+h}\frac{S_C(t\mid X)}{S_C^*(t\mid X)}K_h(t-a) z_{t,h}(a)\theta(t)f_A(t\mid X)\,dt.
\]
Consequently, \(E\{\Gamma_N^*(O,a;h)\mid X\}= F(0,X)\), and averaging over \(X\) gives the first identity in \eqref{eq:mr_moment_identity}. Applying the same calculation to \(S_A(t\mid X) H_D^*(t,X;a,h)\) gives
\[
E\{\Gamma_D^*(O,a;h)\mid X\}=E\{K_h(A-a) z_{A,h}(a)\{ z_{A,h}(a)\}^\top\mid X\},
\]
which proves the second identity. This covers all three configurations in the theorem.

We next control the sample averages as \(h\to0\). For all sufficiently large \(n\), \([a-h,a+h]\subset I\) and \(a+h<\tau\). The working density bound implies \(S_A^*(t\mid x)\ge c(T-\tau)>0\) for \(t\le\tau\). Since \(\| z_{s,h}(a)\|\le\sqrt{2}\) on the kernel support, both working tail functions are uniformly bounded and vanish for \(t>a+h\). Their augmentation terms are also uniformly bounded because
\[
\int_0^\tau\frac{\lambda_C^*(t\mid x)}{S_C^*(t\mid x)}\,dt=\frac{1}{S_C^*(\tau\mid x)}-1\le c^{-1}-1.
\]
The bounds on the working functions and the conditional second moment of \(Y\) ensure a uniformly bounded conditional second moment of \(\xi^*\) on \(I\). Thus the point contributions have second moments of order \(h^{-1}\), since \(\int K_h(t-a)^2p_A(t)\,dt=O(h^{-1})\). It follows that
\[
E\|\Gamma_N^*(O,a;h)\|^2+E\|\Gamma_D^*(O,a;h)\|^2=O(h^{-1}).
\]
Independence of the observations, Chebyshev's inequality, and \eqref{eq:mr_moment_identity} now imply
\[
 N_h^*(a)- N_h(a)=O_p\{(nh)^{-1/2}\},\qquad  D_h^*(a)- D_h(a)=O_p\{(nh)^{-1/2}\}.
\]

Finally, continuity of \(\theta\) and \(p_A\), symmetry of \(K\), and the change of variables \(t=a+hu\) give
\[
 D_h(a)\longrightarrow p_A(a)\begin{pmatrix}1&0\\0&\nu_2\end{pmatrix},\qquad  N_h(a)\longrightarrow p_A(a)\theta(a)e_1,
\]
where \(\nu_2=\int u^2K(u)\,du>0\). Because \(p_A(a)>0\) and \(nh\to\infty\), \( D_h^*(a)\) is invertible with probability tending to one. The continuous mapping theorem therefore yields \(\beta_h^*(a)\xrightarrow{p}(\theta(a),0)^\top\), and taking its first coordinate proves \(\theta_h^*(a)\xrightarrow{p}\theta(a)\). The definition on the singular event does not affect this convergence.
\end{proof}

\subsection{Rate of Convergence}
\begin{proof}[Proof of Theorem~\ref{thm: rate1}]
Let \(\mathbb P\) denote expectation under the true distribution. Write the local Taylor coefficients and residual score as
\[
\beta_h(a)=\{\theta(a),h\theta'(a)\}^{\top},
\qquad
\Psi(O,a;h,\eta)
=\Gamma_N(O,a;h,\eta)
-\Gamma_D(O,a;h,\eta)\beta_h(a).
\]
Evaluation at \(\hat\eta_{-k}\) includes the fitted marginal
functions \(\hat m_{-k}\) and \(\hat p_{A,-k}\).
Define
\begin{align*}
 S_{n,oracle}
&=\mathbb P_n\{\Psi(O,a;h,\bar\eta)\},\\
 R_{n,1}
&=\frac1K\sum_{k=1}^K(\mathbb P_{n,k}-\mathbb P)
\{\Psi(O,a;h,\hat\eta_{-k})-\Psi(O,a;h,\bar\eta)\},\\
 R_{n,2}
&=\frac1K\sum_{k=1}^K\mathbb P
\{\Psi(O,a;h,\hat\eta_{-k})-\Psi(O,a;h,\bar\eta)\}.
\end{align*}
The equal fold sizes give \(K^{-1}\sum_k\mathbb P_{n,k}=\mathbb P_n\).
Since \(\hat{ D}_h(a)\hat\beta_h(a)=\hat{ N}_h(a)\),
on the event that \(\hat{ D}_h(a)\) is invertible,
\begin{equation}
\label{eq:rate_decomposition}
\hat\theta_h(a)-\theta(a)
=e_1^\top\{\hat{ D}_h(a)\}^{-1}
\{ S_{n,oracle}+ R_{n,1}+ R_{n,2}\}.
\end{equation}
We establish the required matrix stability below.

Let
\( D_h(a)=\mathbb P\{\Gamma_D(O,a;h,\bar\eta)\}\).
Condition~\ref{cond:rate:correctness} and the moment identities
\eqref{eq:mr_moment_identity}, established in the proof of
Theorem~\ref{thm:triple_robustness}, imply
\begin{equation}
\label{eq:rate_design_limit}
\begin{aligned}
 D_h(a)
&=\int K(u)
\begin{pmatrix}1&u\\u&u^2\end{pmatrix}p_A(a+hu)\,du\\
&\longrightarrow
p_A(a)\begin{pmatrix}1&0\\0&\nu_2\end{pmatrix},
\qquad \nu_2=\int u^2K(u)\,du>0.
\end{aligned}
\end{equation}
Here we used Conditions~\ref{cond:rate:kernel},
\ref{cond:rate:smoothness}, and \ref{cond:rate:bandwidth}.
In particular, \(\|\{ D_h(a)\}^{-1}\|=O(1)\).

\noindent\textbf{1. Oracle score.}
The same moment identities give
\[
\mathbb E\{\Psi(O,a;h,\bar\eta)\}
=\int K_h(t-a) z_{t,h}(a)
\{\theta(t)-\theta(a)-\theta'(a)(t-a)\}p_A(t)\,dt.
\]
Taylor expansion, the change of variables \(t=a+hu\), and
symmetry of \(K\) yield the local linear bias expansion
\citep{FanGijbels1996}
\begin{equation}
\label{eq:rate_oracle_mean}
\mathbb E\{\Psi(O,a;h,\bar\eta)\}
=\frac12h^2\theta''(a)p_A(a)
\begin{pmatrix}\nu_2\\0\end{pmatrix}
+o(h^2).
\end{equation}

To bound the second moment, define
\begin{align*}
W_1(u,X,Y;\bar\eta)
&=\frac{Y-\bar\mu(X,u)}{\bar f_A(u\mid X)}\bar p_A(u)
+\bar m(u)-\{ z_{u,h}(a)\}^\top\beta_h(a),\\
W_2(s,X,t;\bar\eta)
&=\{\bar m(s)-\{ z_{s,h}(a)\}^\top\beta_h(a)\}
\frac{\bar f_A(s\mid X)}{\bar S_A(t\mid X)}.
\end{align*}
Then \(\Psi(O,a;h,\bar\eta)=T_1+T_2-T_3\), where
\begin{align*}
T_1
&=\frac{\Delta}{\bar S_C(U\mid X)}
K_h(U-a) z_{U,h}(a) W_1(U,X,Y;\bar\eta),\\
T_2
&=\frac{1-\Delta}{\bar S_C(U\mid X)}
\int_U^\infty K_h(s-a) z_{s,h}(a)
W_2(s,X,U;\bar\eta)\,ds,\\
T_3
&=\int_0^U\frac{\bar\lambda_C(t\mid X)}{\bar S_C(t\mid X)}
\left\{\int_t^\infty K_h(s-a) z_{s,h}(a)
W_2(s,X,t;\bar\eta)\,ds\right\}dt.
\end{align*}
Point terms are zero outside the kernel window or when
\(\Delta=0\); tail terms are zero for \(t>a+h\).
Fitted versions of these quantities are obtained by replacing
the limiting functions by their fitted counterparts.

For all sufficiently large \(n\), \([a-h,a+h]\subset I\) and
\(a+h<\tau\). Condition~\ref{cond:rate:density} gives
\[
\inf_{t\in[0,\tau],\,x\in\mathcal X}\bar S_A(t\mid x)
\ge c(T-\tau)>0.
\]
The density and moment bounds imply that \(W_2\) is uniformly
bounded on its integration region and
\[
\sup_{|u-a|\le h,\,x\in\mathcal X}
\mathbb E\{W_1(u,x,Y;\bar\eta)^2\mid X=x,A=u\}\le C_0.
\]
Since \(\| z_{u,h}(a)\|\le\sqrt2\) on the kernel support,
Condition~\ref{cond:rate:censoring} yields
\[
\mathbb E\|T_1\|^2
\le C_0\int K_h(u-a)^2p_A(u)\,du=O(h^{-1}).
\]
The bounded \(L_1\) norm of the kernel gives \(\|T_2\|\le C_0\).
The survival-hazard relation gives
\[
\|T_3\|
\le C_0\int_0^{a+h}
\frac{\bar\lambda_C(t\mid X)}{\bar S_C(t\mid X)}\,dt
=C_0\{\bar S_C(a+h\mid X)^{-1}-1\}\le C_0.
\]
Thus
\[
\mathbb E\|\Psi(O,a;h,\bar\eta)\|^2=O(h^{-1}),\qquad
\mathbb E\|(\mathbb P_n-\mathbb P)\Psi(O,a;h,\bar\eta)\|^2
=O((nh)^{-1}).
\]
Together with \eqref{eq:rate_oracle_mean}, this proves
\begin{equation}
\label{eq:rate_oracle_bound}
\| S_{n,oracle}\|
=O_p\{h^2+(nh)^{-1/2}\}.
\end{equation}

\noindent\textbf{2. Empirical process error.}
For a fixed fold \(k\), condition on its training data
\(\mathcal I_{-k}\) and suppress the fold index on fitted
functions. Put
\[
\Delta\Psi=\Psi(O,a;h,\hat\eta_{-k})-\Psi(O,a;h,\bar\eta).
\]
Let \(\delta_n\) be the sum of the uniform errors relative to
the limiting functions in Condition~\ref{cond:rate:limits}
and the two marginalization errors in
Condition~\ref{cond:rate:marginal-error}.
By Condition~\ref{cond:rate:marginal-definition},
the local differences \(|\hat m-\bar m|\) and
\(|\hat p_A-\bar p_A|\) are bounded by \(\delta_n=o_p(1)\).
On events with probability tending to one, the inverse-weight
bounds imply
\[
\|T_1(\hat\eta)-T_1(\bar\eta)\|
\le C_0\delta_n\Delta |K_h(A-a)|(1+|Y|).
\]
Consequently,
\[
\mathbb E\{\|T_1(\hat\eta)-T_1(\bar\eta)\|^2
\mid\mathcal I_{-k}\}
\le C_0\delta_n^2h^{-1}.
\]
The preceding tail bounds also hold for the fitted functions:
\(\hat S_A(t\mid x)\ge c(T-\tau)\) for \(t\le\tau\), and
\[
\int_0^{a+h}\frac{\hat\lambda_C(t\mid X)}
{\hat S_C(t\mid X)}\,dt
=\hat S_C(a+h\mid X)^{-1}-1\le C_0.
\]
Thus the fitted and limiting versions of \(T_2,T_3\) are
uniformly bounded, and
\[
\mathbb E(\|\Delta\Psi\|^2\mid\mathcal I_{-k})
\le C_0\{\delta_n^2h^{-1}+1\}=o_p(h^{-1}).
\]
Within the evaluation fold, the \(n_k=n/K\) centered observations are
conditionally independent. Hence
\[
\mathbb E\!\left[
\|(\mathbb P_{n,k}-\mathbb P)\Delta\Psi\|^2
\mid\mathcal I_{-k}\right]
\le \frac1{n_k}
\mathbb E(\|\Delta\Psi\|^2\mid\mathcal I_{-k})
=o_p((nh)^{-1}).
\]
Conditional Chebyshev's inequality and fixed \(K\) give
\begin{equation}
\label{eq:rate_empirical_bound}
\| R_{n,1}\|=o_p((nh)^{-1/2}).
\end{equation}

\noindent\textbf{3. Nuisance drift.}
Continue to condition on one training fold. Define
\[
H(X,A,Y;\eta)=K_h(A-a) z_{A,h}(a)W_1(A,X,Y;\eta),
\qquad
Q(t,X;\eta)=\int_t^\infty
K_h(s-a) z_{s,h}(a)W_2(s,X,t;\eta)\,ds.
\]
Write \(\hat H,\hat Q\) and \(\bar H,\bar Q\) for the fitted
and limiting versions. Conditional independent censoring gives
\[
\mathbb P\{T_1(\hat\eta)\}
=\mathbb P\!\left\{
\frac{S_C(A\mid X)}{\hat S_C(A\mid X)}\hat H(X,A,Y)
\right\}.
\]
The working compensated process satisfies
\[
d\hat M_C(t)
=dM_C(t)-I(U\ge t)
\{\hat\lambda_C(t\mid X)-\lambda_C(t\mid X)\}\,dt.
\]
The integral against the true martingale has mean zero under
the bounds already established. Therefore
\[
\mathbb P\{T_2(\hat\eta)-T_3(\hat\eta)\}
=-\mathbb E_X\int_0^{a+h}
S_A(t\mid X)\hat Q(t,X)
\frac{S_C(t\mid X)}{\hat S_C(t\mid X)}
\{\hat\lambda_C(t\mid X)-\lambda_C(t\mid X)\}\,dt.
\]
Using
\[
d\{S_C/\hat S_C\}
=(S_C/\hat S_C)(\hat\lambda_C-\lambda_C)\,dt,
\qquad
S_C(0\mid X)/\hat S_C(0\mid X)=1,
\]
and interchanging the resulting integrals yields
\[
\mathbb P\{\Psi(O,a;h,\hat\eta)\}
=\mathbb P\hat H+J_n(a),
\qquad
J_n(a)=\mathbb E_X\int_0^{a+h}
\Omega_h(t,X;a)\,d\!\left\{\frac{S_C(t\mid X)}
{\hat S_C(t\mid X)}\right\},
\]
where
\[
\Omega_h(t,X;a)=\mathbb E[
I(A\ge t)\{\hat H(X,A,Y)-\hat Q(t,X)\}\mid X].
\]
The corresponding integral for \(\bar\eta\) is zero.
Indeed, if \(\bar S_C=S_C\), its integrating differential
vanishes. Otherwise, Condition~\ref{cond:rate:correctness}
gives \(\bar f_A=f_A\) and \(\bar\mu=\mu\), so the residual
has conditional mean zero and
\[
\mathbb E[I(A\ge t)\{\bar H(X,A,Y)-\bar Q(t,X)\}\mid X]=0.
\]
Thus the single-fold drift is
\(\mathbb P(\hat H-\bar H)+J_n(a)\).

To bound \(J_n(a)\), put
\[
\hat F_h(t,X)=\int_t^T K_h(s-a) z_{s,h}(a)
\{\hat m(s)-\{ z_{s,h}(a)\}^\top\beta_h(a)\}\hat f_A(s\mid X)\,ds.
\]
Conditioning on \((X,A)\) gives the exact decomposition
\begin{align*}
\Omega_h(t,X;a)
={}&\left\{1-\frac{S_A(t\mid X)}{\hat S_A(t\mid X)}\right\}
\hat F_h(t,X)\\
&+\int_t^T K_h(s-a) z_{s,h}(a)
\Bigg[
\{\hat m(s)-\{ z_{s,h}(a)\}^\top\beta_h(a)\}
\{f_A(s\mid X)-\hat f_A(s\mid X)\}\\
&\hspace{38mm}
-\{\hat\mu(X,s)-\mu(X,s)\}
\frac{\hat p_A(s)f_A(s\mid X)}{\hat f_A(s\mid X)}
\Bigg]ds.
\end{align*}
Let \(\operatorname{TV}\) denote total variation on
\([0,a+h]\), using the Euclidean norm for vector-valued
functions. Since
\[
\int_0^T |K_h(s-a)|\,\| z_{s,h}(a)\|\,ds
=\int_{-1}^1|K(u)|\sqrt{1+u^2}\,du<\infty,
\]
the final integral in this decomposition has total variation
\(O_p\{r_n(a)+s_n(a)\}\), uniformly in \(x\). Also,
\[
\sup_{x\in\mathcal X}
\left\{\sup_{0\le t\le a+h}\|\hat F_h(t,x)\|
+\operatorname{TV}\{\hat F_h(\cdot,x)\}\right\}=O_p(1).
\]
The common support and density normalization imply
\[
\sup_{t\in[0,\tau],\,x\in\mathcal X}
|\hat S_A(t\mid x)-S_A(t\mid x)|
\le T\sup_{s\in[0,T],\,x\in\mathcal X}
|\hat f_A(s\mid x)-f_A(s\mid x)|
=O_p\{r_n(a)\}.
\]
Almost everywhere,
\[
\partial_t\left\{\frac{S_A(t\mid x)}{\hat S_A(t\mid x)}\right\}
=\frac{\hat f_A(t\mid x)-f_A(t\mid x)}
{\hat S_A(t\mid x)}
+\frac{\{S_A(t\mid x)-\hat S_A(t\mid x)\}
\hat f_A(t\mid x)}{\hat S_A(t\mid x)^2}.
\]
The lower bound \(\hat S_A(t\mid x)\ge c(T-\tau)\) therefore
gives
\[
\sup_{x\in\mathcal X}\left[
\sup_{0\le t\le a+h}
\left|1-\frac{S_A(t\mid x)}{\hat S_A(t\mid x)}\right|
+\operatorname{TV}\left\{
\frac{S_A(\cdot\mid x)}{\hat S_A(\cdot\mid x)}\right\}
\right]=O_p\{r_n(a)\}.
\]
The product inequality for total variation now yields
\[
\sup_{x\in\mathcal X}\operatorname{TV}\{\Omega_h(\cdot,x;a)\}
=O_p\{r_n(a)+s_n(a)\}.
\]
Because \(\Omega_h(a+h,X;a)=0\) and both censoring survival functions
equal one at zero, integration by parts gives
\[
J_n(a)
=-\mathbb E_X\int_0^{a+h}
\left\{\frac{S_C(t\mid X)}{\hat S_C(t\mid X)}-1\right\}
\,d\Omega_h(t,X;a).
\]
Consequently,
\begin{align*}
\|J_n(a)\|
&\le
\sup_{t\in[0,\tau],\,x\in\mathcal X}
\left|\frac{S_C(t\mid x)}{\hat S_C(t\mid x)}-1\right|
\mathbb E_X[\operatorname{TV}\{\Omega_h(\cdot,X;a)\}]\\
&=O_p[q_n(a)\{r_n(a)+s_n(a)\}].
\end{align*}
This bound uses the survival-function error and does not
require a convergence rate for the censoring hazard.

It remains to bound \(\mathbb P(\hat H-\bar H)\).
Write
\[
e_m(t)=\hat m(t)-\mathbb E_X\{\hat\mu(X,t)\},
\qquad
e_p(t)=\hat p_A(t)-\mathbb E_X\{\hat f_A(t\mid X)\}.
\]
Condition~\ref{cond:rate:marginal-error} gives
\[
\sup_{|t-a|\le h}\{|e_m(t)|+|e_p(t)|\}=O_p(n^{-1/2}).
\]
By Conditions~\ref{cond:rate:correctness} and
\ref{cond:rate:marginal-definition}, either \(\bar\mu=\mu\)
or \(\bar f_A=f_A\), and the corresponding limiting
full-data bias is zero. Conditioning on \((X,A)\) then gives
\begin{align*}
\mathbb P(\hat H-\bar H)
=\int K_h(t-a) z_{t,h}(a)\Bigg[
&\mathbb E_X\!\left\{
(\hat\mu(X,t)-\mu(X,t))
\left(p_A(t)-\hat p_A(t)
\frac{f_A(t\mid X)}{\hat f_A(t\mid X)}\right)
\right\}\\
&+p_A(t)e_m(t)\Bigg]dt.
\end{align*}
The density factor satisfies
\[
p_A(t)-\hat p_A(t)\frac{f_A(t\mid x)}{\hat f_A(t\mid x)}
=
\frac{
p_A(t)\{\hat f_A(t\mid x)-f_A(t\mid x)\}
-f_A(t\mid x)\{\hat p_A(t)-p_A(t)\}}
{\hat f_A(t\mid x)},
\]
where
\[
\sup_{|t-a|\le h}|\hat p_A(t)-p_A(t)|
\le
\sup_{t\in[0,T],\,x\in\mathcal X}
|\hat f_A(t\mid x)-f_A(t\mid x)|
+\sup_{|t-a|\le h}|e_p(t)|.
\]
The uniform bounds and the bounded \(L_1\) norm of
\(K_h(t-a) z_{t,h}(a)\) imply
\[
\|\mathbb P(\hat H-\bar H)\|
=O_p\{r_n(a)s_n(a)+n^{-1/2}\}.
\]
Averaging the single-fold drift over the fixed number of
folds proves
\begin{equation}
\label{eq:rate_drift_bound}
\| R_{n,2}\|
=O_p\{q_n(a)r_n(a)+q_n(a)s_n(a)
+r_n(a)s_n(a)+n^{-1/2}\}.
\end{equation}

\noindent\textbf{4. Matrix stability and conclusion.}
Apply the preceding integration-by-parts argument to the
design-matrix score, replacing the residual-weighted kernel
by \(K_h(t-a) z_{t,h}(a)\{ z_{t,h}(a)\}^\top\).
There is no regression or marginalization term. For each fold,
\begin{align*}
&\left\|
\mathbb P\{\Gamma_D(O,a;h,\hat\eta_{-k})\}
- D_h(a)\right\|\\
&\quad\le C_0
\sup_{t\in[0,\tau],\,x\in\mathcal X}
|\hat S_{C,-k}(t\mid x)-S_C(t\mid x)|
\sup_{t\in[0,T],\,x\in\mathcal X}
|\hat f_{A,-k}(t\mid x)-f_A(t\mid x)|
=o_p(1).
\end{align*}
The last step follows from Conditions~\ref{cond:rate:limits}
and \ref{cond:rate:correctness}: at least one of
\(\bar S_C,\bar f_A\) is correct, and both errors are bounded
in probability. The same point and tail bounds used above
give conditional second moments \(O_p(h^{-1})\) for the
design-matrix score entries. Their centered fold averages are
therefore \(O_p((nh)^{-1/2})\). By
\eqref{eq:rate_design_limit} and
Condition~\ref{cond:rate:bandwidth},
\begin{equation}
\label{eq:rate_matrix_stability}
\|\hat{ D}_h(a)- D_h(a)\|=o_p(1),
\qquad
\|\{\hat{ D}_h(a)\}^{-1}\|=O_p(1),
\qquad
\|\{\hat{ D}_h(a)\}^{-1}-\{ D_h(a)\}^{-1}\|=o_p(1).
\end{equation}
In particular, the empirical matrix is invertible with
probability tending to one.

Combining \eqref{eq:rate_decomposition},
\eqref{eq:rate_oracle_bound}, \eqref{eq:rate_empirical_bound},
and \eqref{eq:rate_drift_bound} gives
\[
|\hat\theta_h(a)-\theta(a)|
=O_p\!\left\{
(nh)^{-1/2}+h^2+q_n(a)r_n(a)
+q_n(a)s_n(a)+r_n(a)s_n(a)
\right\},
\]
where \(n^{-1/2}=o((nh)^{-1/2})\) because \(h\to0\).
This proves the result.
\end{proof}

\subsection{Asymptotic Normality}
\begin{proof}[Proof of Theorem~\ref{thm:asymptotic_normality}]
Use the notation of the preceding proof. By
\eqref{eq:rate_empirical_bound}, \eqref{eq:rate_drift_bound},
and Condition~\ref{cond:clt:products},
\[
\sqrt{nh}\,\| R_{n,1}\|=o_p(1),
\qquad
\sqrt{nh}\,\| R_{n,2}\|=o_p(1).
\]
The marginalization contribution is negligible because
\(\sqrt{nh}\,n^{-1/2}=\sqrt h\to0\).
Furthermore, \eqref{eq:rate_oracle_bound} and the bandwidth
requirement \eqref{eq:clt_bandwidth} imply
\[
\sqrt{nh}\,\| S_{n,oracle}\|
=O_p\{1+\sqrt{nh^5}\}=O_p(1).
\]
Together with \eqref{eq:rate_matrix_stability}, this gives
\[
e_1^\top(\{\hat{ D}_h(a)\}^{-1}-\{ D_h(a)\}^{-1})
\sqrt{nh}\, S_{n,oracle}=o_p(1).
\]
Thus the exact decomposition \eqref{eq:rate_decomposition}
reduces to
\[
\sqrt{nh}\{\hat\theta_h(a)-\theta(a)\}
=e_1^\top\{ D_h(a)\}^{-1}
\sqrt{nh}\, S_{n,oracle}+o_p(1).
\]

By \eqref{eq:rate_design_limit},
\(e_1^\top\{ D_h(a)\}^{-1}\to(p_A(a)^{-1},0)\).
The bias expansion \eqref{eq:rate_oracle_mean} therefore gives
\[
e_1^\top\{ D_h(a)\}^{-1}
\mathbb E\{\Psi(O,a;h,\bar\eta)\}=b_h(a)+o(h^2).
\]
The bandwidth requirement \eqref{eq:clt_bandwidth} makes the
remainder negligible after multiplication by \(\sqrt{nh}\).
Consequently,
\begin{equation}
\label{eq:centered_clt}
\begin{aligned}
&\sqrt{nh}\{\hat\theta_h(a)-\theta(a)-b_h(a)\}\\
&\quad=e_1^\top\{ D_h(a)\}^{-1}
\frac1{\sqrt n}\sum_{i=1}^n
\sqrt h\,[\Psi(O_i,a;h,\bar\eta)
-\mathbb E\{\Psi(O_i,a;h,\bar\eta)\}]
+o_p(1).
\end{aligned}
\end{equation}

We next calculate the covariance of the centered oracle score.
Let
\[
 V_h(a)=\operatorname{Var}\{\sqrt h\,\Psi(O,a;h,\bar\eta)\}.
\]
Since the mean is \(O(h^2)\),
\[
 V_h(a)
=h\,\mathbb E\{\Psi(O,a;h,\bar\eta)\Psi(O,a;h,\bar\eta)^\top\}
+o(1).
\]
Recall \(\Psi=T_1+T_2-T_3\). The rate proof gives
\(\mathbb E\|T_1\|^2=O(h^{-1})\), while \(T_2,T_3\)
are uniformly bounded. Hence, for \(j=2,3\),
\[
\|h\,\mathbb E(T_1T_j^\top)\|
\le h\{\mathbb E\|T_1\|^2\,\mathbb E\|T_j\|^2\}^{1/2}
=O(\sqrt h).
\]
Terms involving only \(T_2,T_3\) are \(O(h)\) after
multiplication by \(h\). Thus
\[
 V_h(a)=h\,\mathbb E(T_1T_1^\top)+o(1).
\]
On \(\{\Delta=1\}\), \(U=A\), and independent censoring gives
\(\mathbb E(\Delta\mid X,A,Y)=S_C(A\mid X)\). It follows that
\[
h\,\mathbb E(T_1T_1^\top)
=h\,\mathbb E\!\left[
\frac{S_C(A\mid X)}{\bar S_C(A\mid X)^2}
K_h(A-a)^2 z_{A,h}(a)\{ z_{A,h}(a)\}^\top
\{\bar\xi(X,A,Y)-\{ z_{A,h}(a)\}^\top\beta_h(a)\}^2
\right].
\]
Uniformly over \(|t-a|\le h\),
\(\theta(t)-\{ z_{t,h}(a)\}^\top\beta_h(a)=O(h^2)\).
The local second-moment and inverse-weight bounds therefore
give, by Cauchy-Schwarz,
\[
\mathbb E\!\left[
\frac{S_C(t\mid X)}{\bar S_C(t\mid X)^2}
\{\bar\xi(X,t,Y)-\{ z_{t,h}(a)\}^\top\beta_h(a)\}^2
\mathrel{\Big|} A=t\right]
=v(t)+O(h^2).
\]
Conditioning on \(A=t\) and setting \(t=a+hu\), we obtain
\[
 V_h(a)
=\int K(u)^2
\begin{pmatrix}1&u\\u&u^2\end{pmatrix}
p_A(a+hu)\{v(a+hu)+O(h^2)\}\,du+o(1).
\]
Condition~\ref{cond:clt:variance} and continuity of \(p_A\)
imply
\[
 V_h(a)\longrightarrow
 V(a)=p_A(a)v(a)
\int K(u)^2\begin{pmatrix}1&u\\u&u^2\end{pmatrix}\,du.
\]

For the triangular-array central limit theorem,
the bounded coefficients of \(\bar\xi\) imply
\(|\bar\xi(X,t,Y)|\le C_0(1+|Y|)\) on the kernel window.
Condition~\ref{cond:clt:moments} yields
\begin{align*}
\mathbb E\|T_1\|^{2+\delta}
&\le C_0\int |K_h(t-a)|^{2+\delta}p_A(t)
\mathbb E\{(1+|Y|)^{2+\delta}\mid A=t\}\,dt\\
&=O(h^{-1-\delta}).
\end{align*}
Since \(T_2,T_3\) are uniformly bounded, the same order
holds for \(\mathbb E\|\Psi(O,a;h,\bar\eta)\|^{2+\delta}\).
Centering and multiplying by \(\sqrt h\) therefore give
\[
\frac1{n^{\delta/2}}
\mathbb E\left\|
\sqrt h\,[\Psi(O,a;h,\bar\eta)-\mathbb E\{\Psi(O,a;h,\bar\eta)\}]
\right\|^{2+\delta}
=O((nh)^{-\delta/2})\longrightarrow0.
\]
This verifies the Lyapunov condition for every fixed linear
projection. The Lindeberg-Feller theorem and the
Cram\'er-Wold device \citep[Chapter~2]{vandervaart1998} imply
\[
\frac1{\sqrt n}\sum_{i=1}^n
\sqrt h\,[\Psi(O_i,a;h,\bar\eta)-\mathbb E\{\Psi(O_i,a;h,\bar\eta)\}]
\xrightarrow{d}\mathcal N_2\{0, V(a)\}.
\]
Finally, \(e_1^\top\{ D_h(a)\}^{-1}\to(p_A(a)^{-1},0)\).
Applying Slutsky's theorem to \eqref{eq:centered_clt} gives
\[
\sqrt{nh}\{\hat\theta_h(a)-\theta(a)-b_h(a)\}
\xrightarrow{d}\mathcal N\{0,\sigma^2(a)\},
\]
where
\[
\sigma^2(a)
=\begin{pmatrix}p_A(a)^{-1}&0\end{pmatrix}
 V(a)
\begin{pmatrix}p_A(a)^{-1}\\0\end{pmatrix}
=\frac{v(a)}{p_A(a)}\int K(u)^2\,du.
\]
\end{proof}

\section{Debiased Multiply Robust Method}
\subsection{Convergence Rate}
\begin{proof}[Proof of Theorem~\ref{thm:debiased_convergence}]
Let $\mathbb P$ denote expectation under the true distribution,
$\mathbb P_n$ the empirical average over all observations,
and $\mathbb P_{n,k}$ the empirical average over evaluation fold $k$.
In an expression involving $\hat\eta_{-k}$, $\mathbb P$ holds
the fitted functions fixed and integrates over a new observation.
Evaluation at $\hat\eta_{-k}$ includes the corresponding fitted
marginal functions.

We first decompose the estimation error.
Write
\[
\beta_h(a)=(\theta(a),h\theta'(a))^T,\qquad
\chi_n(a)=q_n(a)r_n(a)+q_n(a)s_n(a)+r_n(a)s_n(a),
\]
and recall the local linear residual score
\[
\Psi(O,a;h,\eta)
=
\Gamma_N(O,a;h,\eta)
-\Gamma_D(O,a;h,\eta)\beta_h(a).
\]
Insert the limiting score $\Psi(O,a;h,\bar\eta)$ into each fold
average, and split the fitted-minus-limiting difference into
its centered part and its expectation.
This defines the oracle average, empirical-process error,
and nuisance drift:
\begin{align*}
 S_{n,oracle}
&=\mathbb P_n[\Psi(O,a;h,\bar\eta)],\\
 R_{n,1}
&=\frac1K\sum_{k=1}^K(\mathbb P_{n,k}-\mathbb P)
[\Psi(O,a;h,\hat\eta_{-k})-\Psi(O,a;h,\bar\eta)],\\
 R_{n,2}
&=\frac1K\sum_{k=1}^K\mathbb P
[\Psi(O,a;h,\hat\eta_{-k})-\Psi(O,a;h,\bar\eta)].
\end{align*}
The equal fold sizes give
$K^{-1}\sum_k\mathbb P_{n,k}=\mathbb P_n$.
Consequently, the local linear estimating equation implies
\[
\hat\theta_h(a)-\theta(a)
=
e_1^T\{\hat{ D}_h(a)\}^{-1}
\{ S_{n,oracle}+ R_{n,1}+ R_{n,2}\}.
\]

Let
\[
 B_h(a)
=
\mathbb P[\Gamma_B(O,a;h,\bar\eta)].
\]
Under Condition~\ref{cond:rate:correctness}, either the limiting censoring model or
the limiting treatment model is correct, so
\[
 B_h(a)
=
\int K_h(t-a) z_{t,h}^{B}(a)p_A(t)\,dt.
\]
Substituting the bias correction into the preceding identity gives
\begin{equation}
\hat\theta_{\mathrm{DB}}(a)-\theta(a)
=
e_1^T\{\hat{ D}_h(a)\}^{-1}
\{ S_{\mathrm{DB}}+ R_{n,1}+ R_{n,2}\},
\label{eq:debiased_algebraic_trick}
\end{equation}
where
\[
 S_{\mathrm{DB}}
=
 S_{n,oracle}
-\frac12h^2\hat\theta''(a)\hat{ B}_h(a).
\]
The product in the correction satisfies the exact identity
\[
\hat\theta''(a)\hat{ B}_h(a)
=
\theta''(a) B_h(a)
+\{\hat\theta''(a)-\theta''(a)\}\hat{ B}_h(a)
+\theta''(a)(\hat{ B}_h(a)- B_h(a)).
\]
We also split $ S_{n,oracle}$ into its expectation and
centered fluctuation. Combining these two identities gives
\[
 S_{\mathrm{DB}}
= E_1+ E_2- E_3- E_4,
\]
with
\begin{align*}
 E_1
&= S_{n,oracle}-\mathbb E S_{n,oracle},\\
 E_2
&=\mathbb E S_{n,oracle}
-\frac12h^2\theta''(a) B_h(a),\\
 E_3
&=\frac12h^2\{\hat\theta''(a)-\theta''(a)\}
\hat{ B}_h(a),\\
 E_4
&=\frac12h^2\theta''(a)
(\hat{ B}_h(a)- B_h(a)).
\end{align*}
These terms represent, respectively, the oracle fluctuation,
the remaining smoothing bias, the derivative estimation error,
and the empirical quadratic-term moment vector error.

Condition~\ref{cond:db:marginal-error} implies the error bound in
Condition~\ref{cond:rate:marginal-error} on the smaller window $|t-a|\le h$.
Thus the local linear bounds already proved in
Theorem~\ref{thm: rate1} remain applicable:
\[
\| R_{n,1}\|=o_p((nh)^{-1/2}),
\qquad
\| R_{n,2}\|
=O_p\{\chi_n(a)+n^{-1/2}\}.
\]
For the population design matrix
\[
 D_h(a)
=
\int K_h(t-a) z_{t,h}(a)
\big( z_{t,h}(a)\big)^T p_A(t)\,dt,
\]
the same design-matrix score argument gives
\[
\|\hat{ D}_h(a)- D_h(a)\|
=
O_p\{(nh)^{-1/2}+q_n(a)r_n(a)\},
\qquad
\|\hat{ D}_h(a)- D_h(a)\|=o_p(1).
\]
Since $ D_h(a)$ has a positive definite limit,
$\|\{\hat{ D}_h(a)\}^{-1}\|=O_p(1)$.
Moreover, the exact inverse identity yields
\[
\|\{\hat{ D}_h(a)\}^{-1}-\{ D_h(a)\}^{-1}\|
=
O_p\{(nh)^{-1/2}+q_n(a)r_n(a)\}.
\]
It remains to control the four terms in $ S_{\mathrm{DB}}$.

\medskip
\noindent\textbf{1. The oracle fluctuation $ E_1$.}
The second-moment calculation in the proof of
Theorem~\ref{thm: rate1} gives
\[
\mathbb E\|\Psi(O,a;h,\bar\eta)\|^2=O(h^{-1}).
\]
Independence of the observations therefore implies
\[
\mathbb E\| E_1\|^2
\le \frac1n\mathbb E\|\Psi(O,a;h,\bar\eta)\|^2
=O((nh)^{-1}).
\]
By Chebyshev's inequality and the inverse-matrix bound,
\[
\| E_1\|=O_p((nh)^{-1/2}),
\qquad
|e_1^T\{\hat{ D}_h(a)\}^{-1} E_1|
=O_p((nh)^{-1/2}).
\]

\medskip
\noindent\textbf{2. The remaining smoothing bias $ E_2$.}
Under Condition~\ref{cond:rate:correctness}, the oracle estimating equation has the
same expectation as its true full-data counterpart. Hence
\[
\mathbb E S_{n,oracle}
=
\int K_h(t-a) z_{t,h}(a)
\{\theta(t)-\theta(a)-\theta'(a)(t-a)\}p_A(t)\,dt.
\]
The identity
$ z_{t,h}(a)((t-a)/h)^2
= z_{t,h}^{B}(a)$ gives
\[
 E_2
=
\int K_h(t-a) z_{t,h}(a)
R_2(t,a)p_A(t)\,dt,
\]
where
\[
R_2(t,a)
=
\theta(t)-\theta(a)-\theta'(a)(t-a)
-\frac12\theta''(a)(t-a)^2.
\]
Define the modulus of continuity
\[
\omega_2(r,a)=
\sup_{|t-a|\le r}|\theta''(t)-\theta''(a)|.
\]
By Condition~\ref{cond:rate:smoothness}, $\omega_2(r,a)\to0$ as $r\to0$.
Taylor's theorem in integral form implies
\[
|R_2(t,a)|
\le \frac12|t-a|^2\omega_2(|t-a|,a).
\]
Since $|t-a|\le h$ on the kernel support and the normalized
basis is bounded there,
\[
\| E_2\|=O\{h^2\omega_2(h,a)\}.
\]
In particular,
$e_1^T\{ D_h(a)\}^{-1} E_2=o(h^2)$.
Keeping the random inverse error separate, we obtain
\begin{align*}
e_1^T\{\hat{ D}_h(a)\}^{-1} E_2
&=
e_1^T\{ D_h(a)\}^{-1} E_2\\
&\quad+
O_p\!\left[
h^2\omega_2(h,a)
\{(nh)^{-1/2}+q_n(a)r_n(a)\}
\right].
\end{align*}
Since $\omega_2(h,a)\to0$, this completes the bound
\[
e_1^T\{\hat{ D}_h(a)\}^{-1} E_2
=
o(h^2)
+
O_p\!\left(
\frac{h^2}{\sqrt{nh}}+h^2q_n(a)r_n(a)
\right).
\]
Here $o(h^2)$ is the deterministic population projection;
the random inverse error is kept in the $O_p$ term.

For the additional smoothness remark only, suppose that
$\theta$ is four times continuously differentiable near $a$.
Taylor's theorem then gives
\[
R_2(t,a)
=
\frac16\theta'''(a)(t-a)^3
+\frac1{24}\theta^{(4)}(a)(t-a)^4
+o(|t-a|^4).
\]
Using symmetry of $K$ and
$p_A(a+hv)=p_A(a)+hv p_A'(a)+o(h)$ uniformly for $|v|\le1$,
the two components satisfy
\[
 E_2=
\begin{pmatrix}O(h^4)\\O(h^3)\end{pmatrix}.
\]
Because $e_1^T\{ D_h(a)\}^{-1}=(O(1),O(h))$, the population
projection is $O(h^4)=O(h^2b^2)$ under $h/b=O(1)$.
Thus, in this additional smoothness case,
\[
e_1^T\{\hat{ D}_h(a)\}^{-1} E_2
=
O(h^2b^2)
+
O_p\!\left(
\frac{h^2}{\sqrt{nh}}+h^2q_n(a)r_n(a)
\right).
\]
The deterministic smoothing remainder has therefore been
controlled within this term under either smoothness condition.

\medskip
\noindent\textbf{3. The derivative estimation error $ E_3$.}
We first analyze the local quadratic fit that defines
$\hat\theta''(a)$.
Set
\[
\beta_b^{(2)}(a)
=
\left(\theta(a),b\theta'(a),\frac12b^2\theta''(a)\right)^T
\]
and define its residual score by
\[
\Psi^{(2)}(O,a;b,\eta)
=
\Gamma_N^{(2)}(O,a;b,\eta)
-\Gamma_D^{(2)}(O,a;b,\eta)\beta_b^{(2)}(a).
\]
The pooled coefficient estimate is
$\hat\beta_b^{(2)}(a)=(\hat{ D}_b^{(2)}(a))^{-1}
\hat{ N}_b^{(2)}(a)$, so its estimating equation gives
\[
\hat\beta_b^{(2)}(a)-\beta_b^{(2)}(a)
=
(\hat{ D}_b^{(2)}(a))^{-1}
\frac1K\sum_{k=1}^K
\mathbb P_{n,k}[\Psi^{(2)}(O,a;b,\hat\eta_{-k})].
\]
To separate its errors, add and subtract the limiting score
inside each fold average. Center the resulting score difference
with respect to $\mathbb P$. This gives three terms, defined by
\begin{align*}
 S_{n,oracle}^{(2)}
&=\mathbb P_n[\Psi^{(2)}(O,a;b,\bar\eta)],\\
 R_{n,1}^{(2)}
&=\frac1K\sum_{k=1}^K(\mathbb P_{n,k}-\mathbb P)
[\Psi^{(2)}(O,a;b,\hat\eta_{-k})
-\Psi^{(2)}(O,a;b,\bar\eta)],\\
 R_{n,2}^{(2)}
&=\frac1K\sum_{k=1}^K\mathbb P
[\Psi^{(2)}(O,a;b,\hat\eta_{-k})
-\Psi^{(2)}(O,a;b,\bar\eta)].
\end{align*}
Thus the pooled local quadratic estimating equation gives
\[
\hat\theta''(a)-\theta''(a)
=
2b^{-2}e_3^T(\hat{ D}_b^{(2)}(a))^{-1}
\{ S_{n,oracle}^{(2)}
+ R_{n,1}^{(2)}+ R_{n,2}^{(2)}\}.
\]
The oracle term contains the smoothing bias and sampling noise.
Each fold's contribution to the second term has mean zero
conditional on that fold's training sample; the third term
averages the changes in the score expectation.
We bound these three terms separately.

The population quadratic design matrix is
\[
 D_b^{(2)}(a)
=
\int K_b(t-a) z_{t,b}^{(2)}(a)
\big( z_{t,b}^{(2)}(a)\big)^T p_A(t)\,dt.
\]
Writing $\nu_j=\int v^jK(v)\,dv$, we have
\[
 D_b^{(2)}(a)
\longrightarrow
p_A(a)
\begin{pmatrix}
1&0&\nu_2\\
0&\nu_2&0\\
\nu_2&0&\nu_4
\end{pmatrix}.
\]
This limit is positive definite because
$\nu_2>0$ and $\nu_4>\nu_2^2$.
Each design-matrix score entry has a bounded normalized
polynomial factor on $|t-a|\le b$.
Applying the design-matrix argument from
Theorem~\ref{thm: rate1} entry by entry gives
\[
\|\hat{ D}_b^{(2)}(a)- D_b^{(2)}(a)\|
=
O_p\{(nb)^{-1/2}+q_n(a)r_n(a)\},
\qquad
\|\hat{ D}_b^{(2)}(a)- D_b^{(2)}(a)\|=o_p(1).
\]
Here $nb\to\infty$ follows from $nh\to\infty$ and $h/b=O(1)$.
Consequently,
\[
\|(\hat{ D}_b^{(2)}(a))^{-1}\|=O_p(1),
\]
and the inverse difference has the same rate as the matrix
difference.

\emph{(i) The quadratic oracle term.}
By Condition~\ref{cond:rate:correctness},
\[
\mathbb E S_{n,oracle}^{(2)}
=
\int K_b(t-a) z_{t,b}^{(2)}(a)
R_2(t,a)p_A(t)\,dt.
\]
The Taylor bound already obtained for $R_2$ therefore implies
\[
\|\mathbb E S_{n,oracle}^{(2)}\|
=O\{b^2\omega_2(b,a)\}.
\]

To control the fluctuation, write the oracle score explicitly
as a point component and two tail components. Define
\[
\bar Q_b^{(2)}(t,X)
=
\frac1{\bar S_A(t\mid X)}
\int_t^T K_b(s-a) z_{s,b}^{(2)}(a)
\{\bar m(s)-( z_{s,b}^{(2)}(a))^T\beta_b^{(2)}(a)\}
\bar f_A(s\mid X)\,ds
\]
for $0\le t\le a+b$, and set it to zero for $t>a+b$.
Then
\[
\Psi^{(2)}(O,a;b,\bar\eta)
=
\bar T_1^{(2)}+\bar T_2^{(2)}-\bar T_3^{(2)},
\]
where
\begin{align*}
\bar T_1^{(2)}
&=
\frac{\Delta K_b(U-a)}{\bar S_C(U\mid X)}
 z_{U,b}^{(2)}(a)
\left\{
\frac{Y-\bar\mu(X,U)}{\bar f_A(U\mid X)}\bar p_A(U)
+\bar m(U)-( z_{U,b}^{(2)}(a))^T\beta_b^{(2)}(a)
\right\},\\
\bar T_2^{(2)}
&=
\frac{1-\Delta}{\bar S_C(U\mid X)}
\bar Q_b^{(2)}(U,X),\\
\bar T_3^{(2)}
&=
\int_0^U
\frac{\bar\lambda_C(t\mid X)}{\bar S_C(t\mid X)}
\bar Q_b^{(2)}(t,X)\,dt.
\end{align*}
For sufficiently large $n$,
$[a-(h\vee b),a+(h\vee b)]\subset I$ and
$a+(h\vee b)<\tau$.
The inverse-weight bounds and the conditional second-moment
bound for $Y$ imply
\[
\mathbb E\|\bar T_1^{(2)}\|^2
\le C_0\int K_b(t-a)^2p_A(t)\,dt
=O(b^{-1}).
\]
The other two terms are bounded.
Indeed, $\bar S_A(t\mid x)\ge c(T-\tau)$ on $[0,a+b]$,
the integral of $|K_b(s-a)|\,\| z_{s,b}^{(2)}(a)\|$
is bounded uniformly in $b$, and
\[
\int_0^{a+b}
\frac{\bar\lambda_C(t\mid x)}{\bar S_C(t\mid x)}\,dt
=
\bar S_C(a+b\mid x)^{-1}-1\le C_0.
\]
It follows that
\[
\mathbb E\|\Psi^{(2)}(O,a;b,\bar\eta)\|^2=O(b^{-1}),
\]
and hence
\[
\| S_{n,oracle}^{(2)}
-\mathbb E S_{n,oracle}^{(2)}\|
=O_p((nb)^{-1/2}).
\]

\emph{(ii) The quadratic empirical-process error.}
Fix a fold $k$ and condition on its training sample
$\mathcal I_{-k}$.
Define $\hat T_j^{(2)}$ by replacing all limiting functions
in $\bar T_j^{(2)}$ with their fitted versions from this fold,
including $\hat m_{-k}$ and $\hat p_{A,-k}$.
Set
\[
\Delta T_j^{(2)}=\hat T_j^{(2)}-\bar T_j^{(2)},\qquad
\Delta\Psi^{(2)}
=\Delta T_1^{(2)}+\Delta T_2^{(2)}-\Delta T_3^{(2)}.
\]
The fold contribution to $ R_{n,1}^{(2)}$ is
$(\mathbb P_{n,k}-\mathbb P)\Delta\Psi^{(2)}$,
which has conditional mean zero.

Let $\delta_{n,k}$ be the sum of the three fitted-to-limit
uniform errors specified in the theorem's setup, together with
\[
\sup_{|t-a|\le h\vee b}|\hat m_{-k}(t)-\bar m(t)|
+
\sup_{|t-a|\le h\vee b}|\hat p_{A,-k}(t)-\bar p_A(t)|.
\]
Conditions~\ref{cond:rate:marginal-definition} and~\ref{cond:db:marginal-error},
together with uniform convergence of the fitted nuisance functions,
imply $\delta_{n,k}=o_p(1)$.
For the point component, subtraction of the inverse-weighted
expressions gives
\[
\|\Delta T_1^{(2)}\|
\le C_0\delta_{n,k}\Delta |K_b(U-a)|(1+|Y|).
\]
Using the conditional second-moment bound for $Y$,
\[
\mathbb E(\|\Delta T_1^{(2)}\|^2\mid\mathcal I_{-k})
\le C_0\delta_{n,k}^2/b.
\]
Each fitted or limiting tail component is bounded by the same
survival and kernel-integral bounds used above. Therefore
\[
\mathbb E(\|\Delta\Psi^{(2)}\|^2\mid\mathcal I_{-k})
\le C_0\delta_{n,k}^2/b+C_0
=o_p(b^{-1}).
\]
Conditional Chebyshev's inequality, followed by averaging over
the fixed number of folds, yields
\[
\| R_{n,1}^{(2)}\|=o_p((nb)^{-1/2}).
\]

\emph{(iii) The quadratic nuisance drift.}
For this term, apply the integration-by-parts calculation in
Theorem~\ref{thm: rate1} with the weight
$K_b(t-a) z_{t,b}^{(2)}(a)$.
The replacement preserves the bounded kernel $L_1$ norm:
\[
\int |K_b(t-a)|\,\| z_{t,b}^{(2)}(a)\|\,dt
=
\int |K(v)|\sqrt{1+v^2+v^4}\,dv<\infty.
\]
It also preserves the tail support and survival lower bounds.
The censoring contribution is thus
$O_p[q_n(a)\{r_n(a)+s_n(a)\}]$:
the censoring ratio is controlled on $[0,\tau]$, the treatment
survival error by the density error on $[0,T]$, and the outcome
error on $|t-a|\le h\vee b$.
The full-data contribution is
$O_p\{r_n(a)s_n(a)+n^{-1/2}\}$, where the last term comes from
the two marginalization bounds in Condition~\ref{cond:db:marginal-error}.
Consequently,
\[
\| R_{n,2}^{(2)}\|
=O_p\{\chi_n(a)+n^{-1/2}\}.
\]

\emph{Combining the quadratic errors and bounding $ E_3$.}
The first term in the quadratic decomposition has deterministic
mean $\mathbb E S_{n,oracle}^{(2)}$.
Replacing the inverse only in this mean introduces an error of
order
\[
O_p\!\left[
\omega_2(b,a)\{(nb)^{-1/2}+q_n(a)r_n(a)\}
\right].
\]
The centered oracle term, the empirical-process error, and the
nuisance drift have already been bounded above.
After multiplying by $2b^{-2}e_3^T$, their combined contribution
therefore gives
\begin{align*}
\hat\theta''(a)-\theta''(a)
={}&
2b^{-2}e_3^T( D_b^{(2)}(a))^{-1}
\mathbb E S_{n,oracle}^{(2)}\\
&+
O_p\!\left(
\frac1{\sqrt{nb^5}}
+\frac{\chi_n(a)+n^{-1/2}}{b^2}
\right).
\end{align*}
The first line is deterministic and has order
$O\{\omega_2(b,a)\}=o(1)$.
Since $n^{-1/2}b^{-2}=\sqrt b\,(nb^5)^{-1/2}$, the
$n^{-1/2}$ term is absorbed by the displayed stochastic rate.

To control the multiplier $\hat{ B}_h(a)$ in
$ E_3$, observe that the bias score is already contained
in the quadratic design-matrix score:
\[
\Gamma_B(O,a;h,\eta)
=
\begin{pmatrix}
[\Gamma_D^{(2)}(O,a;h,\eta)]_{13}\\
[\Gamma_D^{(2)}(O,a;h,\eta)]_{23}
\end{pmatrix}.
\]
Here the bandwidth argument of the quadratic score is evaluated
at $h$. The two entries have polynomial factors $u^2$ and $u^3$,
including in their tail integrals.
Applying the design-matrix score bound entry by entry thus gives
\[
\|\hat{ B}_h(a)- B_h(a)\|
=
O_p\{(nh)^{-1/2}+q_n(a)r_n(a)\}.
\]
Conditions~\ref{cond:rate:limits} and~\ref{cond:rate:correctness} also make this difference $o_p(1)$.
In particular, $\hat{ B}_h(a)=O_p(1)$.
Combining this result with the inverse-matrix bound,
\begin{align*}
&e_1^T\{\hat{ D}_h(a)\}^{-1}\hat{ B}_h(a)
-e_1^T\{ D_h(a)\}^{-1} B_h(a)\\
&\quad=
e_1^T(\{\hat{ D}_h(a)\}^{-1}-\{ D_h(a)\}^{-1})
 B_h(a)
+
e_1^T\{\hat{ D}_h(a)\}^{-1}
(\hat{ B}_h(a)- B_h(a))\\
&\quad=
O_p\{(nh)^{-1/2}+q_n(a)r_n(a)\}.
\end{align*}

Symmetry of the kernel and continuity of \(p_A\) give
\[
 B_h(a)=\int K(u)(u^2,u^3)^\top p_A(a+hu)\,du
\longrightarrow p_A(a)(\nu_2,0)^\top.
\]
Together with the local linear design limit
\eqref{eq:rate_design_limit} and the preceding consistency
bounds, this proves
\[
e_1^\top\{\hat{ D}_h(a)\}^{-1}\hat{ B}_h(a)
\xrightarrow{p}\nu_2.
\]

Substituting the derivative expansion into the definition of
$ E_3$ now yields
\begin{align*}
e_1^T\{\hat{ D}_h(a)\}^{-1} E_3
={}&
\frac{h^2}{b^2}
\left(e_1^T\{ D_h(a)\}^{-1} B_h(a)\right)
e_3^T( D_b^{(2)}(a))^{-1}
\mathbb E S_{n,oracle}^{(2)}\\
&+
O_p\!\left(
\frac{h^2}{\sqrt{nb^5}}
+\chi_n(a)+\frac{h^2}{\sqrt{nh}}
\right).
\end{align*}
The first line is deterministic.
Its magnitude is $O\{h^2\omega_2(b,a)\}=o(h^2)$.
The stochastic bound uses $h^2/b^2=O(1)$; the error from the
empirical projection is multiplied by $h^2\omega_2(b,a)$
and is absorbed by the second line.
We have therefore completed the bound
\[
e_1^T\{\hat{ D}_h(a)\}^{-1} E_3
=
o(h^2)
+
O_p\!\left(
\frac{h^2}{\sqrt{nb^5}}
+\chi_n(a)+\frac{h^2}{\sqrt{nh}}
\right).
\]

For the additional fourth-order smoothness remark, use the
Taylor expansion of $R_2$ displayed in the analysis of
$ E_2$. Symmetry of $K$ and continuous differentiability
of $p_A$ imply
\[
\mathbb E S_{n,oracle}^{(2)}
=
\begin{pmatrix}O(b^4)\\O(b^3)\\O(b^4)\end{pmatrix}.
\]
The odd-moment entries of the population design matrix are
$O(b)$, so
\[
e_3^T( D_b^{(2)}(a))^{-1}
=(O(1),O(b),O(1)).
\]
It follows that
\[
e_3^T( D_b^{(2)}(a))^{-1}
\mathbb E S_{n,oracle}^{(2)}
=O(b^4).
\]
The deterministic first line of the preceding expansion for
$ E_3$ is therefore $O(h^2b^2)$.
Thus, under this additional smoothness,
\[
e_1^T\{\hat{ D}_h(a)\}^{-1} E_3
=
O(h^2b^2)
+
O_p\!\left(
\frac{h^2}{\sqrt{nb^5}}
+\chi_n(a)+\frac{h^2}{\sqrt{nh}}
\right).
\]
This completes the treatment of both the deterministic and
stochastic contributions of $ E_3$.

\medskip
\noindent\textbf{4. The empirical quadratic-term moment vector error $ E_4$.}
We now separate the quadratic-term moment vector error into its three sources.
Starting from the fold averages defining $\hat{ B}_h(a)$,
add and subtract the limiting score
$\Gamma_B(O,a;h,\bar\eta)$ in each fold.
The averages of this fixed limiting score combine into
$\mathbb P_n[\Gamma_B(O,a;h,\bar\eta)]$.
Subtracting its population expectation gives the oracle
fluctuation. We then center the fitted-minus-limiting score
difference within each fold, leaving its expectation as the
nuisance drift. Thus
\[
\hat{ B}_h(a)- B_h(a)
=R_{B,1}+R_{B,2}+R_{B,3},
\]
where
\begin{align*}
R_{B,1}
&=(\mathbb P_n-\mathbb P)
[\Gamma_B(O,a;h,\bar\eta)],\\
R_{B,2}
&=\frac1K\sum_{k=1}^K(\mathbb P_{n,k}-\mathbb P)
[\Gamma_B(O,a;h,\hat\eta_{-k})
-\Gamma_B(O,a;h,\bar\eta)],\\
R_{B,3}
&=\frac1K\sum_{k=1}^K\mathbb P
[\Gamma_B(O,a;h,\hat\eta_{-k})
-\Gamma_B(O,a;h,\bar\eta)].
\end{align*}
\emph{(i) The oracle fluctuation $R_{B,1}$.}
The bias basis is bounded on $|t-a|\le h$.
Its score has a point component with second moment $O(h^{-1})$
and two bounded tail components, by the same kernel and survival
bounds used for the quadratic score.
Therefore
\[
\mathbb E\|R_{B,1}\|^2
\le \frac1n
\mathbb E\|\Gamma_B(O,a;h,\bar\eta)\|^2
=O((nh)^{-1}),
\]
so $\|R_{B,1}\|=O_p((nh)^{-1/2})$.

\emph{(ii) The empirical-process error $R_{B,2}$.}
Condition on the training sample of a fold.
For the difference between its fitted and limiting bias scores,
the point component has conditional second moment at most
$C_0\delta_{n,k}^2/h$, and the tail components remain bounded.
Thus the full conditional second moment is
\[
C_0\delta_{n,k}^2/h+C_0=o_p(h^{-1}).
\]
The centered fold average has conditional mean zero.
Conditional Chebyshev's inequality consequently gives
$\|R_{B,2}\|=o_p((nh)^{-1/2})$.

\emph{(iii) The nuisance drift $R_{B,3}$.}

To calculate $R_{B,3}$, condition on one training sample
and suppress its fold index. Set
$g_h(s)=K_h(s-a) z_{s,h}^{B}(a)$, and define
\begin{align*}
 J_h(t,X)
={}&
\int_t^T g_h(s)
\{f_A(s\mid X)-\hat f_A(s\mid X)\}\,ds\\
&+
\left\{1-\frac{S_A(t\mid X)}{\hat S_A(t\mid X)}\right\}
\int_t^T g_h(s)\hat f_A(s\mid X)\,ds.
\end{align*}
The censoring expectation identity used in
Theorem~\ref{thm: rate1} gives the exact single-fold drift
\begin{align*}
&\mathbb P[
\Gamma_B(O,a;h,\hat\eta)
-\Gamma_B(O,a;h,\bar\eta)]\\
&\quad=
\mathbb E_X\int_0^{a+h}
 J_h(t,X)\,
d\left\{\frac{S_C(t\mid X)}{\hat S_C(t\mid X)}\right\}\\
&\quad=
-\mathbb E_X\int_0^{a+h}
\left\{\frac{S_C(t\mid X)}{\hat S_C(t\mid X)}-1\right\}
\,d J_h(t,X).
\end{align*}
The limiting score has the correct expectation by Condition~\ref{cond:rate:correctness}.
The boundary terms vanish because
$ J_h(a+h,X)=0$ and
$S_C(0\mid X)/\hat S_C(0\mid X)=1$.

The common support gives
\[
\sup_{t\in[0,\tau],\,x\in\mathcal X}
|\hat S_A(t\mid x)-S_A(t\mid x)|
\le
T\sup_{s\in[0,T],\,x\in\mathcal X}
|\hat f_A(s\mid x)-f_A(s\mid x)|
=O_p\{r_n(a)\}.
\]
Moreover,
\[
\partial_t\!\left(\frac{S_A}{\hat S_A}\right)
=
\frac{\hat f_A-f_A}{\hat S_A}
+
\frac{(S_A-\hat S_A)\hat f_A}{\hat S_A^2}.
\]
Using $\hat S_A\ge c(T-\tau)$, the bounded $L_1$ norm of $g_h$,
and the product inequality for total variation, we obtain
\[
\sup_{x\in\mathcal X}
\operatorname{TV}_{[0,a+h]}
\{ J_h(\cdot,x)\}
=O_p\{r_n(a)\}.
\]
Here $\operatorname{TV}$ denotes total variation in $t$,
using the Euclidean norm for the vector-valued function.
The norm of each fold's drift is therefore at most
\[
\sup_{t\in[0,\tau],\,x\in\mathcal X}
\left|\frac{S_C(t\mid x)}{\hat S_C(t\mid x)}-1\right|
\mathbb E_X[
\operatorname{TV}_{[0,a+h]}\{ J_h(\cdot,X)\}]
=
O_p\{q_n(a)r_n(a)\}.
\]
Averaging over the fixed number of folds gives
$\|R_{B,3}\|=O_p\{q_n(a)r_n(a)\}$.
Consequently,
\[
\|\hat{ B}_h(a)- B_h(a)\|
=
O_p\{(nh)^{-1/2}+q_n(a)r_n(a)\}.
\]
Conditions~\ref{cond:rate:limits} and~\ref{cond:rate:correctness} also imply that this difference is $o_p(1)$.
Since $\theta''(a)$ is fixed and
$\|\{\hat{ D}_h(a)\}^{-1}\|=O_p(1)$, this gives
\[
|e_1^T\{\hat{ D}_h(a)\}^{-1} E_4|
=
O_p\!\left(
\frac{h^2}{\sqrt{nh}}+h^2q_n(a)r_n(a)
\right).
\]
Thus $ E_4$ has no deterministic Taylor remainder.
The displayed stochastic bound is the same under the second-order
and fourth-order smoothness conditions.

\medskip
\noindent\textbf{Summary.}
Under the twice continuous differentiability in Condition~\ref{cond:rate:smoothness},
the four projected terms have now been bounded as follows:
\begin{align*}
|e_1^T\{\hat{ D}_h(a)\}^{-1} E_1|
&=O_p\!\left(\frac1{\sqrt{nh}}\right),\\
|e_1^T\{\hat{ D}_h(a)\}^{-1} E_2|
&=
O_p\!\left(
\frac{h^2}{\sqrt{nh}}+h^2q_n(a)r_n(a)
\right)+o(h^2),\\
|e_1^T\{\hat{ D}_h(a)\}^{-1} E_3|
&=
O_p\!\left(
\frac{h^2}{\sqrt{nb^5}}
+\chi_n(a)+\frac{h^2}{\sqrt{nh}}
\right)+o(h^2),\\
|e_1^T\{\hat{ D}_h(a)\}^{-1} E_4|
&=
O_p\!\left(
\frac{h^2}{\sqrt{nh}}+h^2q_n(a)r_n(a)
\right).
\end{align*}
The $o(h^2)$ terms here are bounds on the deterministic parts
identified within $ E_2$ and $ E_3$.
Under the additional fourth-order smoothness of the remark,
each of these two deterministic remainders is instead
$O(h^2b^2)$; the bounds for $ E_1$ and $ E_4$
are unchanged.

The two local linear nuisance terms satisfy
\[
|e_1^T\{\hat{ D}_h(a)\}^{-1}
( R_{n,1}+ R_{n,2})|
=
O_p\!\left(\frac1{\sqrt{nh}}+\chi_n(a)\right).
\]
Since $h^2/\sqrt{nh}=o((nh)^{-1/2})$ and
$h^2q_n(a)r_n(a)=O\{\chi_n(a)\}$, substituting the preceding
bounds into \eqref{eq:debiased_algebraic_trick} gives
\[
|\hat\theta_{\mathrm{DB}}(a)-\theta(a)|
=
O_p\!\left(
\frac1{\sqrt{nh}}
+\frac{h^2}{\sqrt{nb^5}}
+\chi_n(a)
\right)+o(h^2).
\]
Expanding $\chi_n(a)$ proves the theorem.

Under the additional fourth-order smoothness, the same
summation gives
\[
|\hat\theta_{\mathrm{DB}}(a)-\theta(a)|
=
O_p\!\left(
\frac1{\sqrt{nh}}
+\frac{h^2}{\sqrt{nb^5}}
+h^2b^2
+q_n(a)r_n(a)+q_n(a)s_n(a)+r_n(a)s_n(a)
\right),
\]
which proves the sharper rate in the remark.
\end{proof}

\subsection{Asymptotic Normality of the Debiased Estimator}

\begin{proof}[Proof of Theorem~\ref{thm:debiased_asymptotic_normality}]
Use the notation from the proof of
Theorem~\ref{thm:debiased_convergence}, and write
\[
\chi_n(a)=q_n(a)r_n(a)+q_n(a)s_n(a)+r_n(a)s_n(a),
\qquad
\kappa_n=\frac hb.
\]
Here $\kappa_n>0$ is the actual bandwidth ratio, whereas
$\kappa$ denotes its limit.
In particular, $\kappa_n=O(1)$, and $nb\to\infty$ follows
from $nh\to\infty$.

We first isolate the centered oracle fluctuations.
The exact decomposition from the preceding proof is
\begin{equation}
\label{eq:clt_expansion}
\begin{aligned}
\sqrt{nh}\{\hat\theta_{\mathrm{DB}}(a)-\theta(a)\}
={}&
\sqrt{nh}\,e_1^T\{\hat{ D}_h(a)\}^{-1}
\{ E_1+ E_2- E_3- E_4\\
&\hspace{42mm}+ R_{n,1}+ R_{n,2}\}.
\end{aligned}
\end{equation}
Condition~\ref{cond:clt:products}, the local linear remainder bounds, and
$\sqrt{nh}\,n^{-1/2}=\sqrt h\to0$ give
\[
\sqrt{nh}\,\| R_{n,1}\|=o_p(1),
\qquad
\sqrt{nh}\,\| R_{n,2}\|=o_p(1).
\]
The preceding proof also gives
\[
\| E_2\|=o(h^2),
\qquad
\| E_4\|
=
O_p\{h^2((nh)^{-1/2}+q_n(a)r_n(a))\}.
\]
Thus, by Condition~\ref{cond:clt:products} and the bandwidth requirement~\eqref{eq:clt_bandwidth},
\[
\sqrt{nh}\,\| E_2\|
=
\sqrt{nh^5}\,o(1)=o(1),
\qquad
\sqrt{nh}\,\| E_4\|=o_p(1).
\]
These conclusions use only twice continuous differentiability.

The local linear design matrix has the positive definite limit
\[
 D(a)
=
p_A(a)
\begin{pmatrix}
1&0\\
0&\nu_2
\end{pmatrix}.
\]
Since $\hat{ D}_h(a)\to_p D(a)$,
\[
\{\hat{ D}_h(a)\}^{-1}=\{ D(a)\}^{-1}+o_p(1),
\qquad
e_1^T\{ D(a)\}^{-1}=\bigl(p_A(a)^{-1},0\bigr),
\]
where the matrix remainder is measured in operator norm.
Furthermore,
$\sqrt{nh}\,\| E_1\|=O_p(1)$.
It follows that
\[
\sqrt{nh}\,e_1^T\{\hat{ D}_h(a)\}^{-1} E_1
=
\sqrt{nh}\,e_1^T\{ D(a)\}^{-1} E_1+o_p(1).
\]

We next extract the centered fluctuation in the derivative
estimate. Recall the local quadratic coefficient
\[
\beta_b^{(2)}(a)
=
\left(\theta(a),b\theta'(a),\frac12b^2\theta''(a)\right)^T
\]
and its oracle residual score
\[
\Psi^{(2)}(O,a;b,\bar\eta)
=
\Gamma_N^{(2)}(O,a;b,\bar\eta)
-
\Gamma_D^{(2)}(O,a;b,\bar\eta)\beta_b^{(2)}(a).
\]
Write
\[
 S_{n,oracle}^{(2)}
=
\mathbb P_n\Psi^{(2)}(O,a;b,\bar\eta),
\qquad
 E_1^{(2)}
=
(\mathbb P_n-\mathbb P)\Psi^{(2)}(O,a;b,\bar\eta).
\]
The exact local quadratic expansion established in the
preceding proof is
\[
\begin{aligned}
\hat\theta''(a)-\theta''(a)
=
2b^{-2}e_3^T
\bigl(\hat{ D}_b^{(2)}(a)\bigr)^{-1}
\bigl\{
 E_1^{(2)}
+\mathbb E S_{n,oracle}^{(2)}
+ R_{n,1}^{(2)}
+ R_{n,2}^{(2)}
\bigr\},
\end{aligned}
\]
where $e_3=(0,0,1)^T$.
The two quadratic nuisance terms are the foldwise centered
score difference and its expectation, respectively, defined
in the preceding proof. Their bounds are
\[
\| R_{n,1}^{(2)}\|=o_p((nb)^{-1/2}),
\qquad
\| R_{n,2}^{(2)}\|
=O_p\{\chi_n(a)+n^{-1/2}\}.
\]
The oracle terms satisfy
\[
\| E_1^{(2)}\|=O_p((nb)^{-1/2}),
\qquad
\|\mathbb E S_{n,oracle}^{(2)}\|=o(b^2).
\]

The local quadratic design matrix converges to
\[
 D^{(2)}(a)
=
p_A(a)
\begin{pmatrix}
1&0&\nu_2\\
0&\nu_2&0\\
\nu_2&0&\nu_4
\end{pmatrix},
\]
whose inverse satisfies
\[
e_3^T\bigl( D^{(2)}(a)\bigr)^{-1}
=
\frac{1}{p_A(a)(\nu_4-\nu_2^2)}
(-\nu_2,0,1).
\]
Consequently,
\[
\bigl(\hat{ D}_b^{(2)}(a)\bigr)^{-1}
=
\bigl( D^{(2)}(a)\bigr)^{-1}+o_p(1).
\]
Also,
$\hat{ B}_h(a)\to_p p_A(a)(\nu_2,0)^T$, so
\[
e_1^T\{\hat{ D}_h(a)\}^{-1}\hat{ B}_h(a)
=
\nu_2+o_p(1).
\]

To substitute these limits into $ E_3$, note that the
coefficient multiplying the quadratic score is $h^2/b^2$.
The scaled oracle mean and nuisance terms therefore obey
\begin{align*}
\sqrt{nh}\,\frac{h^2}{b^2}
\|\mathbb E S_{n,oracle}^{(2)}\|
&=
\sqrt{nh^5}\,o(1)=o(1),\\
\sqrt{nh}\,\frac{h^2}{b^2}
\| R_{n,1}^{(2)}\|
&=
\kappa_n^{5/2}o_p(1)=o_p(1),\\
\sqrt{nh}\,\frac{h^2}{b^2}
\| R_{n,2}^{(2)}\|
&=
\kappa_n^2 O_p\{\sqrt{nh}\,\chi_n(a)+\sqrt h\}
=o_p(1).
\end{align*}
Moreover,
\[
\sqrt{nh}\,\frac{h^2}{b^2}
\| E_1^{(2)}\|
=
\kappa_n^{5/2}O_p(1)=O_p(1).
\]
Thus the additive $o_p(1)$ matrix and quadratic-term moment vector errors
contribute only $o_p(1)$ after scaling.
Using the definition
$ E_3=\frac12h^2\{\hat\theta''(a)-\theta''(a)\}
\hat{ B}_h(a)$, we obtain
\[
\sqrt{nh}\,e_1^T\{\hat{ D}_h(a)\}^{-1} E_3
=
\sqrt{nh}\,\frac{h^2}{b^2}\nu_2
e_3^T\bigl( D^{(2)}(a)\bigr)^{-1}
 E_1^{(2)}
+o_p(1).
\]
Combining these bounds in \eqref{eq:clt_expansion} yields
\begin{equation}
\label{eq:debiased_oracle_dominance}
\begin{aligned}
\sqrt{nh}\{\hat\theta_{\mathrm{DB}}(a)-\theta(a)\}
=
\sqrt{nh}\left\{
e_1^T\{ D(a)\}^{-1} E_1
-
\frac{h^2}{b^2}\nu_2
e_3^T\bigl( D^{(2)}(a)\bigr)^{-1}
 E_1^{(2)}
\right\}
+o_p(1).
\end{aligned}
\end{equation}

We now express the two oracle fluctuations through their
uncensored contributions. Define
\[
W(O)
=
\frac{\Delta}{\bar S_C(U\mid X)}
\{\bar\xi(X,U,Y)-\theta(U)\},
\]
with the expression set to zero when $\Delta=0$.
Under Conditions~\ref{cond:rate:correctness} and~\ref{cond:rate:marginal-definition},
\[
\mathbb E\{W(O)\mid A=t\}=0.
\]
To see this, first suppose that $\bar\mu=\mu$.
Then $\bar m=\theta$, and
$\bar\xi(X,t,Y)-\theta(t)$ is a bounded multiple of
$Y-\mu(X,t)$.
Its conditional mean given $(X,A=t)$ is zero.
Conditional independent censoring preserves this property
after multiplication by $\Delta/\bar S_C(A\mid X)$.
If $\bar\mu\ne\mu$, Condition~\ref{cond:rate:correctness} requires
$\bar S_C=S_C$ and $\bar f_A=f_A$.
In this case $\bar p_A=p_A$, the censoring weight cancels
in conditional expectation, and
\[
\begin{aligned}
\mathbb E\{\bar\xi(X,t,Y)\mid A=t\}
&=
\mathbb E_X\{\mu(X,t)-\bar\mu(X,t)\}+\bar m(t)\\
&=\theta(t).
\end{aligned}
\]
This proves the claimed centering in all three cases of
Condition~\ref{cond:rate:correctness}.

For the local linear score, add and subtract $\theta(U)$
inside its uncensored term. The resulting Taylor remainder is
\[
R_h(t,a)
=
\theta(t)-\theta(a)-\theta'(a)(t-a)
=
O(h^2),
\qquad |t-a|\le h.
\]
The score's two augmentation terms are uniformly bounded,
as shown in the preceding rate proof.
Their centered empirical averages are consequently
$O_p(n^{-1/2})$.
The centered uncensored Taylor contribution has second
moment bounded by
\[
C_0\int K_h(t-a)^2 R_h(t,a)^2p_A(t)\,dt
=
O(h^3),
\]
and hence its empirical fluctuation is
$O_p(h^2/\sqrt{nh})$.
Since the entire oracle score is centered in $ E_1$,
these facts imply
\[
e_1^T\{ D(a)\}^{-1} E_1
=
\frac{1}{np_A(a)}
\sum_{i=1}^n K_h(U_i-a)W(O_i)
+
o_p((nh)^{-1/2}).
\]
In particular, only the centered augmentation fluctuations
have been discarded here.

For the quadratic score, the corresponding Taylor remainder is
\[
R_b^{(2)}(t,a)
=
\theta(t)-\theta(a)-\theta'(a)(t-a)
-\frac12\theta''(a)(t-a)^2
=
o(b^2)
\]
uniformly over $|t-a|\le b$.
The projected polynomial
$((t-a)/b)^2-\nu_2$ is bounded on this window.
Thus the centered uncensored Taylor contribution to the
projected quadratic score is $o_p(b^2/\sqrt{nb})$.
The centered augmentation contribution is again
$O_p(n^{-1/2})$.
After multiplication by $\sqrt{nh}\,h^2/b^2$, these bounds give
\[
\sqrt{nh}\,\frac{h^2}{b^2}
o_p\!\left(\frac{b^2}{\sqrt{nb}}\right)
=
h^2\sqrt{\kappa_n}\,o_p(1)=o_p(1),
\]
and
\[
\sqrt{nh}\,\frac{h^2}{b^2}O_p(n^{-1/2})
=
\kappa_n^2 O_p(\sqrt h)=o_p(1).
\]
The inverse-matrix row computed above therefore yields
\[
\begin{aligned}
&\sqrt{nh}\,\frac{h^2}{b^2}\nu_2
e_3^T\bigl( D^{(2)}(a)\bigr)^{-1} E_1^{(2)}
\\
&\quad=
\frac{\sqrt{nh}}{np_A(a)}
\frac{h^2}{b^2}\frac{\nu_2}{\nu_4-\nu_2^2}
\sum_{i=1}^n
K_b(U_i-a)
\left\{\left(\frac{U_i-a}{b}\right)^2-\nu_2\right\}
W(O_i)
+o_p(1).
\end{aligned}
\]
This argument does not require consistency of
$\hat\theta''(a)$ itself.

Subtracting the last two projected processes combines their
sampling variation, including their covariance.
For convenience, put
\[
L(u)=\frac{u^2-\nu_2}{\nu_4-\nu_2^2}K(u).
\]
For the actual ratio $\kappa_n=h/b$, define
\[
K_{\mathrm{DB},n}(v)
=
K(v)-\nu_2\kappa_n^3 L(\kappa_n v).
\]
Indeed, with $v=(t-a)/h$,
\[
\begin{aligned}
&K_h(t-a)
-\frac{h^2}{b^2}
\frac{\nu_2}{\nu_4-\nu_2^2}
K_b(t-a)\left\{\left(\frac{t-a}{b}\right)^2-\nu_2\right\}\\
&\qquad=
\frac1h
\left\{
K(v)
-\frac{\nu_2\kappa_n^3}{\nu_4-\nu_2^2}
(\kappa_n^2v^2-\nu_2)K(\kappa_n v)
\right\}
=
\frac1h K_{\mathrm{DB},n}(v).
\end{aligned}
\]
This is a finite-sample identity.
Consequently, \eqref{eq:debiased_oracle_dominance} becomes
\[
\sqrt{nh}\{\hat\theta_{\mathrm{DB}}(a)-\theta(a)\}
=
\frac1{\sqrt{nh}}
\sum_{i=1}^n
\frac{W(O_i)}{p_A(a)}
K_{\mathrm{DB},n}\!\left(\frac{U_i-a}{h}\right)
+o_p(1).
\]

To pass from the actual equivalent kernel to its limit, first
consider $\kappa>0$.
Continuity and compact support of $K$, together with
$\kappa_n\to\kappa$, give
$\|K_{\mathrm{DB},n}-K_{\mathrm{DB}}\|_2\to0$.
When $\kappa=0$, a change of variables instead gives
\[
\|K_{\mathrm{DB},n}-K\|_2^2
=
\nu_2^2\kappa_n^5\int L(u)^2\,du
\longrightarrow0.
\]
Thus the $L_2$ convergence holds in both cases, with the
definition in \eqref{eq:equivalent_kernel}, and
\[
\int\{K_{\mathrm{DB},n}(u)\}^2\,du
\longrightarrow V_{K,\kappa}.
\]
Although the support of $K_{\mathrm{DB},n}$ need not be uniformly
bounded when $\kappa=0$, its contribution at treatment level
$t$ is supported on $|t-a|\le h\vee b$.
This interval shrinks to $a$.

It remains to evaluate the variance and verify the
triangular-array central limit theorem.
Set
\[
\zeta_{ni}
=
\frac{W(O_i)}{p_A(a)\sqrt h}
K_{\mathrm{DB},n}\!\left(\frac{U_i-a}{h}\right).
\]
The conditional centering of $W$ gives
$\mathbb E\zeta_{ni}=0$.
On $\{\Delta=1\}$, $U=A$, and
$\mathbb E(\Delta\mid X,A,Y)=S_C(A\mid X)$.
Therefore, with $v(t)$ defined in Section~\ref{sec:mr_asymptotics},
\[
\begin{aligned}
\mathbb E\zeta_{ni}^2
&=
\frac{1}{p_A(a)^2h}
\int
K_{\mathrm{DB},n}\!\left(\frac{t-a}{h}\right)^2
p_A(t)v(t)\,dt\\
&=
\frac{1}{p_A(a)^2}
\int\{K_{\mathrm{DB},n}(u)\}^2
p_A(a+hu)v(a+hu)\,du.
\end{aligned}
\]
By continuity of $p_A(t)v(t)$ at $a$,
\[
\begin{aligned}
&\left|
\mathbb E\zeta_{ni}^2
-
\frac{v(a)}{p_A(a)}
\int\{K_{\mathrm{DB},n}(u)\}^2\,du
\right|\\
&\quad\le
\frac{1}{p_A(a)^2}
\sup_{|t-a|\le h\vee b}
|p_A(t)v(t)-p_A(a)v(a)|
\int\{K_{\mathrm{DB},n}(u)\}^2\,du
=o(1).
\end{aligned}
\]
Hence
\[
\mathbb E\zeta_{ni}^2
\longrightarrow
\frac{v(a)}{p_A(a)}V_{K,\kappa}
=
\sigma_{\mathrm{DB}}^2(a).
\]

For the Lyapunov condition, let $\delta>0$ be as in
Condition~\ref{cond:clt:moments}.
The bounded nuisance coefficients and inverse weights imply
\[
|W(O)|\le C_0\Delta(1+|Y|)
\]
on the active kernel window, which is eventually contained in
$I$.
Moreover, since $\kappa_n=O(1)$,
\[
\begin{aligned}
\int|K_{\mathrm{DB},n}(u)|^{2+\delta}\,du
&\le
C_0\int|K(u)|^{2+\delta}\,du\\
&\quad+
C_0\kappa_n^{5+3\delta}
\int|L(u)|^{2+\delta}\,du
=O(1).
\end{aligned}
\]
Condition~\ref{cond:clt:moments} and a change of variables consequently give
\[
\mathbb E|\zeta_{ni}|^{2+\delta}
\le
C_0 h^{-\delta/2}
\int|K_{\mathrm{DB},n}(u)|^{2+\delta}\,du
=
O(h^{-\delta/2}).
\]
It follows that
\[
n^{-\delta/2}\mathbb E|\zeta_{ni}|^{2+\delta}
=
O((nh)^{-\delta/2})\longrightarrow0.
\]
Since $\sigma_{\mathrm{DB}}^2(a)>0$, the Lyapunov condition holds. The triangular-array central limit theorem \citep[Chapter~2]{vandervaart1998} gives
\[
\frac1{\sqrt n}\sum_{i=1}^n\zeta_{ni}
\xrightarrow{d}
\mathcal N\bigl(0,\sigma_{\mathrm{DB}}^2(a)\bigr).
\]

For completeness, the limiting kernel can also be used in
the leading stochastic representation.
The difference obtained by replacing $K_{\mathrm{DB},n}$ with
$K_{\mathrm{DB}}$ has mean zero and variance bounded by
\[
C_0\|K_{\mathrm{DB},n}-K_{\mathrm{DB}}\|_2^2=o(1),
\]
because $p_A(t)v(t)$ is bounded on the combined shrinking
treatment windows.
We therefore also have
\begin{equation}
\label{eq:final_clt_form}
\begin{aligned}
&\sqrt{nh}\{\hat\theta_{\mathrm{DB}}(a)-\theta(a)\}\\
&\quad=
\frac1{\sqrt{nh}}\sum_{i=1}^n
\frac{\Delta_i}{p_A(a)\bar S_C(U_i\mid X_i)}
K_{\mathrm{DB}}\!\left(\frac{U_i-a}{h}\right)
\{\bar\xi(X_i,U_i,Y_i)-\theta(U_i)\}
+o_p(1).
\end{aligned}
\end{equation}
Applying Slutsky's theorem gives
\[
\sqrt{nh}\{\hat\theta_{\mathrm{DB}}(a)-\theta(a)\}
\xrightarrow{d}
\mathcal N\bigl(0,\sigma_{\mathrm{DB}}^2(a)\bigr),
\]
with the variance in \eqref{eq:debiased_variance}.
\end{proof}

\subsection{Variance Estimation and Pointwise Confidence Intervals}
\label{app:db_variance}

We describe the linearized contributions used in the simulations.
The conditional nuisance fits, bandwidths, and fold assignments
are held fixed in this construction. The empirical averages
defining the marginal functions are allowed to vary with the
observations.

Let \(\mathcal I_k\) and \(\mathcal I_{-k}\) denote evaluation
fold \(k\) and its training sample. For positive observation
weights \( w=(w_1,\ldots,w_n)^\top\), define
\[
\hat m_{-k}(t; w)=
\frac{\sum_{j\in\mathcal I_{-k}}w_j\hat\mu_{-k}(X_j,t)}
{\sum_{j\in\mathcal I_{-k}}w_j},\qquad
\hat p_{A,-k}(t; w)=
\frac{\sum_{j\in\mathcal I_{-k}}w_j\hat f_{A,-k}(t\mid X_j)}
{\sum_{j\in\mathcal I_{-k}}w_j}.
\]
At \(w=(1,\ldots,1)^\top\), these are the training-sample
marginal averages used in the implementation.
Evaluate the estimating equations with these marginal
functions and replace each pooled average
\(n^{-1}\sum_i F_i\) by \(n^{-1}\sum_i w_iF_i\).
Apply this replacement to both local systems and to the
empirical quadratic-term moment vector. Let \(\hat\theta_{\mathrm{DB}}(a; w)\)
denote the resulting estimator, using the same bias-correction
formula as in \eqref{eq:debiased_estimator}.
The observation-level linearized contribution is
\begin{equation}
\label{eq:db_full_contribution}
\hat\varphi_{\mathrm{DB},i}(a)=
n\left.\frac{\partial\hat\theta_{\mathrm{DB}}(a; w)}
{\partial w_i}\right|_{w=(1,\ldots,1)^\top}.
\end{equation}
The derivative exists on the event that both fitted design
matrices are invertible, whose probability tends to one.

Formula~\eqref{eq:db_full_contribution} differentiates both
the response vectors and the design matrices. Differentiating
the bias-correction product also retains the contribution
from the empirical bias projection. In particular, the
training-sample marginalization terms follow from
\[
\left.\frac{\partial\hat m_{-k}(t; w)}{\partial w_i}
\right|_{w=(1,\ldots,1)^\top}
=\frac{I(i\in\mathcal I_{-k})}{|\mathcal I_{-k}|}
\{\hat\mu_{-k}(X_i,t)-\hat m_{-k}(t)\},
\]
\[
\left.\frac{\partial\hat p_{A,-k}(t; w)}{\partial w_i}
\right|_{w=(1,\ldots,1)^\top}
=\frac{I(i\in\mathcal I_{-k})}{|\mathcal I_{-k}|}
\{\hat f_{A,-k}(t\mid X_i)-\hat p_{A,-k}(t)\}.
\]
These derivatives are included for every fold in which
observation \(i\) belongs to the training sample.
The conditional nuisance models themselves are not refitted
when taking the derivative. Write
\(\bar\varphi_{\mathrm{DB}}(a)=n^{-1}\sum_{i=1}^n\hat\varphi_{\mathrm{DB},i}(a)\).
The variance estimator is
\begin{equation}
\label{eq:db_variance_estimator}
\hat\sigma_{\mathrm{DB}}^2(a)=\frac hn\sum_{i=1}^n
\{\hat\varphi_{\mathrm{DB},i}(a)-\bar\varphi_{\mathrm{DB}}(a)\}^2.
\end{equation}
The standard error calculated from \eqref{eq:db_variance_estimator} is
\[
\frac{\hat\sigma_{\mathrm{DB}}(a)}{\sqrt{nh}}
=\frac1n\left[\sum_{i=1}^n
\{\hat\varphi_{\mathrm{DB},i}(a)-\bar\varphi_{\mathrm{DB}}(a)\}^2\right]^{1/2}.
\]
For fixed \(\alpha\in(0,1)\), the corresponding pointwise
confidence interval is
\[
\left[\hat\theta_{\mathrm{DB}}(a)-z_{1-\alpha/2}\frac{\hat\sigma_{\mathrm{DB}}(a)}{\sqrt{nh}},\quad
\hat\theta_{\mathrm{DB}}(a)+z_{1-\alpha/2}\frac{\hat\sigma_{\mathrm{DB}}(a)}{\sqrt{nh}}\right],
\]
where \(z_{1-\alpha/2}\) is the corresponding standard normal
quantile. We establish variance consistency and asymptotic
coverage under the conditions of
Theorem~\ref{thm:debiased_asymptotic_normality} below.

\begin{proof}[Proof of consistency of the variance estimator]
Use the notation of the preceding rate and normality proofs.
Write \(\kappa_n=h/b\) and let \(K_{\mathrm{DB},n}\) denote the
equivalent kernel with this actual bandwidth ratio. Define
the oracle leading contribution
\[
\varphi_{\mathrm{DB},i}^{\,0}(a)=
\frac{W(O_i)}{h p_A(a)}
K_{\mathrm{DB},n}\!\left(\frac{U_i-a}{h}\right),
\]
where \(W(O)\) is defined in the proof of
Theorem~\ref{thm:debiased_asymptotic_normality}.
We first establish the empirical squared-error comparison
\begin{equation}
\label{eq:db_contribution_comparison}
\frac hn\sum_{i=1}^n
\{\hat\varphi_{\mathrm{DB},i}(a)-\varphi_{\mathrm{DB},i}^{\,0}(a)\}^2=o_p(1).
\end{equation}

For \(\ell\in\{h,b\}\), the conditional moment bounds and
the bounded treatment density imply
\[
\frac\ell n\sum_{i=1}^n
\Delta_i K_\ell(U_i-a)^2(1+Y_i^2)=O_p(1).
\]
The matrix limits and the fitted-to-limit uniform convergence
therefore allow the point contributions in
\eqref{eq:db_full_contribution} to be replaced by their oracle
counterparts in the scaled empirical squared norm.
The fitted local coefficients also satisfy
\(\hat\beta_h(a)-\beta_h(a)=o_p(1)\) and
\(\hat\beta_b^{(2)}(a)-\beta_b^{(2)}(a)=o_p(1)\), by the bounds
in the rate proof. Thus their fitted polynomial values may
be replaced by \(\theta(t)\) on the corresponding shrinking
windows in the same norm. For the quadratic score, the
scaling satisfies
\[
h\left(\frac hb\right)^4
\frac1n\sum_{i=1}^n\Delta_i K_b(U_i-a)^2(1+Y_i^2)
=\left(\frac hb\right)^5 O_p(1)=O_p(1).
\]
This also covers \(\kappa=0\).

The fitted augmentation terms are uniformly bounded, as in
the rate proof, and hence contribute \(O_p(h)=o_p(1)\) to
the scaled empirical squared norm. The additional
marginalization derivatives have bounded directions because
\(\hat\mu\), \(\hat f_A\), and their empirical averages
are bounded on the local windows. For either bandwidth
\(\ell\), the norm of each such contribution is bounded,
up to a fixed constant, by
\[
\max_k\left\{1+\frac1{|\mathcal I_{-k}|}
\sum_{j\in\mathcal I_k}
\Delta_j|K_\ell(U_j-a)|(1+|Y_j|)\right\}=O_p(1).
\]
Here the final bound follows from the bounded kernel
\(L_1\) norm, the moment conditions, and the fixed number
of folds. These contributions likewise vanish in the scaled
empirical squared norm.

The derivative of the empirical bias projection is multiplied
by \((h/b)^2e_3^\top\hat\beta_b^{(2)}(a)\).
The quadratic coefficient bounds give
\(e_3^\top\hat\beta_b^{(2)}(a)=o_p(1)\).
The projection derivative has scaled empirical squared norm
\(O_p(1)\), by the design and bias-score bounds. Its
contribution therefore vanishes. The remaining local linear
and quadratic point contributions combine into
\(\varphi_{\mathrm{DB},i}^{\,0}(a)\), by the equivalent-kernel
identity in the normality proof. This establishes
\eqref{eq:db_contribution_comparison}.

We next verify consistency of the oracle second moment.
Let \(\rho=1+\min(\delta/2,1)\in(1,2]\), where
\(\delta\) is from Condition~\ref{cond:clt:moments}, and set
\(G_{ni}=h\{\varphi_{\mathrm{DB},i}^{\,0}(a)\}^2\).
The variance calculation in the normality proof gives
\(E(G_{ni})\to\sigma_{\mathrm{DB}}^2(a)\). The same inverse-weight
and kernel bounds give
\[
E|G_{ni}|^\rho
\le C_0 h^{1-\rho}\int|K_{\mathrm{DB},n}(u)|^{2\rho}\,du
=O(h^{1-\rho}).
\]
The integral is bounded since \(\kappa_n=O(1)\).
By the moment inequality for independent centered sums
\citep{vonbahr1965},
\[
E\left|\frac1n\sum_{i=1}^n\{G_{ni}-E(G_{ni})\}\right|^\rho
\le C_0 n^{1-\rho}E|G_{n1}|^\rho
=O\{(nh)^{1-\rho}\}\longrightarrow0.
\]
Consequently,
\(h\mathbb P_n\{(\varphi_{\mathrm{DB}}^{\,0})^2\}
\xrightarrow{p}\sigma_{\mathrm{DB}}^2(a)\).
Conditional centering of \(W\) also gives
\(E\{\varphi_{\mathrm{DB},i}^{\,0}(a)\}=0\), so
\(\sqrt h\,\mathbb P_n\varphi_{\mathrm{DB}}^{\,0}=O_p(n^{-1/2})\).
By \eqref{eq:db_contribution_comparison} and
Cauchy-Schwarz, replacing the oracle contributions by
\(\hat\varphi_{\mathrm{DB},i}(a)\), and then centering by
\(\bar\varphi_{\mathrm{DB}}(a)\), changes the scaled second moment
by \(o_p(1)\). Thus
\(\hat\sigma_{\mathrm{DB}}^2(a)\xrightarrow{p}\sigma_{\mathrm{DB}}^2(a)>0\).
The claimed pointwise coverage follows from
Theorem~\ref{thm:debiased_asymptotic_normality} and
Slutsky's theorem.
\end{proof}

\subsection{Joint Asymptotic Normality}

\begin{proof}[Proof of Theorem~\ref{thm:joint_normality}]
The case $m=1$ follows from
Theorem~\ref{thm:debiased_asymptotic_normality}.
Henceforth suppose that $m\ge2$.

For each $a_j\in\mathcal A^*$, define the single-observation
leading term
\[
\zeta_{ni,j}
=
\frac{\Delta_i}{p_A(a_j)\bar S_C(U_i\mid X_i)\sqrt h}
K_{\mathrm{DB}}\!\left(\frac{U_i-a_j}{h}\right)
\{\bar\xi(X_i,U_i,Y_i)-\theta(U_i)\},
\]
with the expression set to zero when $\Delta_i=0$.
The representation \eqref{eq:final_clt_form} gives
\[
\sqrt{nh}\{\hat\theta_{\mathrm{DB}}(a_j)-\theta(a_j)\}
=
\frac1{\sqrt n}\sum_{i=1}^n\zeta_{ni,j}
+\varepsilon_{n,j},
\qquad
\varepsilon_{n,j}=o_p(1).
\]
Since $m$ is fixed, for every $\epsilon>0$,
\[
P\!\left(
\max_{1\le j\le m}|\varepsilon_{n,j}|>\epsilon
\right)
\le
\sum_{j=1}^m
P(|\varepsilon_{n,j}|>\epsilon)
\longrightarrow0.
\]
Thus the remainders are jointly negligible.
This argument does not require independence across evaluation
points or across their fitted nuisance estimators.

Fix a nonzero vector $c=(c_1,\ldots,c_m)^T\in\mathbb R^m$,
and put
\[
Z_{ni}(c)=\sum_{j=1}^m c_j\zeta_{ni,j}.
\]
Then
\[
\sum_{j=1}^m c_j\sqrt{nh}
\{\hat\theta_{\mathrm{DB}}(a_j)-\theta(a_j)\}
=
\frac1{\sqrt n}\sum_{i=1}^n Z_{ni}(c)+o_p(1).
\]
The limiting nuisance functions are fixed.
Consequently, for each $n$, the variables $Z_{ni}(c)$ are
independent and identically distributed across observations
$i$, even though the original estimators share training data.
The conditional centering established in the preceding
pointwise proof gives
\[
\mathbb E\zeta_{ni,j}=0,
\qquad
\mathbb E Z_{ni}(c)=0.
\]

We next evaluate the cross-point covariances.
Since the evaluation points are fixed and distinct,
\[
d_{\mathcal A}
=
\min_{j\ne\ell}|a_j-a_\ell|>0.
\]
The limiting equivalent kernel $K_{\mathrm{DB}}$ has compact support
for every fixed $\kappa\in[0,\infty)$.
This includes $\kappa=0$, for which $K_{\mathrm{DB}}=K$.
Choose a finite constant $C_K$ such that
$\operatorname{supp}(K_{\mathrm{DB}})\subseteq[-C_K,C_K]$.
Because $h\to0$, for all sufficiently large $n$,
\[
2C_Kh<d_{\mathcal A}.
\]
The kernel windows around distinct points are then disjoint:
\[
K_{\mathrm{DB}}\!\left(\frac{U_i-a_j}{h}\right)
K_{\mathrm{DB}}\!\left(\frac{U_i-a_\ell}{h}\right)
=0,
\qquad j\ne\ell.
\]
It follows that
\[
\mathbb E(\zeta_{ni,j}\zeta_{ni,\ell})=0,
\qquad
\operatorname{Cov}(\zeta_{ni,j},\zeta_{ni,\ell})=0
\]
for all sufficiently large $n$.
Here we use the fixed limiting-kernel representation already
justified in \eqref{eq:final_clt_form}; its validity includes
the case $\kappa=0$.

The same second-moment calculation as in the pointwise proof
gives
\[
\mathbb E\zeta_{ni,j}^2
=
\frac{1}{p_A(a_j)^2}
\int \{K_{\mathrm{DB}}(u)\}^2
p_A(a_j+hu)v(a_j+hu)\,du
\longrightarrow \sigma_{\mathrm{DB}}^2(a_j).
\]
Thus
\[
\begin{aligned}
\operatorname{Var}\{Z_{ni}(c)\}
&=
\sum_{j=1}^m c_j^2\,\mathbb E\zeta_{ni,j}^2\\
&\longrightarrow
\sum_{j=1}^m c_j^2\sigma_{\mathrm{DB}}^2(a_j)
=
c^T\Sigma_{\mathrm{DB}}c.
\end{aligned}
\]
The limiting variance is positive for nonzero $c$.

To establish a central limit theorem for this linear
combination, we verify the Lyapunov condition.
If the moment exponents in Condition~\ref{cond:clt:moments} differ across
evaluation points, choose a common $\delta>0$ no larger
than any of them.
Such a choice is possible because $m$ is finite.
For all sufficiently large $n$, the kernel window around
$a_j$ lies in its corresponding fixed neighborhood.
The bounded nuisance coefficients and inverse weights,
together with Condition~\ref{cond:clt:moments}, therefore imply
\[
\mathbb E|\zeta_{ni,j}|^{2+\delta}
\le
C_j h^{-\delta/2}
\int |K_{\mathrm{DB}}(u)|^{2+\delta}\,du
=
O(h^{-\delta/2}).
\]
Using the inequality
$|\sum_{j=1}^m x_j|^{2+\delta}
\le m^{1+\delta}\sum_{j=1}^m|x_j|^{2+\delta}$, we obtain
\[
\begin{aligned}
\mathbb E|Z_{ni}(c)|^{2+\delta}
&\le
m^{1+\delta}\sum_{j=1}^m
|c_j|^{2+\delta}
\mathbb E|\zeta_{ni,j}|^{2+\delta}\\
&=
O(h^{-\delta/2}).
\end{aligned}
\]
Consequently,
\[
n^{-\delta/2}\mathbb E|Z_{ni}(c)|^{2+\delta}
=
O((nh)^{-\delta/2})
\longrightarrow0.
\]
Together with the positive variance limit, this verifies the
Lyapunov condition, and hence
\[
\frac1{\sqrt n}\sum_{i=1}^n Z_{ni}(c)
\xrightarrow{d}
\mathcal N(0,c^T\Sigma_{\mathrm{DB}}c).
\]
Slutsky's theorem gives the same limit for the corresponding
linear combination of the scaled estimators.

Since $c$ was arbitrary, the Cram\'er-Wold device \citep[Chapter~2]{vandervaart1998} yields
\[
\sqrt{nh}
\begin{pmatrix}
\hat\theta_{\mathrm{DB}}(a_1)-\theta(a_1)\\
\vdots\\
\hat\theta_{\mathrm{DB}}(a_m)-\theta(a_m)
\end{pmatrix}
\xrightarrow{d}
\mathcal N_m( 0,\Sigma_{\mathrm{DB}}).
\]
The diagonal covariance matrix then implies independence of
the limiting Gaussian coordinates.
\end{proof}

\subsection{Uniform Inference and Confidence Bands}
\label{app:uniform_inference}

We first strengthen the pointwise expansion to a uniform one
and compare the fitted full-equation contributions with their
leading oracle terms. Write \(\kappa_n=h/b\) and let
\(K_{\mathrm{DB},n}\) be the kernel in \eqref{eq:equivalent_kernel}
with \(\kappa\) replaced by \(\kappa_n\). Retain \(\bar\xi\)
from the pointwise analysis and define
\[
W_{\mathcal{A}_C}(O)=\frac{\Delta}{\bar S_C(U\mid X)}
\{\bar\xi(X,U,Y)-\theta(U)\}I(U\in \mathcal{A}_C),\qquad
\varphi_{\mathrm{DB},i}^{\,0}(a)=\frac{W_{\mathcal{A}_C}(O_i)}{h p_A(a)}
K_{\mathrm{DB},n}\!\left(\frac{U_i-a}{h}\right).
\]
The first expression is defined as zero when \(\Delta=0\)
or \(U\notin \mathcal{A}_C\). Under the two-correct-components condition,
conditional centering gives \(E\{W_{\mathcal{A}_C}(O)\mid A=t\}=0\)
for \(t\in \mathcal{A}_C\); hence the oracle contributions have mean zero.
All active kernel windows are contained in \(\mathcal{A}_C\) for
sufficiently large \(n\). Here \(W_{\mathcal{A}_C}(O)=W(O)I(U\in \mathcal{A}_C)\), with \(W\) as defined in the pointwise proof. Hence the displayed kernel-weighted oracle contribution agrees with its pointwise counterpart for all sufficiently large \(n\).

\begin{lemma}[Uniform expansion and approximation of the fitted contributions]
\label{lem:uniform_expansion}
Under the inherited pointwise conditions and
Conditions~\ref{cond:uniform:regularity}-\ref{cond:uniform:bandwidth-grid},
the estimator and the contributions defined in
\eqref{eq:db_full_contribution} satisfy
\begin{equation}
\label{eq:uniform_linear_expansion}
\sup_{a\in \mathcal{A}_0}\left|
\hat\theta_{\mathrm{DB}}(a)-\theta(a)
-\frac1n\sum_{i=1}^n\varphi_{\mathrm{DB},i}^{\,0}(a)
\right|=o_p\{(nh\log n)^{-1/2}\},
\end{equation}
\begin{equation}
\label{eq:uniform_contribution_comparison}
\sup_{a\in \mathcal{A}_0}\left[
\frac hn\sum_{i=1}^n
\{\hat\varphi_{\mathrm{DB},i}(a)-\varphi_{\mathrm{DB},i}^{\,0}(a)\}^2
\right]^{1/2}=o_p\{(\log n)^{-1}\}.
\end{equation}
The second statement also holds after replacing
\(\hat\varphi_{\mathrm{DB},i}(a)\) by
\(\hat\varphi_{\mathrm{DB},i}(a)-\bar\varphi_{\mathrm{DB}}(a)\).
\end{lemma}

\begin{proof}
Write \(L_n=\log n\) and \(\chi_n=q_nr_n+q_ns_n+r_ns_n\),
retain \(\delta_n\) from Condition~\ref{cond:uniform:nuisance}, and let
\(\epsilon_n\) be the maximum, over all folds, of the two
marginalization errors on \(\mathcal{A}_C\). Thus
\(\epsilon_n=O_p(n^{-1/2})\), and the uniform errors of
\(\hat m_{-k}\) and \(\hat p_{A,-k}\) relative to their
coherent limits are bounded by \(C_0\delta_n+\epsilon_n\).
All suprema in this proof are over \(a\in \mathcal{A}_0\), unless
their domains are specified otherwise. The constants in the
bounds below are uniform over the fixed number of folds.
We use the actual ratio \(\kappa_n=h/b\) throughout.

For \(\ell\in\{h,b\}\), set
\[
H_{n,\ell}=\sqrt{\frac{L_n}{n\ell}}
+\frac{L_n}{n^{1-1/q}\ell}.
\]
The bandwidth conditions and \(\kappa_n=O(1)\) imply
\[
H_{n,\ell}=O\!\left(\sqrt{\frac{L_n}{n\ell}}\right),
\qquad H_{n,b}=O(\sqrt{\kappa_n}H_{n,h}).
\]
We first justify the entropy bounds used for the censored
scores. Put \(K_j(u)=u^jK(u)\), \(0\le j\le4\), and
\(K_{j,\ell}(s-a)=\ell^{-1}K_j((s-a)/\ell)\).
For any fitted or limiting nuisance functions \(\eta\),
Fubini's theorem writes a design-matrix tail contribution as
\[
\int K_{j,\ell}(s-a)V_\eta(s,O)\,ds,
\]
where, on the active integration region,
\[
V_\eta(s,O)=f_{A,\eta}(s\mid X)
\left\{
\frac{(1-\Delta)I(U\le s)}{S_{A,\eta}(U\mid X)S_{C,\eta}(U\mid X)}
-\int_0^{U\wedge s}
\frac{\lambda_{C,\eta}(t\mid X)}{S_{A,\eta}(t\mid X)S_{C,\eta}(t\mid X)}\,dt
\right\}.
\]
The first term is evaluated only when \(U\le s\).
A numerator tail has the same form with the additional
factor \(m_\eta(s)\). For all sufficiently large \(n\),
the active values of \(s\) belong to \(\mathcal{A}_C\subset(0,\tau)\).
The inherited bounds give
\(S_{A,\eta}(t\mid x)\ge c(T-\tau)\) for \(t\le\tau\),
and
\[
\int_0^\tau\frac{\lambda_{C,\eta}(t\mid x)}{S_{C,\eta}(t\mid x)}\,dt
=S_{C,\eta}(\tau\mid x)^{-1}-1\le C_0.
\]
Consequently, \(V_\eta\) and \(m_\eta V_\eta\) are uniformly
bounded. This argument applies to fitted hazards through
their survival-hazard relation; it requires no convergence
rate for the hazards themselves.

Each \(K_j\) is compactly supported and Lipschitz. Its
translation bound in \(L_1\) shows that each tail contribution
\(T_a(O)\) has a constant envelope and satisfies
\[
|T_a(O)-T_{a'}(O)|
\le C_0\min\{|a-a'|/\ell,1\}.
\]
Hence its covering number, under \(L_2(Q)\) for any probability
measure \(Q\), is bounded by \(C_0/(\ell\varepsilon)\) at
radius \(C_0\varepsilon\). The same polynomial covering bound,
possibly with a different fixed exponent, holds for the finite
sums and products with the bounded Taylor coefficients used
in the residual equations. Smoothness of the fitted functions
in the treatment argument is not needed: they enter the
integral as bounded weights, and only the kernel is translated.

The point contributions are finite sums of translated
\(K_{j,\ell}\) multiplied by fixed observation functions
after conditioning on a training sample. They have polynomial
covering bounds, squared moments of order \(\ell^{-1}\),
and envelopes bounded by
\(C_0\ell^{-1}\Delta(1+|Y|)I(U\in \mathcal{A}_C)\).
The difference between a fitted and limiting point
contribution has an additional factor
\(\delta_n+\epsilon_n\). The moment and bandwidth conditions
make the logarithmic entropy terms \(O(L_n)\).
Apply the local maximal inequality of
\citet[Theorem~2.1]{vanderVaartWellner2011} to the bounded
tail classes and the finite-moment VC-type inequality of
\citet[Corollary~5.1]{Chernozhukov2014} to the point classes,
conditionally within each evaluation fold. Use the error
for that fold's training sample first, and then bound it by
the maximum across folds. These inequalities give the following
bounds. A centered point-score average is \(O_p(H_{n,\ell})\);
a centered fitted-minus-limiting point-score average is
\(O_p\{(\delta_n+\epsilon_n)H_{n,\ell}\}\); and a centered
tail-score average, or a centered difference of two such
averages, is \(O_p(\sqrt{L_n/n})\).
The fixed number of folds permits these bounds to be combined
without requiring independence between different fold averages.

The integration-by-parts argument in the pointwise rate proof
also holds uniformly on \(\mathcal{A}_0\). Its constants depend only on
the common survival lower bounds and the bounded \(L_1\)
norms of the kernels \(K_j\). Thus the residual-score drift
is \(O_p(\chi_n+\epsilon_n)\), and the design-matrix score drift
is \(O_p(q_nr_n)\). In particular, if
\(A_{n,\ell}\) denotes the uniform error of a fitted design
matrix relative to its finite-bandwidth population matrix,
then
\[
A_{n,\ell}=O_p(H_{n,\ell}+q_nr_n).
\]
The population matrices and their inverses differ from their
local kernel-moment limits by \(O(\ell)\), since \(p_A\)
is continuously differentiable on \(\mathcal{A}_C\). Their smallest
eigenvalues are uniformly bounded away from zero. All fitted
matrices are therefore invertible throughout \(\mathcal{A}_0\), with
probability tending to one, and their inverses are uniformly
bounded on this event.

We next obtain a uniform expansion using the local Taylor
coefficient vectors from the pointwise proof, with bandwidths
\(h\) and \(b\). The uniform
H\"older condition on \(\theta''\) gives
\[
\theta(t)-\theta(a)-\theta'(a)(t-a)
-\tfrac12\theta''(a)(t-a)^2=O(|t-a|^{2+\gamma})
\]
uniformly for the active \((a,t)\). Consequently, the mean
of the local linear residual score is
\(\tfrac12h^2\theta''(a) B_h(a)+O(h^{2+\gamma})\),
and that of the local quadratic residual score is
\(O(b^{2+\gamma})\).
Here \( B_h(a)\) is the population quadratic-term moment vector from
the preceding proof. The quadratic-term moment vector and design bounds give
\[
\sup_a\left|e_1^\top\{\hat{ D}_h(a)\}^{-1}
\hat{ B}_h(a)-\nu_2\right|=O_p(h+A_{n,h}),
\]
\[
\sup_a\|\hat\beta_b^{(2)}(a)-\beta_b^{(2)}(a)\|
=O_p(b^{2+\gamma}+H_{n,b}+\chi_n+\epsilon_n).
\]
In particular,
\(\sup_a|e_3^\top\hat\beta_b^{(2)}(a)|
=O_p(b^2+H_{n,b}+\chi_n+\epsilon_n)\).

Apply the exact debiasing decomposition from the pointwise
proof. The leading mean terms cancel, leaving deterministic
bias of order \(h^{2+\gamma}+h^2b^\gamma\).
Replacing the population quadratic-term moment vector by its empirical version
introduces \(O_p(h^2A_{n,h})\).
On the active kernel windows, the oracle point residual can be written as \(W_{\mathcal{A}_C}(O)\)
plus the difference between \(\theta(U)\) and its Taylor
polynomial. The centered averages of these latter point terms
are bounded by \(h^2H_{n,h}\) and
\(b^{2+\gamma}H_{n,b}\), respectively. The centered oracle
tails are \(O_p(\sqrt{L_n/n})\).

Replacing the empirical inverse matrices and bias projection
in the remaining point fluctuations by their limiting
kernel-moment values introduces an error bounded in probability
by a constant times
\[
(h+A_{n,h})H_{n,h}
+\kappa_n^2(b+A_{n,b}+h+A_{n,h})H_{n,b}.
\]
The resulting two leading point terms combine exactly into
\(n^{-1}\sum_i\varphi_{\mathrm{DB},i}^{\,0}(a)\), with the actual
ratio \(\kappa_n\). To check the case \(\kappa_n\to0\), note
the weighted identities
\[
\kappa_n^{5/2}b=h\kappa_n^{3/2},\qquad
\kappa_n^2b^2=h^2,\qquad
\kappa_n^{5/2}b^{2+\gamma}=h^2b^\gamma\sqrt{\kappa_n}.
\]
Combining the preceding bounds, using \(q_nr_n\le\chi_n\),
gives
\[
\begin{split}
&\sup_a\left|\hat\theta_{\mathrm{DB}}(a)-\theta(a)
-\frac1n\sum_i\varphi_{\mathrm{DB},i}^{\,0}(a)\right|\\
&\quad=O_p\!\left\{
h^2(h^\gamma+b^\gamma)+\chi_n+\epsilon_n
+(\delta_n+h)H_{n,h}+\sqrt{L_n/n}+H_{n,h}^2
\right\}.
\end{split}
\]
After multiplication by \(\sqrt{nhL_n}\), the first term
vanishes by the stated bias condition, the second by the
product-rate condition, and the nuisance empirical term by
\(\delta_nL_n=o_p(1)\). The remaining terms are bounded by
orders \(\sqrt{hL_n}\), \(hL_n\), \(\sqrt hL_n\), and
\(L_n^{3/2}/\sqrt{nh}\), all tending to zero.
Here \(nh^5=O(1)\) makes \(h\) decrease at least
polynomially, and the lower bandwidth bound controls the last
term. This proves the first assertion.

To compare the fitted contributions with their oracle terms,
we use the derivative definition \eqref{eq:db_full_contribution}.
For a local system,
differentiating \(D^{-1}N\) gives
\(D^{-1}(dN-dD\,D^{-1}N)\). Thus its contribution retains
the point and tail residual scores, together with the
training-sample marginalization derivatives described in
Appendix~\ref{app:db_variance}. Differentiating the debiasing
product also retains the derivative of the empirical bias
projection.

The kernel maximal inequality applied to squared point
envelopes gives
\[
\sup_a\frac\ell n\sum_i
\Delta_iK_\ell(U_i-a)^2(1+Y_i^2)=O_p(1),
\qquad \ell\in\{h,b\}.
\]
Indeed, the empirical error here is bounded by
\(O_p[\sqrt{L_n/(n\ell)}+
L_n/\{n^{1-2/q}\ell\}]\), which tends to zero.
The same argument for the first moments gives uniformly
bounded empirical averages of
\(\Delta|K_\ell(U-a)|(1+|Y|)\).

The fitted tail contributions are uniformly bounded.
The marginalization derivative for any training observation
is bounded, uniformly in that observation and in \(a\),
by a constant times
\[
\max_k\left\{1+\frac1{|\mathcal I_{-k}|}
\sum_{j\in\mathcal I_k}
\Delta_j|K_\ell(U_j-a)|(1+|Y_j|)\right\}=O_p(1).
\]
This follows directly from the normalized training-average
derivatives and the inherited bounds on the fitted regression
and treatment density. Both the tails and marginalization
derivatives therefore have order \(O_p(\sqrt h)\) in the
scaled empirical norm
\(\{h n^{-1}\sum_i(\cdot)^2\}^{1/2}\).
The factor \(\kappa_n^2=O(1)\) preserves this bound for the
quadratic contribution.

For the point terms, the local linear fit can be replaced
by the limiting residual \(W_{\mathcal{A}_C}\) at a cost bounded in this
norm by
\[
O_p(\delta_n+\epsilon_n+h+A_{n,h}+h^2+H_{n,h}+\chi_n).
\]
This bound accounts separately for fitted nuisance functions,
the inverse design matrix, and the difference between the
fitted polynomial and \(\theta(U)\).
The corresponding quadratic bound, after the debiasing
factor, is
\[
O_p\!\left[\kappa_n^{5/2}
\{\delta_n+\epsilon_n+b+A_{n,b}
+b^{2+\gamma}+H_{n,b}+\chi_n+h+A_{n,h}\}\right].
\]
The factor \(\sqrt{\kappa_n}\) comes from the scaled
empirical norm of a bandwidth-\(b\) point contribution;
the other factor \(\kappa_n^2\) comes from bias correction.
The terms \(h+A_{n,h}\) include estimation of the bias
projection. The weighted identities above reduce this bound
to the same order as the local linear bound.

The contribution from differentiating the bias projection
has scaled empirical norm \(O_p(1)\), because it consists
of fitted design and bias scores multiplied by uniformly
bounded matrices and coefficients. Its multiplier is
\(\kappa_n^2e_3^\top\hat\beta_b^{(2)}(a)\), whose uniform size is
\[
O_p\{\kappa_n^2(b^2+H_{n,b}+\chi_n+\epsilon_n)\}
=O_p(h^2+H_{n,h}+\chi_n+\epsilon_n).
\]
It therefore also satisfies the required approximation rate.
Once these terms are removed, the two limiting point
contributions are precisely \(\varphi_{\mathrm{DB},i}^{\,0}(a)\).
We have proved the more explicit bound
\[
\sup_a\left[\frac hn\sum_i
\{\hat\varphi_{\mathrm{DB},i}(a)-\varphi_{\mathrm{DB},i}^{\,0}(a)\}^2
\right]^{1/2}
=O_p(\delta_n+\epsilon_n+\sqrt h+H_{n,h}+\chi_n).
\]
Each term on the right is \(o_p(L_n^{-1})\) under the
stated conditions, proving the second assertion.

Finally, the oracle contributions have mean zero, and the
same kernel maximal bound gives
\[
\sup_a\left|\frac{\sqrt h}{n}
\sum_i\varphi_{\mathrm{DB},i}^{\,0}(a)\right|
=O_p(\sqrt hH_{n,h})=o_p(L_n^{-1}).
\]
Cauchy-Schwarz and the preceding empirical norm comparison
therefore show that subtracting
\(\bar\varphi_{\mathrm{DB}}(a)\) preserves the asserted rate.
\end{proof}

\begin{proof}[Proof of Theorem~\ref{thm:uniform_band}]
Write \(P\) for expectation under the true distribution and put \(L_n=\log n\). Use the oracle contributions \(\varphi_{\mathrm{DB}}^{\,0}(a)\) from Lemma~\ref{lem:uniform_expansion}, with the actual ratio \(\kappa_n=h/b\). Conditional independent censoring, coherent marginal limits, and the two-correct-components condition give
\[
E\{W_{\mathcal{A}_C}(O)\mid A=t\}=0,\qquad E\{W_{\mathcal{A}_C}(O)^2\mid A=t\}=v(t),\qquad t\in \mathcal{A}_C.
\]
For example, if \(\bar S_C=S_C\), the first identity follows from the conditional centering of the full-data pseudo-outcome. Otherwise both \(\bar\mu=\mu\) and \(\bar f_A=f_A\), so that \(\bar m=\theta\), and the residual has mean zero given \((X,A)\), including after multiplication by \(S_C/\bar S_C\). The second identity follows by conditioning on \((X,A,Y)\). Thus \(P\varphi_{\mathrm{DB}}^{\,0}(a)=0\).

For Gaussian approximation, use the finite-bandwidth variance \(\sigma_{\mathrm{DB},n}^2(a)=hP\{\varphi_{\mathrm{DB}}^{\,0}(a)^2\}\). Direct calculation gives
\[
\sigma_{\mathrm{DB},n}^2(a)=\frac1{p_A(a)^2}\int\{K_{\mathrm{DB},n}(u)\}^2p_A(a+hu)v(a+hu)\,du.
\]
The integral only involves treatment values within distance \(h\vee b\) of \(a\). Uniform continuity of \(p_Av\), together with \(K_{\mathrm{DB},n}\to K_{\mathrm{DB}}\) in \(L_2\), gives
\[
\sup_{a\in \mathcal{A}_0}|\sigma_{\mathrm{DB},n}^2(a)-\sigma_{\mathrm{DB}}^2(a)|\longrightarrow0,\qquad
0<c\le\inf_{a\in \mathcal{A}_0}\sigma_{\mathrm{DB},n}(a)\le\sup_{a\in \mathcal{A}_0}\sigma_{\mathrm{DB},n}(a)\le C_0
\]
for all sufficiently large \(n\).

Define \(f_{n,a}=\sqrt h\,\varphi_{\mathrm{DB}}^{\,0}(a)/\sigma_{\mathrm{DB},n}(a)\) and \(\mathcal F_n=\{f_{n,a},-f_{n,a}:a\in \mathcal{A}_0\}\). Then \(Pf_{n,a}=0\) and \(Pf_{n,a}^2=1\). To verify regularity of this class, write \(L(u)=(u^2-\nu_2)K(u)/(\nu_4-\nu_2^2)\), so that \(K_{\mathrm{DB},n}(u)=K(u)-\nu_2\kappa_n^3L(\kappa_nu)\). Compact support and Lipschitz continuity imply that \(K,L\) have bounded variation and square-integrable weak derivatives. Since \(\kappa_n=O(1)\),
\[
\sup_n\{\|K_{\mathrm{DB},n}\|_2+\|(K_{\mathrm{DB},n})'\|_2+\|K_{\mathrm{DB},n}\|_\infty\}<\infty.
\]
These bounds also cover \(\kappa=0\); for example, the derivative of the correction term has \(L_2\) norm \(\nu_2\kappa_n^{7/2}\|L'\|_2\).

The translation inequality in \(L_2\), the bounds on \(p_Av\), and normalization in \(L_2(P)\) therefore imply
\begin{equation}
\label{eq:uniform_metric_bound}
\|f_{n,u}-f_{n,v}\|_{P,2}\le C_0\left(\frac{|u-v|}{h}\wedge1\right),\qquad u,v\in \mathcal{A}_0.
\end{equation}
Here \(f_{n,a}\) is the \(L_2(P)\)-normalized version of \(h^{-1/2}W_{\mathcal{A}_C}(O)K_{\mathrm{DB},n}((U-a)/h)\); the reverse triangle inequality controls the change in its normalizing norm.

The translation and dilation classes of bounded-variation kernels are VC type \citep{GineGuillou2002}. Finite sums, multiplication by the fixed function \(W_{\mathcal{A}_C}\), and bounded scalar factors preserve polynomial covering bounds. Hence \(\mathcal F_n\) is VC type with constants independent of \(n\). It is pointwise measurable by continuity in \(a\). An envelope is \(F_n(O)=C_0h^{-1/2}\Delta(1+|Y|)I(U\in \mathcal{A}_C)\). Condition~\ref{cond:uniform:moments} and kernel localization give
\[
\|F_n\|_{P,q}=O(h^{-1/2}),\qquad
\sup_{f\in\mathcal F_n}P|f|^3=O(h^{-1/2}),\qquad
\sup_{f\in\mathcal F_n}P|f|^4=O(h^{-1}).
\]

The metric bound \eqref{eq:uniform_metric_bound} gives \(L_2(P)\) covering numbers at most \(C_0/(h\epsilon)\). The Gaussian entropy integral is finite \citep[Section~2.2]{vanderVaartWellner1996}, so there exists a separable centered Gaussian process \(Z_n^0\) with covariance \(P(f_{n,u}f_{n,v})\), and
\[
E\sup_{a\in \mathcal{A}_0}|Z_n^0(a)|=O\{\sqrt{\log(1/h)}\}=O(\sqrt{L_n}).
\]
The bandwidth conditions and \(nh^5=O(1)\) imply \(\log(1/h)\asymp L_n\).

Let \(\mathbb G_n(a)=n^{-1/2}\sum_i f_{n,a}(O_i)\). Applying \citet[Corollary~2.2]{Chernozhukov2014} to \(\mathcal F_n\), with its probability parameter denoted by \(\zeta_n\), couples the empirical and Gaussian suprema with error
\[
O_p\!\left\{
\frac{L_n}{\sqrt{\zeta_n n^{1-2/q}h}}
+\frac{L_n^{3/4}}{\sqrt{\zeta_n}(nh)^{1/4}}
+\frac{L_n^{2/3}}{\zeta_n^{1/3}(nh)^{1/6}}
\right\}.
\]
Condition~\ref{cond:uniform:bandwidth} permits a deterministic \(\zeta_n\downarrow0\) sufficiently slowly that this error is \(o_p(L_n^{-1/2})\), with coupling failure probability tending to zero. Indeed, the second, fourth, and sixth powers of the three terms after multiplication by \(\sqrt{L_n}\) are controlled by \(L_n^3/(\zeta_n n^{1-2/q}h)\), \(L_n^5/(\zeta_n^2nh)\), and \(L_n^7/(\zeta_n^2nh)\), respectively.

All coordinates of \(Z_n^0\) have variance one. Gaussian anti-concentration \citep[Theorem~3]{Chernozhukov2015}, applied to finite sets with their negatives and then to a countable dense subset, gives
\begin{equation}
\label{eq:uniform_anticoncentration}
\sup_{t\in\mathbb R}P\left\{\left|\sup_{a\in \mathcal{A}_0}|Z_n^0(a)|-t\right|\le\epsilon\right\}\le C_0\epsilon\sqrt{L_n}.
\end{equation}
The same bound holds for the grid maximum. The coupling and \eqref{eq:uniform_anticoncentration}, equivalently \citet[Lemma~2.4]{Chernozhukov2014}, yield a Kolmogorov approximation of \(\sup_a|\mathbb G_n(a)|\) by \(\sup_a|Z_n^0(a)|\). By \eqref{eq:uniform_linear_expansion}, the same approximation holds for
\[
T_n^0=\sup_{a\in \mathcal{A}_0}\frac{\sqrt{nh}|\hat\theta_{\mathrm{DB}}(a)-\theta(a)|}{\sigma_{\mathrm{DB},n}(a)}.
\]
In particular, \(T_n^0=O_p(\sqrt{L_n})\).

We next replace the continuous index set by the grid. By \eqref{eq:uniform_metric_bound}, the Gaussian entropy bound for increments \citep[Section~2.2]{vanderVaartWellner1996} gives
\[
E\sup_{|u-v|\le\omega_n}|Z_n^0(u)-Z_n^0(v)|
\le C_0\frac{\omega_n}{h}\sqrt{\log(C_0/\omega_n)}.
\]
The supremum is restricted to \(u,v\in \mathcal{A}_0\). Since a grid of \(N_n\) points covering a nondegenerate interval has \(\omega_n\ge c/N_n\), its polynomial size implies \(\log(C_0/\omega_n)=O(L_n)\). Condition~\ref{cond:uniform:grid} and Markov's inequality therefore give
\[
\left|\sup_{a\in \mathcal{A}_0}|Z_n^0(a)|-\max_{a\in \mathcal{A}_n}|Z_n^0(a)|\right|=o_p(L_n^{-1/2}).
\]
Together with \eqref{eq:uniform_anticoncentration}, this validates replacement of the continuous Gaussian supremum by its grid maximum.

To control studentization and the oracle empirical covariance, apply the VC-type maximal inequality of \citet[Corollary~5.1]{Chernozhukov2014} to the product class \(\{f_{n,u}f_{n,v}:u,v\in \mathcal{A}_0\}\). Its envelope is \(F_n^2\), and
\[
\sup_{u,v}P(f_{n,u}f_{n,v})^2=O(h^{-1}),\qquad
\|F_n^2\|_{P,2}=O(h^{-1}),\qquad
\left\|\max_{i\le n}F_n(O_i)^2\right\|_2=O(n^{2/q}/h).
\]
The last inequality follows from \(q\ge4\) and the \(L_q\) bound on \(F_n\). Thus
\begin{equation}
\label{eq:uniform_oracle_covariance}
\sup_{u,v}|(\mathbb P_n-P)(f_{n,u}f_{n,v})|
=O_p\!\left\{\sqrt{\frac{L_n}{nh}}+\frac{L_n}{n^{1-2/q}h}\right\}=o_p(L_n^{-2}).
\end{equation}
Applying the same inequality to \(\mathcal F_n\) gives \(\sup_a|\mathbb P_nf_{n,a}|=O_p(\sqrt{L_n/n})\). Centering the oracle contributions therefore adds at most \(O_p(L_n/n)=o_p(L_n^{-2})\) to \eqref{eq:uniform_oracle_covariance}.

By \eqref{eq:uniform_contribution_comparison} and the lower variance bound,
\[
\sup_{a\in \mathcal{A}_0}\left[\frac1n\sum_i\left\{\frac{\sqrt h\,\hat\varphi_{\mathrm{DB},i}(a)}{\sigma_{\mathrm{DB},n}(a)}-f_{n,a}(O_i)\right\}^2\right]^{1/2}=o_p(L_n^{-1}).
\]
Empirical centering is a contraction in this squared norm. Cauchy-Schwarz, the diagonal case of \eqref{eq:uniform_oracle_covariance}, and the variance definition \eqref{eq:db_variance_estimator} consequently give
\begin{equation}
\label{eq:uniform_studentization}
\sup_{a\in \mathcal{A}_0}\left|\frac{\hat\sigma_{\mathrm{DB}}(a)}{\sigma_{\mathrm{DB},n}(a)}-1\right|=o_p(L_n^{-1}).
\end{equation}
In particular, the estimated standard deviations are uniformly positive with probability tending to one. Write \(T_n=\sup_{a\in \mathcal{A}_0}\sqrt{nh}|\hat\theta_{\mathrm{DB}}(a)-\theta(a)|/\hat\sigma_{\mathrm{DB}}(a)\). Since \(T_n^0=O_p(\sqrt{L_n})\), \eqref{eq:uniform_studentization} implies \(|T_n-T_n^0|=o_p(L_n^{-1/2})\). Studentization therefore preserves the preceding distributional approximations.

It remains to validate the feasible multiplier process. On the same multiplier probability space as \eqref{eq:uniform_multiplier}, define the oracle empirical multiplier process
\[
Z_n^\circ(a)=\frac1{\sqrt n}\sum_i G_i\{f_{n,a}(O_i)-\mathbb P_nf_{n,a}\},\qquad a\in \mathcal{A}_n.
\]
Its conditional covariance differs uniformly from \(P(f_{n,u}f_{n,v})\) by \(o_p(L_n^{-2})\), by \eqref{eq:uniform_oracle_covariance} and the centering bound. Apply \citet[Theorem~2]{Chernozhukov2015} to the \(2N_n\) coordinates obtained by adjoining negatives. If the maximum covariance discrepancy is \(\Delta_n\), the Kolmogorov comparison bound is
\[
C_0\Delta_n^{1/3}\{1\vee\log(2N_n/\Delta_n)\}^{2/3}=o_p(1),
\]
with bound zero when \(\Delta_n=0\). This follows from \(N_n=O(n^{d_0})\), \(\Delta_nL_n^2=o_p(1)\), and \(\Delta_n\log^2(1/\Delta_n)=o_p(1)\).

We finally replace the oracle coefficients by the fitted full-equation contributions. Equations~\eqref{eq:uniform_contribution_comparison} and \eqref{eq:uniform_studentization} imply
\[
\sup_{a\in \mathcal{A}_0}\left[\frac1n\sum_i\left\{
\frac{\sqrt h\{\hat\varphi_{\mathrm{DB},i}(a)-\bar\varphi_{\mathrm{DB}}(a)\}}{\hat\sigma_{\mathrm{DB}}(a)}
-f_{n,a}(O_i)+\mathbb P_nf_{n,a}\right\}^2\right]^{1/2}=o_p(L_n^{-1}).
\]
Conditional on the data, the difference \(Z_n^G-Z_n^\circ\) is a Gaussian vector because the same multipliers are used. The Gaussian maximal inequality gives
\[
E_G\max_{a\in \mathcal{A}_n}|Z_n^G(a)-Z_n^\circ(a)|
\le C_0\sqrt{\log(2N_n)}\,o_p(L_n^{-1})=o_p(L_n^{-1/2}),
\]
where \(E_G\) denotes expectation over the multipliers given the data. Conditional Markov's inequality, oracle Gaussian comparison, and \eqref{eq:uniform_anticoncentration} now validate the fitted multiplier maximum.

Let \(F_n^0(t)=P\{\max_{a\in \mathcal{A}_n}|Z_n^0(a)|\le t\}\), a deterministic distribution function, and let \(\hat F_n(t)\) be the conditional distribution function of \(\max_{a\in \mathcal{A}_n}|Z_n^G(a)|\). We have proved
\[
\sup_t|P(T_n\le t)-F_n^0(t)|\longrightarrow0,\qquad
\sup_t|\hat F_n(t)-F_n^0(t)|=o_p(1).
\]

To evaluate coverage at the random critical value, choose deterministic \(\eta_n\downarrow0\) such that \(P\{\sup_t|\hat F_n(t)-F_n^0(t)|>\eta_n\}\to0\). Denote the \(s\)-quantile of \(F_n^0\) by \(q_n^0(s)\). The Gaussian maxima have continuous distribution functions, since their coordinates have positive variances. On the event of uniform distance at most \(\eta_n\),
\[
q_n^0(1-\alpha-\eta_n)\le c_{1-\alpha,n}\le q_n^0(1-\alpha+\eta_n).
\]
The two deterministic thresholds and the approximation of the distribution of \(T_n\) therefore give
\[
1-\alpha-\eta_n-o(1)\le P(T_n\le c_{1-\alpha,n})\le1-\alpha+\eta_n+o(1).
\]
Consequently, \(P(T_n\le c_{1-\alpha,n})\to1-\alpha\), which is the simultaneous coverage assertion for \eqref{eq:uniform_band}.
\end{proof}

\section{Additional simulation results}
\label{sec:sim-supplement}

Tables~\ref{tab:sim-main-linear} and~\ref{tab:sim-main-nonlinear}
report bias, RMSE, pointwise coverage, and mean interval length at
$a=0.8,1.3,1.8$. Both bandwidth rules and all three sample sizes
are included. PI estimates are repeated across bandwidths because
PI has no second-stage local smoothing; dashes indicate that
its confidence intervals are not evaluated.

For the main experiment, the outcome regression is fitted by
Nystr\"om Gaussian kernel ridge regression using uncensored training
observations. The treatment density is a Bernstein/beta mixture,
and the censoring hazard is a positive Bernstein expansion;
their covariate-dependent coefficients are learned by neural
networks using observed-data likelihoods. Conventional procedures
refit the outcome and treatment learners within complete cases.
Tuning uses an internal validation sample within each training
fold, followed by refitting. PI averages out-of-fold outcome
predictions over all individuals, including censored individuals.

\begin{table}[ht]
\centering
\caption{Main experiment: Setting 1 (linear).}
\label{tab:sim-main-linear}
\begingroup
\renewcommand{\baselinestretch}{1}
\fontsize{8.5}{10}\selectfont
\setlength{\tabcolsep}{1.6pt}
\renewcommand{\arraystretch}{1.10}
\begin{tabular}{@{}rl*{12}{r}@{}}
\toprule
\multirow{2}{*}{$n$} & \multirow{2}{*}{Method} & \multicolumn{4}{c}{$a=0.8$} & \multicolumn{4}{c}{$a=1.3$} & \multicolumn{4}{c}{$a=1.8$} \\
\cmidrule(lr){3-6}\cmidrule(lr){7-10}\cmidrule(lr){11-14}
 & & Bias & RMSE & CP (\%) & Length & Bias & RMSE & CP (\%) & Length & Bias & RMSE & CP (\%) & Length \\
\midrule
\multicolumn{14}{l}{\textit{Panel A: Regular bandwidth, $h=2n^{-1/5}$}} \\
\multirow{5}{*}{2000} & Conv. & -0.696 & 0.701 & 0.0 & 0.354 & -0.772 & 0.779 & 0.0 & 0.406 & -0.842 & 0.851 & 0.0 & 0.468 \\
 & Conv.DB & -0.696 & 0.703 & 0.0 & 0.387 & -0.773 & 0.781 & 0.0 & 0.448 & -0.842 & 0.853 & 0.0 & 0.525 \\
 & PI & -0.002 & 0.074 & \textemdash & \textemdash & -0.002 & 0.084 & \textemdash & \textemdash & -0.004 & 0.098 & \textemdash & \textemdash \\
 & MR & 0.000 & 0.079 & 95.2 & 0.318 & -0.004 & 0.094 & 94.9 & 0.378 & 0.000 & 0.117 & 94.9 & 0.462 \\
 & MRDB & 0.000 & 0.089 & 95.4 & 0.357 & -0.005 & 0.109 & 94.6 & 0.432 & 0.000 & 0.137 & 94.5 & 0.538 \\
\addlinespace[3pt]
\multirow{5}{*}{4000} & Conv. & -0.698 & 0.700 & 0.0 & 0.256 & -0.769 & 0.773 & 0.0 & 0.294 & -0.846 & 0.850 & 0.0 & 0.339 \\
 & Conv.DB & -0.698 & 0.701 & 0.0 & 0.282 & -0.769 & 0.773 & 0.0 & 0.327 & -0.846 & 0.851 & 0.0 & 0.383 \\
 & PI & -0.002 & 0.054 & \textemdash & \textemdash & 0.002 & 0.060 & \textemdash & \textemdash & -0.001 & 0.068 & \textemdash & \textemdash \\
 & MR & 0.001 & 0.058 & 95.7 & 0.231 & 0.002 & 0.071 & 94.8 & 0.275 & -0.002 & 0.085 & 95.1 & 0.335 \\
 & MRDB & 0.001 & 0.066 & 95.1 & 0.261 & 0.002 & 0.083 & 94.2 & 0.317 & -0.003 & 0.100 & 94.4 & 0.394 \\
\addlinespace[3pt]
\multirow{5}{*}{8000} & Conv. & -0.698 & 0.700 & 0.0 & 0.185 & -0.773 & 0.775 & 0.0 & 0.213 & -0.844 & 0.846 & 0.0 & 0.246 \\
 & Conv.DB & -0.698 & 0.700 & 0.0 & 0.205 & -0.773 & 0.775 & 0.0 & 0.239 & -0.843 & 0.846 & 0.0 & 0.280 \\
 & PI & -0.002 & 0.038 & \textemdash & \textemdash & 0.000 & 0.042 & \textemdash & \textemdash & -0.003 & 0.049 & \textemdash & \textemdash \\
 & MR & 0.000 & 0.042 & 95.2 & 0.168 & -0.002 & 0.052 & 93.9 & 0.201 & -0.001 & 0.062 & 95.2 & 0.245 \\
 & MRDB & -0.001 & 0.048 & 95.1 & 0.191 & -0.002 & 0.061 & 94.7 & 0.233 & 0.000 & 0.074 & 95.2 & 0.291 \\
\addlinespace[4pt]
\multicolumn{14}{l}{\textit{Panel B: Undersmoothed bandwidth, $h=n^{-1/4}$}} \\
\multirow{5}{*}{2000} & Conv. & -0.695 & 0.704 & 0.0 & 0.446 & -0.770 & 0.782 & 0.0 & 0.525 & -0.844 & 0.859 & 0.1 & 0.622 \\
 & Conv.DB & -0.695 & 0.707 & 0.0 & 0.520 & -0.766 & 0.782 & 0.2 & 0.620 & -0.847 & 0.867 & 1.1 & 0.743 \\
 & PI & -0.002 & 0.074 & \textemdash & \textemdash & -0.002 & 0.084 & \textemdash & \textemdash & -0.004 & 0.098 & \textemdash & \textemdash \\
 & MR & 0.000 & 0.106 & 95.6 & 0.423 & -0.003 & 0.133 & 95.1 & 0.521 & -0.002 & 0.170 & 94.6 & 0.645 \\
 & MRDB & 0.000 & 0.126 & 95.6 & 0.506 & 0.001 & 0.161 & 95.1 & 0.632 & -0.005 & 0.208 & 94.3 & 0.785 \\
\addlinespace[3pt]
\multirow{5}{*}{4000} & Conv. & -0.696 & 0.702 & 0.0 & 0.332 & -0.767 & 0.773 & 0.0 & 0.390 & -0.844 & 0.853 & 0.0 & 0.466 \\
 & Conv.DB & -0.696 & 0.703 & 0.0 & 0.391 & -0.767 & 0.776 & 0.0 & 0.464 & -0.844 & 0.856 & 0.0 & 0.562 \\
 & PI & -0.002 & 0.054 & \textemdash & \textemdash & 0.002 & 0.060 & \textemdash & \textemdash & -0.001 & 0.068 & \textemdash & \textemdash \\
 & MR & 0.002 & 0.081 & 95.5 & 0.316 & 0.003 & 0.101 & 93.6 & 0.390 & 0.000 & 0.128 & 94.5 & 0.489 \\
 & MRDB & 0.003 & 0.097 & 95.6 & 0.382 & 0.003 & 0.124 & 94.5 & 0.476 & 0.002 & 0.160 & 94.3 & 0.603 \\
\addlinespace[3pt]
\multirow{5}{*}{8000} & Conv. & -0.700 & 0.702 & 0.0 & 0.246 & -0.772 & 0.776 & 0.0 & 0.292 & -0.844 & 0.849 & 0.0 & 0.347 \\
 & Conv.DB & -0.700 & 0.704 & 0.0 & 0.293 & -0.772 & 0.778 & 0.0 & 0.351 & -0.845 & 0.852 & 0.0 & 0.422 \\
 & PI & -0.002 & 0.038 & \textemdash & \textemdash & 0.000 & 0.042 & \textemdash & \textemdash & -0.003 & 0.049 & \textemdash & \textemdash \\
 & MR & -0.002 & 0.060 & 95.6 & 0.237 & 0.000 & 0.077 & 93.9 & 0.294 & -0.002 & 0.097 & 94.1 & 0.370 \\
 & MRDB & -0.002 & 0.073 & 95.8 & 0.288 & -0.001 & 0.094 & 94.2 & 0.361 & -0.003 & 0.120 & 94.1 & 0.459 \\
\bottomrule
\end{tabular}
\par\vspace{4pt}
\begin{minipage}{\textwidth}
\fontsize{8.5}{10}\selectfont \textit{Notes:} Each entry uses $R=1000$ repetitions. Bias is signed Monte Carlo bias relative to $\theta(a)$; RMSE is root mean squared error. CP is the coverage percentage of nominal $95\%$ pointwise intervals; Length is their mean length. Conv. and Conv.DB denote Conventional and ConventionalDB. The pilot bandwidth is $b=1.25h$. PI has no second-stage bandwidth: its point estimates are repeated in both panels, and its interval entries are not evaluated (\textemdash). The reported doses are the endpoints and midpoint of $[0.8,1.8]$.
\end{minipage}
\endgroup
\end{table}

\begin{table}[H]
\centering
\caption{Main experiment: Setting 2 (nonlinear).}
\label{tab:sim-main-nonlinear}
\begingroup
\renewcommand{\baselinestretch}{1}
\fontsize{8.5}{10}\selectfont
\setlength{\tabcolsep}{1.6pt}
\renewcommand{\arraystretch}{1.10}
\begin{tabular}{@{}rl*{12}{r}@{}}
\toprule
\multirow{2}{*}{$n$} & \multirow{2}{*}{Method} & \multicolumn{4}{c}{$a=0.8$} & \multicolumn{4}{c}{$a=1.3$} & \multicolumn{4}{c}{$a=1.8$} \\
\cmidrule(lr){3-6}\cmidrule(lr){7-10}\cmidrule(lr){11-14}
 & & Bias & RMSE & CP (\%) & Length & Bias & RMSE & CP (\%) & Length & Bias & RMSE & CP (\%) & Length \\
\midrule
\multicolumn{14}{l}{\textit{Panel A: Regular bandwidth, $h=2n^{-1/5}$}} \\
\multirow{5}{*}{2000} & Conv. & -0.469 & 0.478 & 0.1 & 0.357 & -0.545 & 0.555 & 0.1 & 0.409 & -0.613 & 0.625 & 0.1 & 0.472 \\
 & Conv.DB & -0.695 & 0.702 & 0.0 & 0.388 & -0.773 & 0.782 & 0.0 & 0.449 & -0.838 & 0.850 & 0.0 & 0.525 \\
 & PI & 0.025 & 0.088 & \textemdash & \textemdash & 0.015 & 0.106 & \textemdash & \textemdash & 0.027 & 0.113 & \textemdash & \textemdash \\
 & MR & 0.232 & 0.247 & 19.4 & 0.321 & 0.227 & 0.248 & 37.0 & 0.383 & 0.232 & 0.264 & 50.2 & 0.467 \\
 & MRDB & 0.005 & 0.092 & 95.3 & 0.357 & -0.003 & 0.113 & 95.3 & 0.434 & 0.005 & 0.145 & 93.1 & 0.538 \\
\addlinespace[3pt]
\multirow{5}{*}{4000} & Conv. & -0.531 & 0.535 & 0.0 & 0.257 & -0.604 & 0.609 & 0.0 & 0.295 & -0.678 & 0.683 & 0.0 & 0.341 \\
 & Conv.DB & -0.703 & 0.707 & 0.0 & 0.282 & -0.777 & 0.781 & 0.0 & 0.327 & -0.851 & 0.856 & 0.0 & 0.384 \\
 & PI & 0.013 & 0.062 & \textemdash & \textemdash & 0.000 & 0.073 & \textemdash & \textemdash & 0.016 & 0.076 & \textemdash & \textemdash \\
 & MR & 0.172 & 0.182 & 17.3 & 0.232 & 0.171 & 0.185 & 31.8 & 0.277 & 0.172 & 0.192 & 47.6 & 0.339 \\
 & MRDB & -0.001 & 0.065 & 94.8 & 0.261 & -0.003 & 0.081 & 95.2 & 0.318 & -0.002 & 0.102 & 94.8 & 0.396 \\
\addlinespace[3pt]
\multirow{5}{*}{8000} & Conv. & -0.565 & 0.567 & 0.0 & 0.186 & -0.638 & 0.640 & 0.0 & 0.213 & -0.708 & 0.712 & 0.0 & 0.247 \\
 & Conv.DB & -0.695 & 0.698 & 0.0 & 0.205 & -0.769 & 0.771 & 0.0 & 0.239 & -0.840 & 0.843 & 0.0 & 0.281 \\
 & PI & 0.007 & 0.040 & \textemdash & \textemdash & -0.005 & 0.046 & \textemdash & \textemdash & 0.013 & 0.053 & \textemdash & \textemdash \\
 & MR & 0.132 & 0.139 & 13.3 & 0.168 & 0.132 & 0.141 & 26.2 & 0.202 & 0.134 & 0.150 & 44.8 & 0.247 \\
 & MRDB & 0.001 & 0.050 & 94.2 & 0.191 & 0.000 & 0.059 & 94.3 & 0.233 & 0.003 & 0.079 & 93.4 & 0.291 \\
\addlinespace[4pt]
\multicolumn{14}{l}{\textit{Panel B: Undersmoothed bandwidth, $h=n^{-1/4}$}} \\
\multirow{5}{*}{2000} & Conv. & -0.666 & 0.675 & 0.0 & 0.446 & -0.746 & 0.759 & 0.1 & 0.526 & -0.812 & 0.828 & 0.4 & 0.626 \\
 & Conv.DB & -0.690 & 0.702 & 0.0 & 0.520 & -0.774 & 0.791 & 0.1 & 0.620 & -0.839 & 0.862 & 1.9 & 0.748 \\
 & PI & 0.025 & 0.088 & \textemdash & \textemdash & 0.015 & 0.106 & \textemdash & \textemdash & 0.027 & 0.113 & \textemdash & \textemdash \\
 & MR & 0.033 & 0.112 & 94.9 & 0.423 & 0.025 & 0.141 & 93.7 & 0.522 & 0.031 & 0.179 & 92.5 & 0.648 \\
 & MRDB & 0.009 & 0.125 & 96.2 & 0.505 & -0.003 & 0.170 & 92.8 & 0.630 & 0.004 & 0.217 & 92.9 & 0.790 \\
\addlinespace[3pt]
\multirow{5}{*}{4000} & Conv. & -0.683 & 0.688 & 0.0 & 0.331 & -0.758 & 0.765 & 0.0 & 0.391 & -0.833 & 0.842 & 0.0 & 0.465 \\
 & Conv.DB & -0.701 & 0.708 & 0.0 & 0.391 & -0.777 & 0.786 & 0.0 & 0.466 & -0.851 & 0.863 & 0.0 & 0.560 \\
 & PI & 0.013 & 0.062 & \textemdash & \textemdash & 0.000 & 0.073 & \textemdash & \textemdash & 0.016 & 0.076 & \textemdash & \textemdash \\
 & MR & 0.019 & 0.083 & 93.8 & 0.316 & 0.016 & 0.103 & 93.3 & 0.391 & 0.016 & 0.130 & 93.0 & 0.491 \\
 & MRDB & 0.002 & 0.098 & 94.3 & 0.381 & -0.002 & 0.123 & 94.8 & 0.477 & -0.002 & 0.161 & 93.0 & 0.603 \\
\addlinespace[3pt]
\multirow{5}{*}{8000} & Conv. & -0.682 & 0.685 & 0.0 & 0.247 & -0.754 & 0.757 & 0.0 & 0.292 & -0.828 & 0.833 & 0.0 & 0.349 \\
 & Conv.DB & -0.695 & 0.699 & 0.0 & 0.294 & -0.766 & 0.771 & 0.0 & 0.351 & -0.843 & 0.851 & 0.0 & 0.424 \\
 & PI & 0.007 & 0.040 & \textemdash & \textemdash & -0.005 & 0.046 & \textemdash & \textemdash & 0.013 & 0.053 & \textemdash & \textemdash \\
 & MR & 0.015 & 0.064 & 93.6 & 0.237 & 0.015 & 0.076 & 94.7 & 0.295 & 0.015 & 0.101 & 93.6 & 0.372 \\
 & MRDB & 0.002 & 0.075 & 93.9 & 0.288 & 0.003 & 0.092 & 94.7 & 0.363 & -0.001 & 0.125 & 93.8 & 0.461 \\
\bottomrule
\end{tabular}
\par\vspace{4pt}
\begin{minipage}{\textwidth}
\fontsize{8.5}{10}\selectfont \textit{Notes:} Each entry uses $R=1000$ repetitions. Bias is signed Monte Carlo bias relative to $\theta(a)$; RMSE is root mean squared error. CP is the coverage percentage of nominal $95\%$ pointwise intervals; Length is their mean length. Conv. and Conv.DB denote Conventional and ConventionalDB. The pilot bandwidth is $b=1.25h$. PI has no second-stage bandwidth: its point estimates are repeated in both panels, and its interval entries are not evaluated (\textemdash). The reported doses are the endpoints and midpoint of $[0.8,1.8]$.
\end{minipage}
\endgroup
\end{table}

The stability experiment uses the same data-generating processes, sample
sizes, dose grid, bandwidths, and number of replications as the
main experiment. The correct outcome model has regressors
$(1,A,A^2,X,AX)$; the treatment and censoring models are truncated
exponential and exponential models, respectively, with log-rates
linear in $X$. To misspecify $\mu$, we omit $X$ and $AX$ from
the outcome model; to misspecify $f_A$ or $S_C$, we omit $X$
from the corresponding log-rate. Only one model is misspecified
at a time, with the other two correctly specified. All parameters are
estimated using two-fold cross-fitting, and the three scenarios
share the same generated sample within each replication.

Tables~\ref{tab:sim-stability-regular}
and~\ref{tab:sim-stability-undersmooth} report results at
$a=0.8,1.3,1.8$ for $n=2000,4000,8000$ in both settings,
based on $1000$ Monte Carlo replications per configuration.
The bandwidths follow the main experiment:
$h=2n^{-1/5}$ for the regular rule and $h=n^{-1/4}$
for undersmoothing, with $b=1.25h$ in both cases.
Under the regular bandwidth, MRDB retains small bias
and generally well-calibrated intervals under each
single-model misspecification, supporting the
multiple-robustness property. Its RMSE and mean interval
length generally decrease as the sample size increases.

\begingroup
\renewcommand{\baselinestretch}{1}
\fontsize{9}{10.5}\selectfont
\setlength{\tabcolsep}{3pt}
\renewcommand{\arraystretch}{1.12}
\setlength{\LTleft}{\fill}
\setlength{\LTright}{\fill}
\setlength{\LTpre}{6pt}
\setlength{\LTpost}{4pt}
\begin{longtable}{@{}rl*{8}{r}@{}}
\caption{Stability experiment at representative doses with the regular bandwidth.}\label{tab:sim-stability-regular}\\
\toprule
\multirow{2}{*}{$n$} & \multirow{2}{*}{Misspecified} & \multicolumn{4}{c}{MR} & \multicolumn{4}{c}{MRDB} \\
\cmidrule(lr){3-6}\cmidrule(lr){7-10}
 & & Bias & RMSE & CP (\%) & Length & Bias & RMSE & CP (\%) & Length \\
\midrule
\endfirsthead
\multicolumn{10}{c}{\tablename\ \thetable\ (continued)}\\
\toprule
\multirow{2}{*}{$n$} & \multirow{2}{*}{Misspecified} & \multicolumn{4}{c}{MR} & \multicolumn{4}{c}{MRDB} \\
\cmidrule(lr){3-6}\cmidrule(lr){7-10}
 & & Bias & RMSE & CP (\%) & Length & Bias & RMSE & CP (\%) & Length \\
\midrule
\endhead
\midrule
\multicolumn{10}{r@{}}{\footnotesize Continued on next page}\\
\endfoot
\bottomrule
\addlinespace[4pt]
\multicolumn{10}{@{}c@{}}{\makebox[0pt][c]{\parbox[t]{\textwidth}{\fontsize{8.5}{10}\selectfont \textit{Notes:} Each entry uses $R=1000$ repetitions. The Misspecified column identifies the only misspecified model; the other two models are correctly specified. Bias is the signed pointwise Monte Carlo bias, not RMSBias. RMSE, CP, and Length are pointwise root mean squared error, nominal $95\%$ interval coverage (\%), and mean interval length, respectively. The bandwidth is $h=2n^{-1/5}$, with $b=1.25h$.}}} \\
\endlastfoot
\multicolumn{10}{l}{\textit{Panel A: Setting 1 (linear), $a=0.8$}} \\*
\multirow{3}{*}{2000} & $\mu$ & 0.007 & 0.145 & 96.6 & 0.639 & 0.007 & 0.187 & 96.1 & 0.804 \\*
 & $f_A$ & -0.001 & 0.080 & 96.1 & 0.317 & -0.002 & 0.089 & 95.6 & 0.355 \\*
 & $S_C$ & -0.001 & 0.079 & 94.7 & 0.309 & -0.001 & 0.087 & 95.4 & 0.344 \\
\addlinespace[3pt]
\multirow{3}{*}{4000} & $\mu$ & 0.001 & 0.112 & 97.3 & 0.479 & 0.001 & 0.147 & 96.2 & 0.605 \\*
 & $f_A$ & 0.001 & 0.059 & 95.2 & 0.230 & 0.002 & 0.067 & 94.7 & 0.261 \\*
 & $S_C$ & 0.002 & 0.058 & 94.9 & 0.224 & 0.002 & 0.065 & 95.1 & 0.252 \\
\addlinespace[3pt]
\multirow{3}{*}{8000} & $\mu$ & -0.002 & 0.087 & 95.8 & 0.361 & -0.004 & 0.114 & 95.7 & 0.458 \\*
 & $f_A$ & 0.000 & 0.042 & 95.4 & 0.168 & 0.000 & 0.048 & 95.9 & 0.192 \\*
 & $S_C$ & 0.000 & 0.041 & 95.4 & 0.163 & 0.000 & 0.046 & 95.5 & 0.185 \\
\addlinespace[5pt]
\multicolumn{10}{l}{\textit{Panel B: Setting 1 (linear), $a=1.3$}} \\*
\multirow{3}{*}{2000} & $\mu$ & 0.011 & 0.226 & 95.8 & 0.984 & 0.008 & 0.295 & 95.8 & 1.222 \\*
 & $f_A$ & 0.001 & 0.093 & 95.9 & 0.377 & 0.001 & 0.106 & 95.9 & 0.430 \\*
 & $S_C$ & 0.001 & 0.090 & 94.9 & 0.357 & 0.000 & 0.100 & 95.5 & 0.402 \\
\addlinespace[3pt]
\multirow{3}{*}{4000} & $\mu$ & 0.007 & 0.164 & 97.6 & 0.734 & 0.013 & 0.210 & 97.0 & 0.916 \\*
 & $f_A$ & 0.003 & 0.070 & 93.8 & 0.274 & 0.003 & 0.080 & 94.2 & 0.315 \\*
 & $S_C$ & 0.002 & 0.067 & 93.7 & 0.259 & 0.002 & 0.075 & 94.3 & 0.295 \\
\addlinespace[3pt]
\multirow{3}{*}{8000} & $\mu$ & -0.002 & 0.133 & 95.9 & 0.550 & -0.001 & 0.172 & 95.4 & 0.690 \\*
 & $f_A$ & -0.003 & 0.051 & 95.7 & 0.201 & -0.003 & 0.059 & 95.7 & 0.233 \\*
 & $S_C$ & -0.003 & 0.049 & 95.1 & 0.189 & -0.003 & 0.056 & 95.1 & 0.217 \\
\addlinespace[5pt]
\multicolumn{10}{l}{\textit{Panel C: Setting 1 (linear), $a=1.8$}} \\*
\multirow{3}{*}{2000} & $\mu$ & 0.019 & 0.368 & 95.3 & 1.563 & 0.017 & 0.466 & 95.2 & 1.905 \\*
 & $f_A$ & 0.002 & 0.117 & 95.0 & 0.451 & 0.002 & 0.136 & 94.1 & 0.524 \\*
 & $S_C$ & 0.001 & 0.109 & 94.3 & 0.414 & 0.002 & 0.123 & 94.6 & 0.474 \\
\addlinespace[3pt]
\multirow{3}{*}{4000} & $\mu$ & -0.004 & 0.270 & 96.9 & 1.164 & -0.004 & 0.343 & 96.3 & 1.434 \\*
 & $f_A$ & 0.002 & 0.086 & 94.4 & 0.331 & 0.002 & 0.102 & 95.2 & 0.388 \\*
 & $S_C$ & 0.001 & 0.079 & 94.8 & 0.302 & 0.001 & 0.090 & 95.2 & 0.349 \\
\addlinespace[3pt]
\multirow{3}{*}{8000} & $\mu$ & 0.002 & 0.199 & 96.5 & 0.870 & 0.005 & 0.255 & 96.4 & 1.077 \\*
 & $f_A$ & -0.001 & 0.062 & 95.4 & 0.243 & 0.001 & 0.072 & 95.3 & 0.288 \\*
 & $S_C$ & 0.000 & 0.056 & 94.9 & 0.221 & 0.001 & 0.065 & 95.1 & 0.258 \\
\addlinespace[5pt]
\multicolumn{10}{l}{\textit{Panel D: Setting 2 (nonlinear), $a=0.8$}} \\*
\multirow{3}{*}{2000} & $\mu$ & 0.233 & 0.279 & 69.3 & 0.636 & 0.007 & 0.204 & 95.5 & 0.801 \\*
 & $f_A$ & 0.224 & 0.238 & 20.2 & 0.319 & -0.003 & 0.089 & 95.1 & 0.354 \\*
 & $S_C$ & 0.225 & 0.238 & 18.0 & 0.311 & -0.003 & 0.086 & 94.9 & 0.343 \\
\addlinespace[3pt]
\multirow{3}{*}{4000} & $\mu$ & 0.176 & 0.207 & 71.4 & 0.479 & 0.001 & 0.142 & 96.2 & 0.606 \\*
 & $f_A$ & 0.175 & 0.184 & 16.9 & 0.232 & 0.001 & 0.066 & 95.1 & 0.261 \\*
 & $S_C$ & 0.175 & 0.184 & 15.1 & 0.225 & 0.001 & 0.064 & 95.1 & 0.252 \\
\addlinespace[3pt]
\multirow{3}{*}{8000} & $\mu$ & 0.137 & 0.160 & 70.0 & 0.361 & 0.007 & 0.108 & 97.0 & 0.458 \\*
 & $f_A$ & 0.134 & 0.142 & 13.8 & 0.169 & 0.003 & 0.052 & 94.4 & 0.192 \\*
 & $S_C$ & 0.134 & 0.141 & 11.7 & 0.164 & 0.003 & 0.050 & 94.8 & 0.185 \\
\addlinespace[5pt]
\multicolumn{10}{l}{\textit{Panel E: Setting 2 (nonlinear), $a=1.3$}} \\*
\multirow{3}{*}{2000} & $\mu$ & 0.230 & 0.319 & 88.9 & 0.979 & -0.003 & 0.294 & 95.6 & 1.219 \\*
 & $f_A$ & 0.227 & 0.246 & 34.1 & 0.379 & -0.002 & 0.108 & 94.8 & 0.429 \\*
 & $S_C$ & 0.225 & 0.243 & 30.0 & 0.360 & -0.003 & 0.103 & 95.0 & 0.402 \\
\addlinespace[3pt]
\multirow{3}{*}{4000} & $\mu$ & 0.181 & 0.251 & 85.7 & 0.734 & 0.008 & 0.225 & 96.4 & 0.920 \\*
 & $f_A$ & 0.175 & 0.189 & 30.9 & 0.276 & 0.001 & 0.082 & 94.6 & 0.316 \\*
 & $S_C$ & 0.175 & 0.188 & 27.1 & 0.261 & 0.001 & 0.079 & 93.3 & 0.295 \\
\addlinespace[3pt]
\multirow{3}{*}{8000} & $\mu$ & 0.134 & 0.185 & 87.7 & 0.550 & 0.004 & 0.167 & 96.3 & 0.693 \\*
 & $f_A$ & 0.133 & 0.143 & 26.6 & 0.202 & 0.001 & 0.062 & 94.0 & 0.234 \\*
 & $S_C$ & 0.133 & 0.142 & 21.6 & 0.190 & 0.001 & 0.057 & 94.0 & 0.218 \\
\addlinespace[5pt]
\multicolumn{10}{l}{\textit{Panel F: Setting 2 (nonlinear), $a=1.8$}} \\*
\multirow{3}{*}{2000} & $\mu$ & 0.229 & 0.435 & 95.0 & 1.563 & 0.001 & 0.473 & 94.8 & 1.915 \\*
 & $f_A$ & 0.223 & 0.250 & 52.1 & 0.458 & -0.007 & 0.132 & 95.6 & 0.526 \\*
 & $S_C$ & 0.222 & 0.246 & 47.6 & 0.418 & -0.007 & 0.120 & 95.0 & 0.474 \\
\addlinespace[3pt]
\multirow{3}{*}{4000} & $\mu$ & 0.174 & 0.321 & 94.8 & 1.163 & -0.010 & 0.346 & 94.9 & 1.433 \\*
 & $f_A$ & 0.175 & 0.195 & 46.9 & 0.334 & 0.000 & 0.100 & 95.3 & 0.389 \\*
 & $S_C$ & 0.176 & 0.193 & 38.9 & 0.304 & 0.002 & 0.089 & 94.8 & 0.349 \\
\addlinespace[3pt]
\multirow{3}{*}{8000} & $\mu$ & 0.125 & 0.241 & 94.0 & 0.866 & -0.008 & 0.266 & 94.5 & 1.076 \\*
 & $f_A$ & 0.130 & 0.144 & 46.6 & 0.245 & -0.001 & 0.071 & 95.5 & 0.289 \\*
 & $S_C$ & 0.131 & 0.144 & 37.4 & 0.222 & -0.001 & 0.067 & 94.3 & 0.258 \\
\end{longtable}
\endgroup

Table~\ref{tab:sim-stability-undersmooth} gives the corresponding
undersmoothing results. Undersmoothing improves MR's coverage
in the nonlinear setting, while increasing interval lengths.
Under the regular bandwidth, MRDB retains small bias and
coverage generally close to the nominal level under each
of the three misspecification scenarios.

\begingroup
\renewcommand{\baselinestretch}{1}
\fontsize{9}{10.5}\selectfont
\setlength{\tabcolsep}{3pt}
\renewcommand{\arraystretch}{1.12}
\setlength{\LTleft}{\fill}
\setlength{\LTright}{\fill}
\setlength{\LTpre}{6pt}
\setlength{\LTpost}{4pt}
\begin{longtable}{@{}rl*{8}{r}@{}}
\caption{Stability experiment at representative doses with the undersmoothed bandwidth.}\label{tab:sim-stability-undersmooth}\\
\toprule
\multirow{2}{*}{$n$} & \multirow{2}{*}{Misspecified} & \multicolumn{4}{c}{MR} & \multicolumn{4}{c}{MRDB} \\
\cmidrule(lr){3-6}\cmidrule(lr){7-10}
 & & Bias & RMSE & CP (\%) & Length & Bias & RMSE & CP (\%) & Length \\
\midrule
\endfirsthead
\multicolumn{10}{c}{\tablename\ \thetable\ (continued)}\\
\toprule
\multirow{2}{*}{$n$} & \multirow{2}{*}{Misspecified} & \multicolumn{4}{c}{MR} & \multicolumn{4}{c}{MRDB} \\
\cmidrule(lr){3-6}\cmidrule(lr){7-10}
 & & Bias & RMSE & CP (\%) & Length & Bias & RMSE & CP (\%) & Length \\
\midrule
\endhead
\midrule
\multicolumn{10}{r@{}}{\footnotesize Continued on next page}\\
\endfoot
\bottomrule
\addlinespace[4pt]
\multicolumn{10}{@{}c@{}}{\makebox[0pt][c]{\parbox[t]{\textwidth}{\fontsize{8.5}{10}\selectfont \textit{Notes:} Each entry uses $R=1000$ repetitions. The Misspecified column identifies the only misspecified model; the other two models are correctly specified. Bias is the signed pointwise Monte Carlo bias, not RMSBias. RMSE, CP, and Length are pointwise root mean squared error, nominal $95\%$ interval coverage (\%), and mean interval length, respectively. The bandwidth is $h=n^{-1/4}$, with $b=1.25h$.}}} \\
\endlastfoot
\multicolumn{10}{l}{\textit{Panel A: Setting 1 (linear), $a=0.8$}} \\*
\multirow{3}{*}{2000} & $\mu$ & -0.009 & 0.249 & 96.7 & 1.053 & -0.015 & 0.327 & 96.1 & 1.344 \\*
 & $f_A$ & -0.003 & 0.106 & 94.7 & 0.421 & -0.004 & 0.125 & 94.9 & 0.503 \\*
 & $S_C$ & -0.002 & 0.102 & 95.0 & 0.404 & -0.003 & 0.120 & 95.5 & 0.480 \\
\addlinespace[3pt]
\multirow{3}{*}{4000} & $\mu$ & 0.001 & 0.206 & 95.3 & 0.809 & 0.001 & 0.266 & 94.4 & 1.033 \\*
 & $f_A$ & 0.002 & 0.082 & 94.4 & 0.316 & 0.001 & 0.099 & 94.0 & 0.380 \\*
 & $S_C$ & 0.002 & 0.079 & 94.9 & 0.303 & 0.002 & 0.095 & 94.4 & 0.363 \\
\addlinespace[3pt]
\multirow{3}{*}{8000} & $\mu$ & -0.004 & 0.156 & 95.0 & 0.625 & 0.001 & 0.202 & 94.4 & 0.800 \\*
 & $f_A$ & -0.001 & 0.060 & 96.1 & 0.237 & -0.001 & 0.074 & 94.6 & 0.288 \\*
 & $S_C$ & -0.001 & 0.057 & 94.9 & 0.227 & -0.001 & 0.070 & 95.1 & 0.274 \\
\addlinespace[5pt]
\multicolumn{10}{l}{\textit{Panel B: Setting 1 (linear), $a=1.3$}} \\*
\multirow{3}{*}{2000} & $\mu$ & 0.005 & 0.398 & 94.4 & 1.572 & 0.014 & 0.514 & 93.6 & 1.988 \\*
 & $f_A$ & 0.005 & 0.133 & 94.7 & 0.517 & 0.006 & 0.163 & 94.4 & 0.625 \\*
 & $S_C$ & 0.003 & 0.123 & 94.6 & 0.480 & 0.003 & 0.148 & 94.2 & 0.576 \\
\addlinespace[3pt]
\multirow{3}{*}{4000} & $\mu$ & 0.005 & 0.303 & 94.4 & 1.205 & -0.002 & 0.394 & 93.2 & 1.529 \\*
 & $f_A$ & 0.004 & 0.098 & 95.2 & 0.388 & 0.004 & 0.120 & 94.8 & 0.473 \\*
 & $S_C$ & 0.004 & 0.090 & 95.1 & 0.360 & 0.004 & 0.109 & 94.8 & 0.435 \\
\addlinespace[3pt]
\multirow{3}{*}{8000} & $\mu$ & -0.009 & 0.235 & 95.1 & 0.931 & -0.014 & 0.309 & 94.1 & 1.184 \\*
 & $f_A$ & -0.003 & 0.074 & 94.5 & 0.294 & -0.004 & 0.093 & 94.8 & 0.361 \\*
 & $S_C$ & -0.002 & 0.069 & 95.0 & 0.271 & -0.003 & 0.086 & 94.8 & 0.331 \\
\addlinespace[5pt]
\multicolumn{10}{l}{\textit{Panel C: Setting 1 (linear), $a=1.8$}} \\*
\multirow{3}{*}{2000} & $\mu$ & -0.032 & 0.624 & 90.6 & 2.321 & -0.041 & 0.793 & 89.8 & 2.886 \\*
 & $f_A$ & 0.006 & 0.169 & 92.8 & 0.628 & 0.007 & 0.207 & 93.1 & 0.765 \\*
 & $S_C$ & 0.005 & 0.151 & 93.8 & 0.574 & 0.005 & 0.181 & 93.9 & 0.697 \\
\addlinespace[3pt]
\multirow{3}{*}{4000} & $\mu$ & -0.019 & 0.462 & 92.9 & 1.835 & -0.021 & 0.596 & 92.6 & 2.305 \\*
 & $f_A$ & 0.000 & 0.127 & 94.0 & 0.482 & 0.000 & 0.154 & 94.3 & 0.593 \\*
 & $S_C$ & 0.000 & 0.112 & 94.1 & 0.432 & -0.002 & 0.134 & 95.0 & 0.529 \\
\addlinespace[3pt]
\multirow{3}{*}{8000} & $\mu$ & -0.010 & 0.352 & 94.4 & 1.423 & -0.006 & 0.448 & 94.2 & 1.799 \\*
 & $f_A$ & 0.002 & 0.091 & 95.4 & 0.368 & 0.002 & 0.114 & 94.4 & 0.456 \\*
 & $S_C$ & 0.003 & 0.082 & 95.3 & 0.327 & 0.003 & 0.102 & 94.9 & 0.403 \\
\addlinespace[5pt]
\multicolumn{10}{l}{\textit{Panel D: Setting 2 (nonlinear), $a=0.8$}} \\*
\multirow{3}{*}{2000} & $\mu$ & 0.030 & 0.278 & 94.3 & 1.045 & 0.000 & 0.357 & 93.7 & 1.333 \\*
 & $f_A$ & 0.023 & 0.108 & 94.3 & 0.419 & -0.005 & 0.127 & 94.9 & 0.501 \\*
 & $S_C$ & 0.025 & 0.104 & 94.6 & 0.403 & -0.002 & 0.121 & 95.5 & 0.479 \\
\addlinespace[3pt]
\multirow{3}{*}{4000} & $\mu$ & 0.007 & 0.200 & 95.0 & 0.811 & -0.018 & 0.261 & 94.8 & 1.036 \\*
 & $f_A$ & 0.019 & 0.082 & 94.6 & 0.316 & -0.002 & 0.096 & 95.4 & 0.381 \\*
 & $S_C$ & 0.019 & 0.079 & 94.5 & 0.303 & -0.001 & 0.091 & 95.0 & 0.363 \\
\addlinespace[3pt]
\multirow{3}{*}{8000} & $\mu$ & 0.021 & 0.160 & 94.9 & 0.626 & 0.010 & 0.207 & 94.6 & 0.800 \\*
 & $f_A$ & 0.015 & 0.067 & 92.7 & 0.238 & 0.000 & 0.079 & 93.2 & 0.289 \\*
 & $S_C$ & 0.015 & 0.064 & 92.5 & 0.227 & 0.000 & 0.075 & 92.4 & 0.274 \\
\addlinespace[5pt]
\multicolumn{10}{l}{\textit{Panel E: Setting 2 (nonlinear), $a=1.3$}} \\*
\multirow{3}{*}{2000} & $\mu$ & 0.005 & 0.401 & 94.2 & 1.568 & -0.021 & 0.508 & 93.8 & 1.986 \\*
 & $f_A$ & 0.025 & 0.133 & 95.2 & 0.516 & -0.002 & 0.159 & 94.3 & 0.623 \\*
 & $S_C$ & 0.025 & 0.124 & 95.1 & 0.481 & -0.001 & 0.145 & 95.1 & 0.577 \\
\addlinespace[3pt]
\multirow{3}{*}{4000} & $\mu$ & 0.018 & 0.314 & 94.2 & 1.214 & 0.003 & 0.404 & 93.1 & 1.543 \\*
 & $f_A$ & 0.019 & 0.104 & 93.4 & 0.388 & 0.000 & 0.125 & 93.9 & 0.474 \\*
 & $S_C$ & 0.019 & 0.098 & 92.6 & 0.360 & -0.001 & 0.116 & 93.9 & 0.436 \\
\addlinespace[3pt]
\multirow{3}{*}{8000} & $\mu$ & 0.016 & 0.226 & 96.6 & 0.936 & -0.001 & 0.291 & 96.3 & 1.191 \\*
 & $f_A$ & 0.016 & 0.077 & 94.5 & 0.295 & 0.003 & 0.093 & 94.2 & 0.362 \\*
 & $S_C$ & 0.016 & 0.070 & 95.1 & 0.272 & 0.003 & 0.084 & 94.9 & 0.331 \\
\addlinespace[5pt]
\multicolumn{10}{l}{\textit{Panel F: Setting 2 (nonlinear), $a=1.8$}} \\*
\multirow{3}{*}{2000} & $\mu$ & -0.044 & 0.627 & 91.3 & 2.345 & -0.088 & 0.800 & 89.1 & 2.932 \\*
 & $f_A$ & 0.018 & 0.167 & 94.7 & 0.632 & -0.008 & 0.209 & 94.1 & 0.771 \\*
 & $S_C$ & 0.017 & 0.147 & 93.9 & 0.574 & -0.011 & 0.180 & 94.8 & 0.698 \\
\addlinespace[3pt]
\multirow{3}{*}{4000} & $\mu$ & -0.044 & 0.476 & 92.1 & 1.824 & -0.083 & 0.624 & 91.7 & 2.291 \\*
 & $f_A$ & 0.019 & 0.126 & 94.2 & 0.481 & 0.000 & 0.155 & 94.4 & 0.591 \\*
 & $S_C$ & 0.020 & 0.110 & 94.7 & 0.433 & 0.000 & 0.132 & 95.5 & 0.530 \\
\addlinespace[3pt]
\multirow{3}{*}{8000} & $\mu$ & -0.017 & 0.361 & 94.3 & 1.419 & -0.031 & 0.467 & 91.9 & 1.792 \\*
 & $f_A$ & 0.011 & 0.094 & 94.8 & 0.369 & -0.003 & 0.115 & 94.9 & 0.457 \\*
 & $S_C$ & 0.012 & 0.084 & 94.5 & 0.327 & -0.002 & 0.101 & 95.0 & 0.402 \\
\end{longtable}
\endgroup

\begin{figure}[!htb]
\centering
\linespread{1}\selectfont
\includegraphics[width=\textwidth,height=0.70\textheight,keepaspectratio]{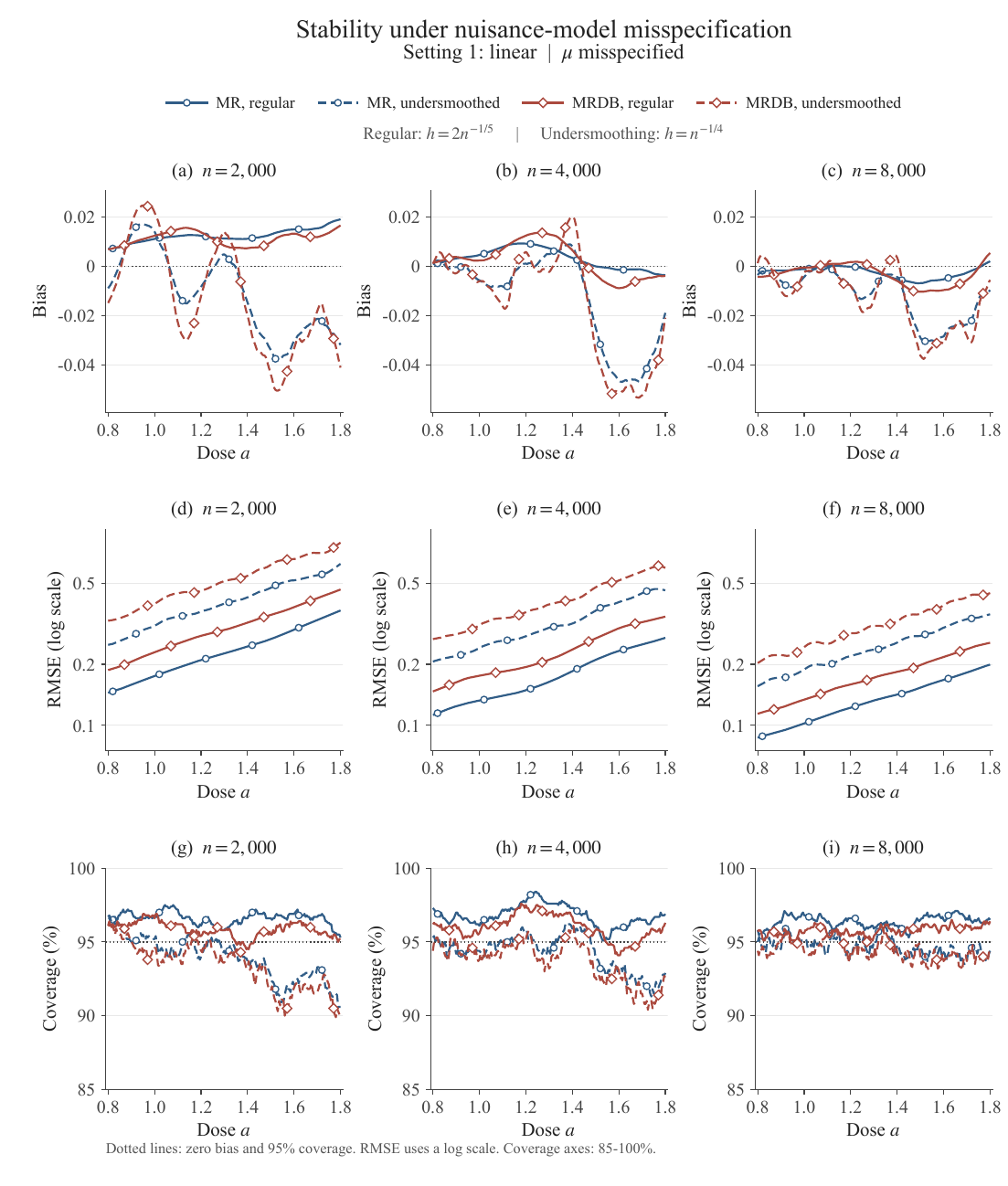}
\caption{Setting 1 (linear): only the outcome model $\mu$ is misspecified.}
\label{fig:sim-stability-linear-mu}
\end{figure}

\begin{figure}[!htb]
\centering
\linespread{1}\selectfont
\includegraphics[width=\textwidth,height=0.70\textheight,keepaspectratio]{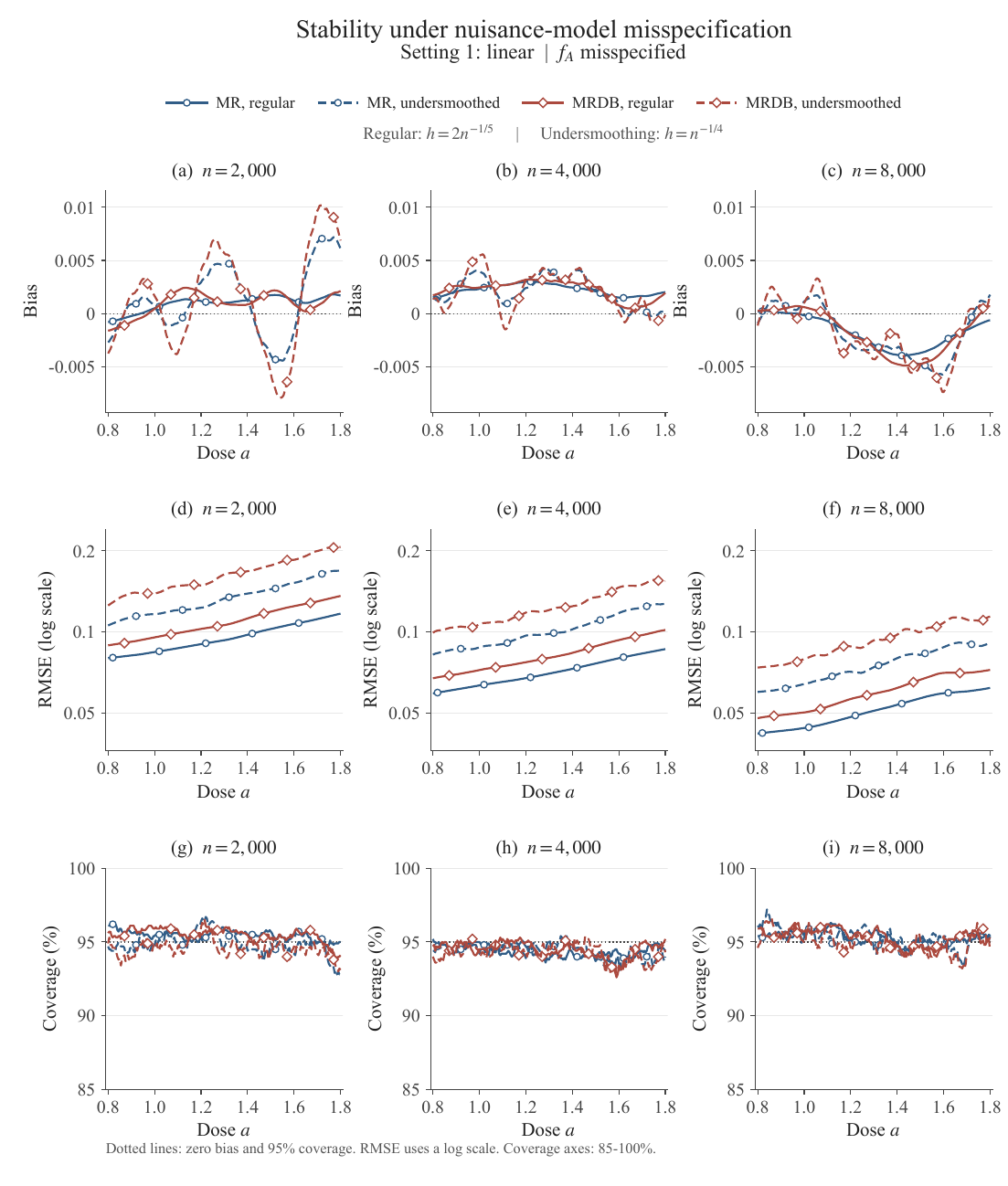}
\caption[Setting 1 (linear): only the treatment density is misspecified.]
{Setting 1 (linear): only the treatment density $f_A$ is misspecified.}
\label{fig:sim-stability-linear-a}
\end{figure}

\begin{figure}[!htb]
\centering
\linespread{1}\selectfont
\includegraphics[width=\textwidth,height=0.70\textheight,keepaspectratio]{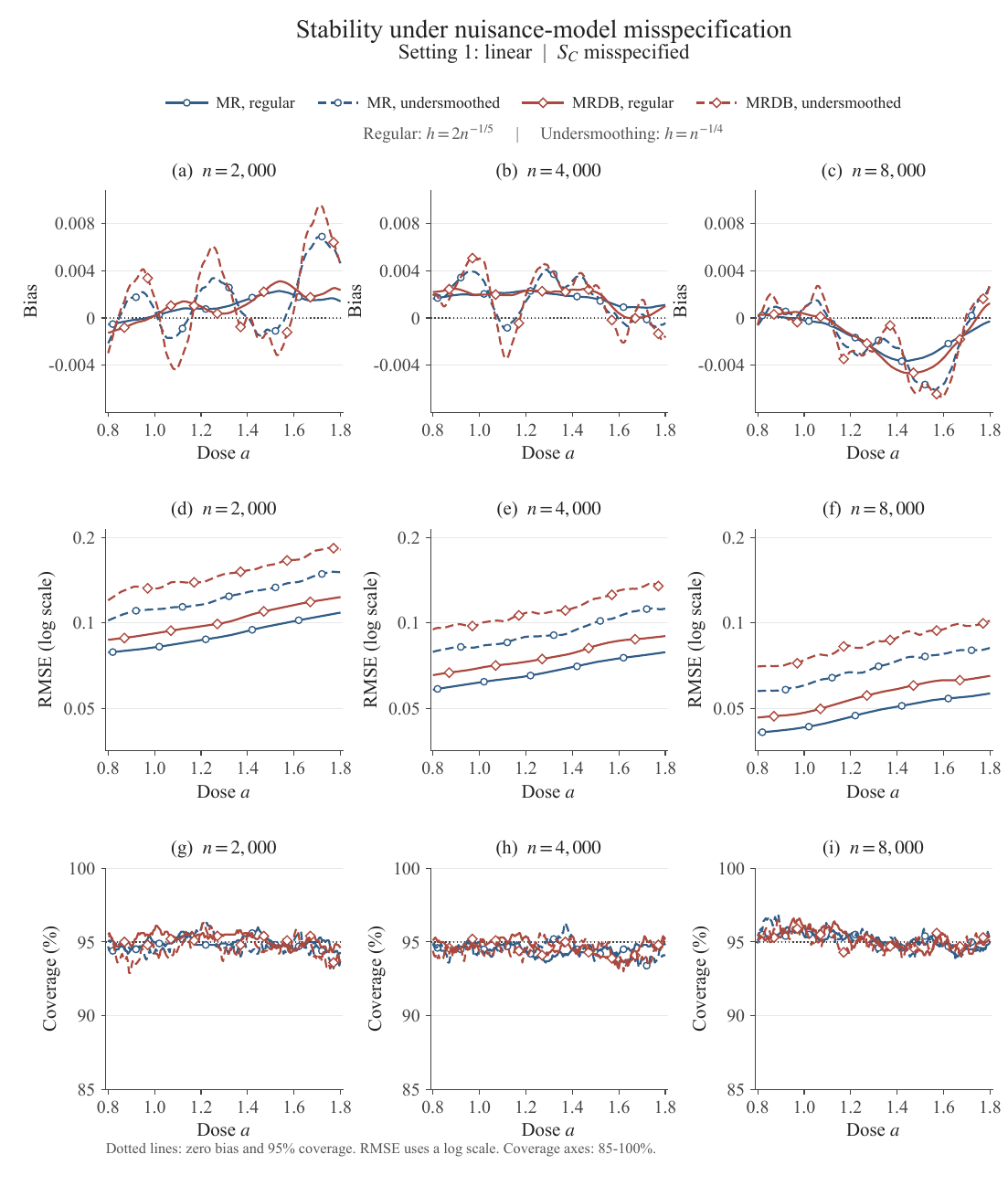}
\caption[Setting 1 (linear): only the censoring model is misspecified.]
{Setting 1 (linear): only the censoring model $S_C$ is misspecified.}
\label{fig:sim-stability-linear-c}
\end{figure}

\begin{figure}[!htb]
\centering
\linespread{1}\selectfont
\includegraphics[width=\textwidth,height=0.70\textheight,keepaspectratio]{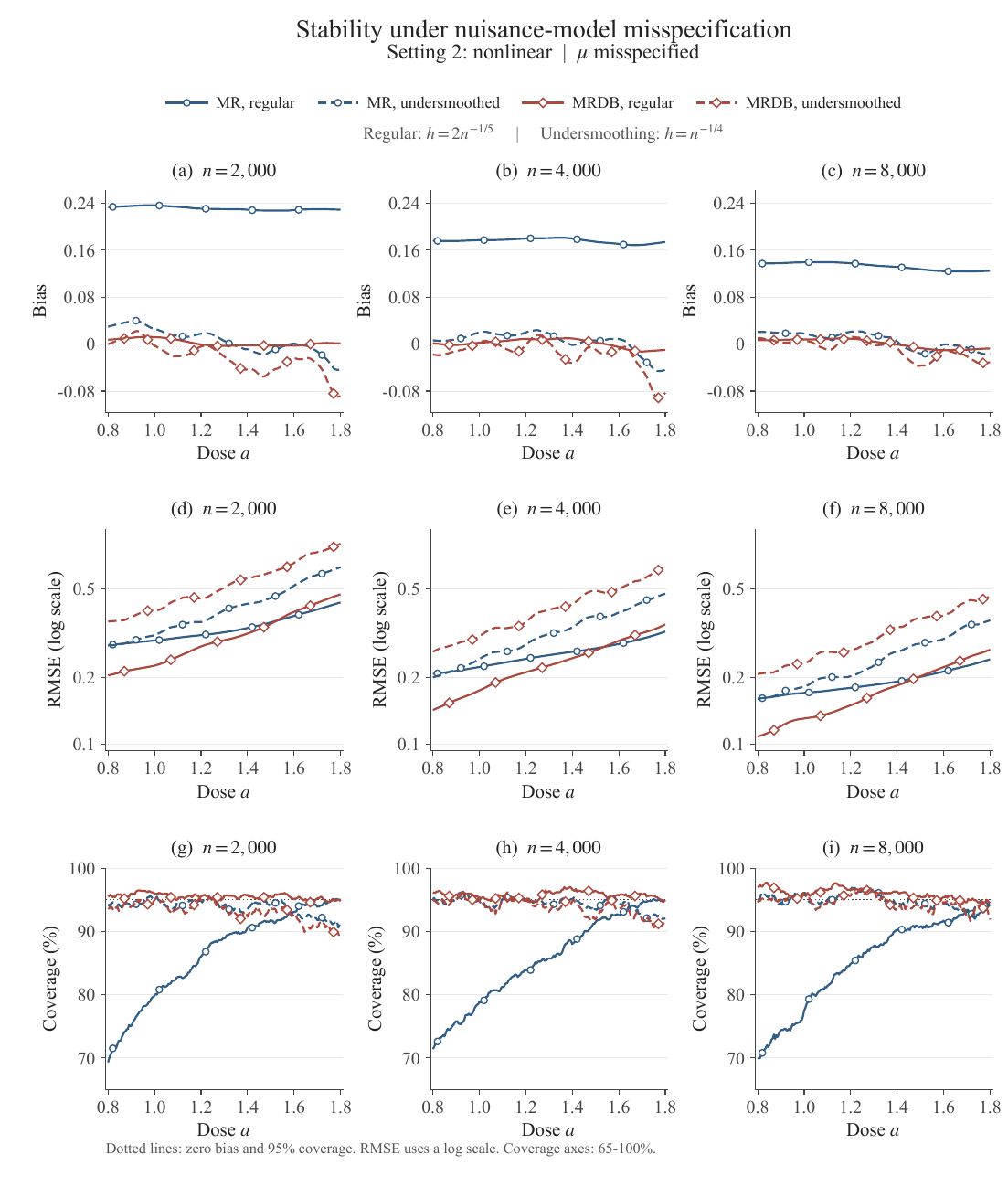}
\caption{Setting 2 (nonlinear): only the outcome model $\mu$ is misspecified.}
\label{fig:sim-stability-nonlinear-mu}
\end{figure}

\begin{figure}[!htb]
\centering
\linespread{1}\selectfont
\includegraphics[width=\textwidth,height=0.70\textheight,keepaspectratio]{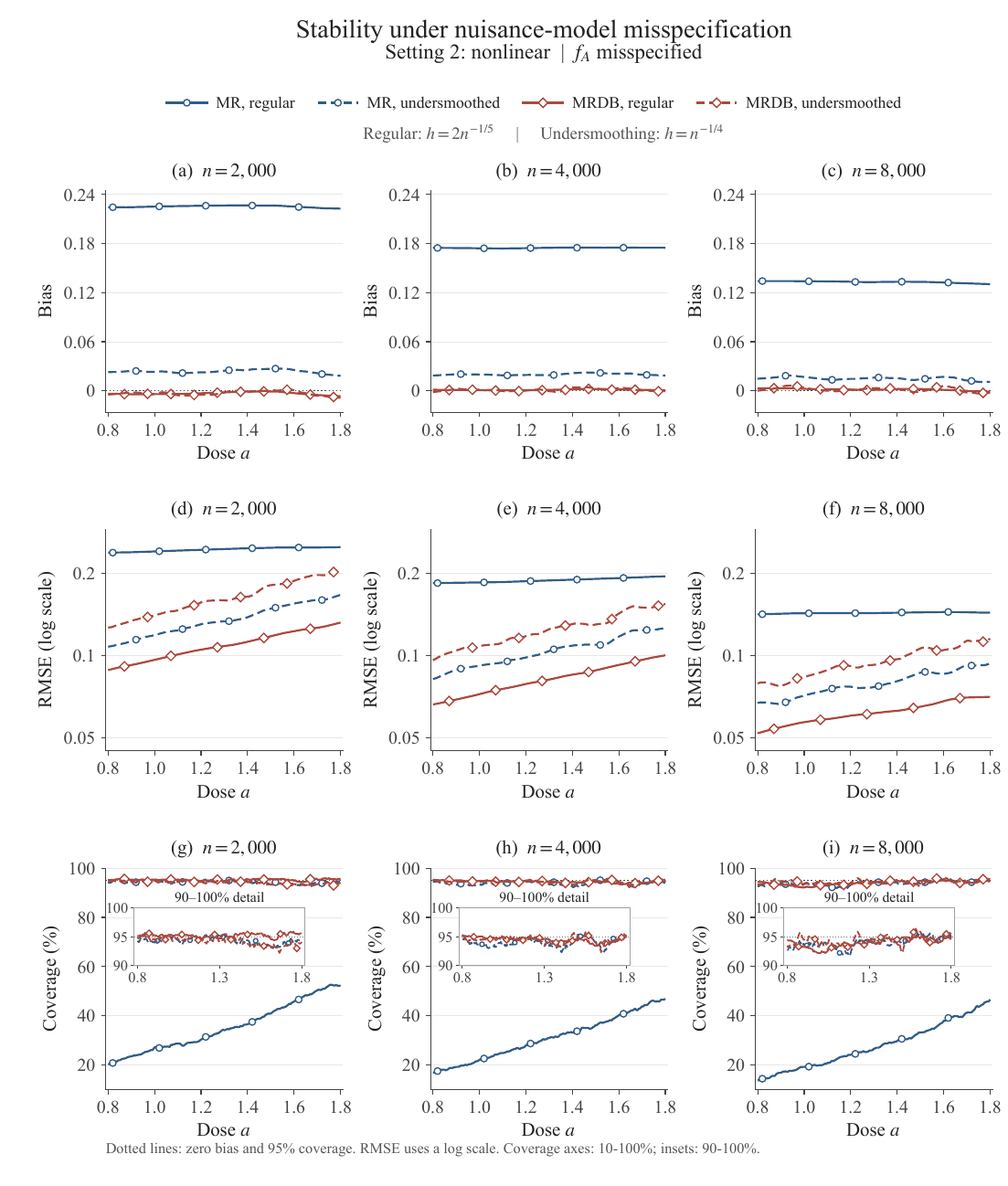}
\caption[Setting 2 (nonlinear): only the treatment density is misspecified.]
{Setting 2 (nonlinear): only the treatment density $f_A$ is misspecified.}
\label{fig:sim-stability-nonlinear-a}
\end{figure}

\begin{figure}[!htb]
\centering
\linespread{1}\selectfont
\includegraphics[width=\textwidth,height=0.70\textheight,keepaspectratio]{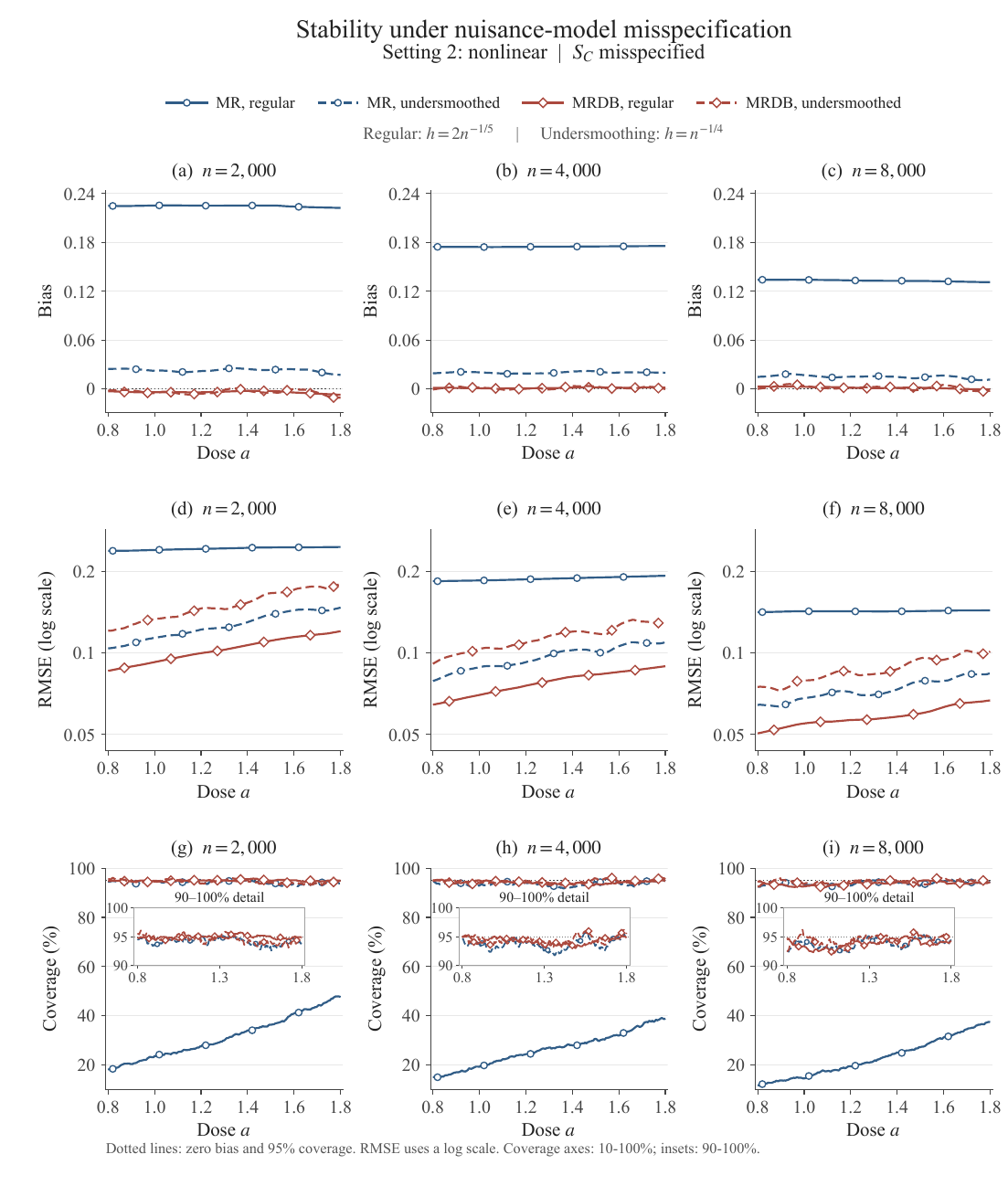}
\caption[Setting 2 (nonlinear): only the censoring model is misspecified.]
{Setting 2 (nonlinear): only the censoring model $S_C$ is misspecified.}
\label{fig:sim-stability-nonlinear-c}
\end{figure}

Figures~\ref{fig:sim-stability-linear-mu}-\ref{fig:sim-stability-nonlinear-c}
show the full-dose results. Rows display bias, RMSE, and
pointwise coverage; columns correspond to $n=2000,4000,8000$.
Solid and dashed curves denote regular bandwidths and
undersmoothing, respectively. RMSE axes are logarithmic.
In Setting 1, MRDB remains approximately unbiased across
the three misspecification scenarios.

In Setting 2, MRDB corrects the smoothing bias that remains in
regular-bandwidth MR despite two correctly specified nuisance
models. Its full-dose coverage is generally close to the nominal
level under the regular bandwidth. With a misspecified outcome
model, undersmoothed MRDB exhibits some finite-sample
undercoverage near the upper end of the dose interval.

\clearpage

\clearpage


\section{Efficient Influence Function with an Always Observed Outcome}
\label{app:always_observed_outcome}

Consider the observation scheme $O^+=(X,U,\Delta,Y)$, in which $Y$ is observed even when $A$ is censored. We retain Assumptions~\ref{ass: consistency}-\ref{ass: positivity} and the support $\mathcal A=[0,T]$. Apart from these conditions, the distribution of $Z$ and the conditional distribution of $C$ given $X$ are unrestricted within a dominated model. The integrated functional remains $\Phi=\int_{\mathcal A}m(a)p_A(a)\,da$, with full-data efficient influence function $\varphi^F$ given in \eqref{eq:full_eif}.

The censoring assumption implies $C\perp\!\!\!\perp A\mid X,Y$ and $P(C\ge t\mid X,Y)=S_C(t\mid X)$. Thus, $M_C$ remains a censoring martingale when the initial information includes $(X,Y)$. In addition to $V(t,X)=E\{\varphi^F(Z)\mid A\ge t,X\}$, define
\[
V_Y(t,X,Y)=E\{\varphi^F(Z)\mid A\ge t,X,Y\}.
\]
Set $V_Y$ to zero wherever $P(A\ge t\mid X,Y)=0$, and use predictable versions in martingale integrals. Conditioning on $Y$ generally changes this tail expectation, since the outcome residual in $\varphi^F$ no longer has conditional mean zero.

\begin{theorem}[Efficient Influence Function with an Always Observed Outcome]
\label{thm:eif_always_observed}
Under the model specified above, suppose that $E[\{\varphi^F(Z)\}^2]<\infty$. The efficient influence function of $\Phi$ based on $O^+$ is
\begin{equation}
\label{eq:eif_always_observed}
\varphi^{O,+}(O^+)=\frac{\Delta\varphi^F(X,U,Y)}{S_C(U\mid X)}+\int_0^T\frac{V_Y(t,X,Y)}{S_C(t\mid X)}\,dM_C(t),
\end{equation}
where the first term is zero when $\Delta=0$. Let $\varphi^O(O)$ denote the influence function in \eqref{eq:observed_eif}, evaluated using $O=(X,U,\Delta,\Delta Y)$. Then
\begin{equation}
\label{eq:eif_always_observed_gain}
\begin{split}
&\operatorname{Var}\{\varphi^O(O)\}-\operatorname{Var}\{\varphi^{O,+}(O^+)\}\\
&\quad=E\!\left[\int_0^A\{V_Y(t,X,Y)-V(t,X)\}^2\frac{\lambda_C(t\mid X)}{S_C(t\mid X)}\,dt\right]\ge0.
\end{split}
\end{equation}
The variance reduction is strict when the expectation on the right-hand side is positive.
\end{theorem}

\begin{proof}[Proof of Theorem~\ref{thm:eif_always_observed}]
We verify the pathwise derivative and membership in the observed-data tangent space, following the standard characterization of efficient influence functions \citep{bickel1993,tsiatis2006}. To avoid repeated tail-integral calculations, define
\[
\widetilde S_Y(Z)=\frac{\varphi^F(Z)}{S_C(A\mid X)}-\int_0^A\frac{V_Y(t,X,Y)\lambda_C(t\mid X)}{S_C(t\mid X)}\,dt.
\]
The inverse-survival identity gives $\int_0^A\lambda_C(t\mid X)/S_C(t\mid X)\,dt=S_C(A\mid X)^{-1}-1\le c_2^{-1}-1$. Conditional Jensen's inequality and Fubini's theorem therefore yield
\[
E\!\left[\int_0^A V_Y(t,X,Y)^2\frac{\lambda_C(t\mid X)}{S_C(t\mid X)}\,dt\right]
\le(c_2^{-1}-1)E[\{\varphi^F(Z)\}^2].
\]
By Cauchy-Schwarz, the second moment of the integral in $\widetilde S_Y$ is at most $(c_2^{-1}-1)^2E[\{\varphi^F(Z)\}^2]$. Positivity then implies $\widetilde S_Y\in L_2(P_Z)$. The same bound and the martingale isometry ensure that the integral in \eqref{eq:eif_always_observed} is square-integrable.

We first show that the candidate in \eqref{eq:eif_always_observed} satisfies
\begin{equation}
\label{eq:eif_always_observed_preimage}
E\{\widetilde S_Y(Z)\mid O^+\}=\varphi^{O,+}(O^+).
\end{equation}
When $\Delta=1$, this follows directly by expanding $dM_C(t)$, since $A=U$ and $Z$ is observed. When $\Delta=0$ and $U=u$, the conditional law of $A$ is its law given $(X,Y,A>u)$. On this event, the inverse-survival identity gives
\[
\begin{aligned}
\widetilde S_Y(Z)
={}&\frac{\varphi^F(Z)}{S_C(u\mid X)}-\int_0^u\frac{V_Y(t,X,Y)\lambda_C(t\mid X)}{S_C(t\mid X)}\,dt\\
&+\int_u^A\frac{\{\varphi^F(Z)-V_Y(t,X,Y)\}\lambda_C(t\mid X)}{S_C(t\mid X)}\,dt.
\end{aligned}
\]
The last integral has conditional mean zero given $(X,Y,A>u)$, because $E[I(A\ge t)\{\varphi^F(Z)-V_Y(t,X,Y)\}\mid X,Y]=0$ for $t\ge u$. Hence
\[
E\{\widetilde S_Y(Z)\mid X,Y,A>u\}
=\frac{V_Y(u,X,Y)}{S_C(u\mid X)}-\int_0^u\frac{V_Y(t,X,Y)\lambda_C(t\mid X)}{S_C(t\mid X)}\,dt,
\]
which is the candidate on censored observations. The conventions $A>t$ and $A\ge t$ agree almost surely at censoring times and almost everywhere in these integrals because $C\mid(X,Y)$ has a density. Thus, \eqref{eq:eif_always_observed_preimage} holds without requiring $A\mid(X,Y)$ to have a density.

Next, conditional on $Z$, the inverse-probability term in \eqref{eq:eif_always_observed} has mean $\varphi^F(Z)$ and the martingale term has mean zero. Consequently,
\[
E\{\varphi^{O,+}(O^+)\mid Z\}=\varphi^F(Z),\qquad E\{\varphi^{O,+}(O^+)\}=0.
\]
To handle variation in the censoring distribution, expand the martingale and condition on $(X,Y,C=c)$. The definition of $V_Y$ gives
\[
\begin{aligned}
&E\{\varphi^{O,+}(O^+)\mid X,Y,C=c\}\\
&\quad=E\!\left[\left.\varphi^F(Z)\left\{\frac{1}{S_C(A\wedge c\mid X)}-\int_0^{A\wedge c}\frac{\lambda_C(t\mid X)}{S_C(t\mid X)}\,dt\right\}\right|X,Y\right]\\
&\quad=E\{\varphi^F(Z)\mid X,Y\}.
\end{aligned}
\]
Since $C\perp\!\!\!\perp Y\mid X$, averaging over $Y$ yields $E\{\varphi^{O,+}(O^+)\mid X,C\}=E\{\varphi^F(Z)\mid X\}$.

For a regular parametric submodel, let $S_Z(Z)$ and $S_{\mathrm{cen}}(C,X)$ denote the full-data and conditional censoring scores, with $E(S_Z)=0$ and $E(S_{\mathrm{cen}}\mid X)=0$. The observed-data score is $S_+(O^+)=E\{S_Z(Z)+S_{\mathrm{cen}}(C,X)\mid O^+\}$. The preceding conditional identities imply
\[
E\{\varphi^{O,+}(O^+)S_+(O^+)\}=E\{\varphi^F(Z)S_Z(Z)\}
=\left.\frac{\partial}{\partial\epsilon}\Phi(P_{Z,\epsilon})\right|_{\epsilon=0}.
\]
Thus, the candidate represents the pathwise derivative for every regular submodel.

It remains to verify that it belongs to the closed observed-data tangent space. By \eqref{eq:eif_always_observed_preimage} and the mean-zero property above, $E(\widetilde S_Y)=0$. Since the full-data model is nonparametric, $\widetilde S_Y\in L_2^0(P_Z)$ can be approximated by bounded mean-zero full-data scores. Holding the conditional censoring distribution fixed, their observed-data scores are the corresponding conditional expectations given $O^+$. Conditional-expectation contraction and \eqref{eq:eif_always_observed_preimage} show that these scores converge in $L_2(P_{O^+})$ to $\varphi^{O,+}$. This proves tangent-space membership and hence efficiency.

Finally, $\varphi^O(O)$ remains an influence function in the present model, since $O$ is a deterministic function of $O^+$. Its difference from the efficient influence function is
\[
\varphi^O(O)-\varphi^{O,+}(O^+)=\int_0^T\frac{V(t,X)-V_Y(t,X,Y)}{S_C(t\mid X)}\,dM_C(t).
\]
This difference is orthogonal to the observed-data tangent space and therefore to $\varphi^{O,+}$. The variance difference is consequently the second moment of this martingale integral. Applying the martingale isometry gives \eqref{eq:eif_always_observed_gain}.
\end{proof}

 \putbib[bibfile_supp] 
\end{bibunit}	
    
\end{document}